\documentclass[oneside,english]{amsart}
\usepackage[T1]{fontenc}
\usepackage[latin9]{inputenc}
\usepackage{geometry}
\usepackage{amstext}
\usepackage{amsthm}
\usepackage{amssymb}
\usepackage{wasysym}

\makeatletter
\numberwithin{equation}{section}
\numberwithin{figure}{section}
\theoremstyle{plain}
\newtheorem{thm}{\protect\theoremname}[section]
\theoremstyle{remark}
\newtheorem{rem}[thm]{\protect\remarkname}
\theoremstyle{definition}
\newtheorem{defn}[thm]{\protect\definitionname}
\theoremstyle{plain}
\newtheorem{prop}[thm]{\protect\propositionname}
\theoremstyle{plain}
\newtheorem{lem}[thm]{\protect\lemmaname}
\theoremstyle{plain}
\newtheorem{cor}[thm]{\protect\corollaryname}

\makeatother

\usepackage{babel}
\providecommand{\corollaryname}{Corollary}
\providecommand{\definitionname}{Definition}
\providecommand{\lemmaname}{Lemma}
\providecommand{\propositionname}{Proposition}
\providecommand{\remarkname}{Remark}
\providecommand{\theoremname}{Theorem}

\begin{document}
\title{Kähler-Ricci Tangent Flows are Infinitesimally Algebraic}
\author{Max Hallgren}
\begin{abstract}
We show that any tangent cone of a singular shrinking Kähler-Ricci
soliton is a normal affine algebraic variety. Moreover, the regular
set of such a tangent cone in the metric sense coincides with the
regular set in the algebraic sense. Along the way, we give a parabolic
proof of Hörmander's $L^{2}$ estimate, which can be used to solve
the $\overline{\partial}$-equation on any singular shrinking Kähler-Ricci
soliton.
\end{abstract}

\maketitle

\section{Introduction}

Recent years have seen fruitful interactions between the structure
theory of singular geometric limits in Riemannian geometry and $L^{2}$
methods in complex geometry. Foundational to this interaction is the
description of noncollapsed Riemannian manifolds with a lower bound
on Ricci curvature and their Gromov-Hausdorff limits, known as Ricci
limit spaces. This was developed in \cite{cheegercolding1,cheegernaberdim4,cheegernaberquant},
where it was shown that noncollapsed limits of Riemannian manifolds
with bounded Ricci curvature are infinitesimally conical, with a singular
set of Minkowski codimension four. In \cite{donaldsun1}, this structure
theory was combined with Hörmander's $L^{2}$ method to establish
the partial $C^{0}$ estimate for Kähler-Einstein manifolds. Adapting
this to the to the more general setting of conical Kähler-Einstein
metrics in \cite{CDSI,CDSII,CDSIII,TianKStable} enabled a resolution
of the Yau-Tian-Donaldson conjecture. Similar techniques have also
been used to show uniqueness of metric tangent cones for polarized
limits of Kähler-Einstein metrics \cite{donaldsun2}, and have enabled
an algebraic description \cite{liuniquecones} of these cones related
to local notions of K-Stability. In addition, these methods have been
pushed to the setting of noncollapsed Kähler Ricci limit spaces in
\cite{szek1,szek2}, where it was shown that any tangent cone of such
a limit is an affine algebraic variety. These results are essential
for connecting the metric geometry of these limits with notions of
stability in algebraic geometry.

The interplay between metric geometry and $L^{2}$ methods has been
extended to the Kähler-Ricci flow in certain settings. In particular,
\cite{tianzhen,chenwang2a,chenwang2b} established the partial $C^{0}$
estimate for Fano Kähler-Ricci flows, which was used in \cite{chensunwang}
to give a Ricci flow proof of the Yau-Tian-Donaldson conjecture. The
key input from Riemannian geometry is the structure theory for limits
of Ricci flows with bounded scalar curvature developed in \cite{bam1,bam2,chenwang2a,chenwang2b},
which mirrors the theory of noncollapsed limits of Riemannian manifolds
with bounded Ricci curvature. Building on this, it was shown in \cite{szekdegen,lihanunique}
that the tangent flow is unique and algebraically determined, which
has led to examples of new Kähler-Ricci flow behavior \cite{tianGman}.

In \cite{bamlergen1,bamlergen2,bamlergen3}, a parabolic analogue
of the theory of Ricci limit spaces was developed. In particular,
a notion of convergence for flows of metric measure spaces was introduced
in \cite{bamlergen2}, where it was shown that noncollapsed Ricci
flows form a precompact collection with respect to this convergence.
In \cite{bamlergen3}, a partial regularity theory for the limits
of noncollapsed Ricci flows was developed. It was shown that such
limits are infinitesimally self-similar, and have singularities of
codimension four. This structure theory was combined with geometric
estimates for projective Kähler-Ricci flows in \cite{JSTI} to prove
a partial $C^{0}$ estimate near certain points called Ricci vertices.
Using this, it was shown in \cite{JSTI,JSTII} that Gromov-Hausdorff
limits of Kähler-Ricci flow based at Ricci vertices are infinitesimally
conical, continuous in time, and that their time slices are normal
analytic varieties. Away from Ricci vertices (whose locations are
so far difficult to determine in general), or without the assumption
of projectivity, there is little known for general limits of Kähler-Ricci
flow that is not present in the Riemannian setting. In \cite{complexsplitting},
the Kähler structure was used to improve the description of the singular
sets and tangent flows of singularity models, analogous to results
for Ricci limit spaces proved in \cite{cheegercoldingtian},\cite{gangliucone}.
However, the conjectural picture -- at least in the projective case
(see \cite{JSTI}) -- is much stronger: any Kähler-Ricci tangent
flow should be a normal analytic variety. 

In light of this, our focus in this paper is to improve the structure
theory of Kähler-Ricci tangent flows. We will prove a Hörmander-type
estimate on any Kähler-Ricci tangent flow, and apply this estimate
to show that such a flow is infinitesimally a normal affine algebraic
variety.

More concretely, suppose $(M,(\widetilde{g}_{t})_{t\in[0,T)},J)$
is a compact solution to the Kähler-Ricci flow

\[
\partial_{t}\widetilde{g}_{t}=-Rc(\widetilde{g}_{t})
\]
which develops a singularity at time $T<\infty$. Fix a conjugate
heat kernel $(\nu_{p,T;t})_{t\in[0,T)}$ based at the singular time
(for details, see Section 2.6 of \cite{bamlergen3}), and a sequence
$t_{i}\nearrow T$. By Bamler's compactness theory \cite{bamlergen2,bamlergen3},
we can pass to a subsequence so that the rescaled flows $g_{i,t}:=(T-t_{i})^{-1}\widetilde{g}_{T+(T-t_{i})t}$
with reference conjugate heat flows $\nu_{t}^{i}:=\nu_{p,T;T+(T-t_{i})t}$
$\mathbb{F}$-converge within some correspondence (see Section 5.1
of \cite{bamlergen2} for definitions) to a metric soliton $(\mathcal{X},(\nu_{t})_{t\in(-\infty,0)})$,
uniformly on compact time intervals:
\[
(M_{i},(g_{i,t})_{t\in[(T-t_{i})^{-1}T,0)},(\nu_{t}^{i})_{t\in[(T-t_{i})^{-1}T,0)})\xrightarrow[i\to\infty]{\mathbb{F}}(\mathcal{X},(\nu_{t})_{t\in(-\infty,0)}).
\]
Moreover, $\mathcal{X}$ is the metric flow corresponding to a gradient
Kähler-Ricci soliton with singularities of codimension four, and $d\nu_{t}=(2\pi\tau)^{-n}e^{-f}dg_{t}$
on the regular part $(\mathcal{R},g)$ of $\mathcal{X}$, where $f$
is the soliton potential function:
\[
Rc(g)+\nabla\overline{\nabla}f=\frac{1}{\tau}g,\,\,\,\,\,\,\,\,\,\,\,\,\,\,\nabla\nabla f=0.
\]
Given any $t_{0}\in(-\infty,0)$, $x_{0}\in\mathcal{X}_{t_{0}}$,
and $\lambda_{i}\nearrow\infty$, the sequence $(\mathcal{X}_{t_{0}},\lambda_{i}d_{t_{0}},x_{0})$
subconverges in the pointed Gromov-Hausdroff sense to a metric cone
$(C(Y),d_{C(Y)})$. Let $\mathcal{R}_{C(Y)}$ denote the regular set
of $C(Y)$ as a metric space, which is equipped with a Ricci-flat
Kähler metric $g_{C(Y)}$. Let $\mathcal{O}_{C(Y)}$ be the sheaf
of continuous functions on $C(Y)$ whose restriction to $\mathcal{R}_{C(Y)}$
is holomorphic. Our main theorem roughly states that the complex structure
of the regular part $\mathcal{R}_{C(Y)}$ extends to all of $C(Y)$. 
\begin{thm}
\label{thm1} The locally ringed space $(C(Y),\mathcal{O}_{C(Y)})$
is isomorphic to a normal affine algebraic variety, and $\mathcal{R}_{C(Y)}$
is biholomorphically identified with the regular points of the variety. 
\end{thm}

\begin{rem}
From \cite{complexsplitting}, it is known that the vector field $J\nabla f$
generates a one-parameter family of isometries of $C(Y)$. We will
show that there is a real torus $\mathbb{T}$ containing this one-parameter
family, which acts effectively and isometrically on $C(Y)$, such
that the embedding given by Theorem \ref{thm1} can be taken to be
$\mathbb{T}$-equivariant with respect to a linear torus action on
$\mathbb{C}^{N}$. 
\end{rem}

\begin{rem}
It is not important that $C(Y)$ arises as a tangent cone of $\mathcal{X}$.
While this is a natural way such Ricci flat cones may arise, Theorem
\ref{thm1} applies to any $\mathbb{F}$-limit of noncollapsed (at
all scales) Kähler-Ricci flows which happens to be a Ricci-flat cone
(defined at least for all negative time). 
\end{rem}

\begin{rem}
It is unclear whether $g_{C(Y)}$ is the restriction of a singular
Calabi-Yau metric on $C(Y)$ (as in \cite{singKE}), or even if $C(Y)$
has log-terminal singularities. 
\end{rem}

The general strategy for proving Theorem \ref{thm1} is similar to
that of \cite{szek2}: namely, we establish a Hörmander-type estimate
for the $\overline{\partial}$-operator on a singular soliton by arguing
on a sequence of smooth approximants. However, the failure of elliptic
estimates on time slices of a Ricci flow necessitates a new approach
to proving the $\overline{\partial}$-estimates, one which is more
adapted to the parabolic setting. We let $\mathcal{R}$ denote the
regular set of $\mathcal{X}$. 
\begin{thm}
\label{thm2} Suppose $t_{1}\in(-\infty,0)$, $v\in C_{c}^{\infty}(\mathcal{R}_{t_{1}},\mathbb{C})$,
$\varphi\in C_{c}^{\infty}(\mathcal{R}_{t_{1}})$ satisfy $\text{supp}(\varphi)\subset\subset\mathcal{R}_{t_{1}}$,
\[
Rc(\omega_{t_{1}})+\sqrt{-1}\partial\overline{\partial}(f_{t_{1}}+\varphi_{t_{1}})\geq\beta\omega_{t_{1}}
\]
in the sense of currents for some constant $\beta>0$, as well as
\[
\int_{\mathcal{R}_{t_{1}}}|\overline{\partial}v|^{2}e^{-\varphi}d\nu_{t_{1}}<\infty.
\]
 Then there exists $u\in C^{\infty}(\mathcal{R}_{t_{1}},\mathbb{C})$
satisfying $\overline{\partial}u=\overline{\partial}v$ and
\[
\int_{\mathcal{R}_{t_{1}}}|u|^{2}e^{-\varphi}d\nu_{t_{1}}\leq\frac{1}{\beta}\int_{\mathcal{R}_{t_{1}}}|\overline{\partial}v|^{2}e^{-\varphi}d\nu_{t_{1}}.
\]
\end{thm}

\begin{rem}
Note that this theorem also applies if $\mathcal{X}$ is replaced
with the static metric flow corresponding to any of its tangent cones. 
\end{rem}

We now outline the proof of Theorem \ref{thm2}, which is the main
technical result needed to prove that $C(Y)$ is an analytic variety.
Theorem \ref{thm2} follows if we can show that for any $\eta\in\mathcal{A}_{c}^{0,1}(\mathcal{R}_{t_{1}})$,
we have
\[
\left|\int_{\mathcal{R}_{t_{1}}}\langle\eta,\overline{\overline{\partial}v}\rangle e^{-\varphi}d\nu_{t_{1}}\right|\leq\left(\int_{\mathcal{R}_{t_{1}}}\frac{1}{\beta}|\overline{\partial}v|^{2}e^{-\varphi}d\nu_{t_{1}}\right)^{\frac{1}{2}}\left(\int_{\mathcal{R}_{t_{1}}}|\overline{\partial}_{f+\varphi}^{\ast}\eta|^{2}e^{-\varphi}d\nu_{t_{1}}\right)^{\frac{1}{2}}.
\]
For simplicity, we temporarily restrict our attention to the case
$\varphi\equiv0$, so that $\beta=\frac{1}{|t_{1}|}$. We write $\eta=\eta^{(1)}+\eta^{(2)}$,
where $\eta^{(1)}\in\ker(\overline{\partial})$ and $\eta^{(2)}$
is $L^{2}(\mathcal{R}_{t_{1}},d\nu_{t_{1}})$-orthogonal to $\ker(\overline{\partial})$.
If we knew the singular set had locally finite $(2n-4)$-dimensional
Minkowski content, and that relative volume comparison held on $\mathcal{X}_{t_{1}}$,
we could construct a good enough cutoff function to use
\[
\Delta|\eta^{(1)}|^{2}\geq|\nabla\eta^{(1)}|^{2}-C|\eta^{(1)}|^{2}
\]
and standard arguments for elliptic estimates as in Section 2 of \textit{\cite{szek2}
}\textit{\emph{to prove that $\eta^{(1)}$ is locally bounded near
the singular set, and then establish the desired identity by the appropriate
Kodaira-Nakano identity and integration by parts on $\mathcal{R}_{t_{0}}$.
To make do without this, we instead consider $\eta\in\mathcal{A}_{c}^{0,1}(\mathcal{R}_{t_{0}})$
for some $t_{0}<t_{1}$, and we work on the sequence $(M,(g_{i,t}))$
of smooth Kähler-Ricci flows approximating $\mathcal{X}$, taking
$\eta_{i,t_{0}}\in\mathcal{A}^{0,1}(M)$ and $v_{i,t_{1}}\in C^{\infty}(M,\mathbb{C})$
which converge smoothly to $\eta$ and $v$ respectively as $i\to\infty.$
We then let $(\eta_{i,t}^{(1)})_{t\in[t_{0},t_{1}]}$ be the forward
$\Delta_{\overline{\partial}}$ heat flow starting of $\eta_{i,t_{0}}^{(1)}$,
where $\eta_{i,t_{0}}^{(1)}$ is a weighted $L^{2}$-orthogonal projection
of $\eta_{i,t_{0}}$ to $\ker(\overline{\partial})$, and we let $(v_{i,t})_{t\in[t_{0},t_{1}]}$
solve a backwards $\Delta_{\overline{\partial},f_{i}}$-heat flow,
where $(2\pi\tau)^{-n}e^{-f_{i}}dg_{i,t}=\nu_{t}^{i}$ are the conjugate
heat kernels converging to $(\nu_{t})$. The flow for $\eta_{i}^{(1)}$
preserves the condition $\overline{\partial}\eta_{i}^{(1)}=0$, and
almost preserves the condition that $\eta_{i}^{(2)}$ is orthogonal
to $\ker(\overline{\partial})$, using the fact that this region of
our flow is almost-selfsimilar in the sense of \cite{bamlergen3}.
Again using almost-selfsimilarity, we can show that 
\[
t\mapsto\tau\int_{M}\langle\eta_{i,t}^{(1)},\overline{\overline{\partial}v}_{i,t}\rangle d\nu_{t}^{i}
\]
is almost constant. Moreover, the flow $(\eta_{i,t}^{(1)})$ regularizes
$\eta_{i,t_{0}}^{(1)}$, and $\eta_{i,t_{1}}^{(1)}$ is in fact locally
bounded; this would not be possible if we instead attempted to use
elliptic estimates as in \cite{szek2}, since we lack lower bounds
on the Ricci curvature. We can thus integrate by parts the appropriate
Kodaira-Nakano formula to get
\[
\left|\int_{M}\langle\eta_{i,t_{1}}^{(1)},\overline{\partial}v_{i,t_{1}}\rangle d\nu_{t_{1}}^{i}\right|\apprle\left(\int_{\text{supp}(\chi_{i})}|\overline{\partial}_{f_{i}}^{\ast}\eta_{i,t_{1}}|^{2}d\nu_{t_{1}}^{i}\right)^{\frac{1}{2}}\left(|t_{1}|\int_{M}|\overline{\partial}v_{i,t_{1}}|^{2}d\nu_{t_{1}}^{i}\right)^{\frac{1}{2}},
\]
where $\chi_{i}$ are appropriate cutoff functions. It remains to
show that 
\[
|t_{1}|^{2}\int_{\text{supp}(\chi_{i})}|\overline{\partial}_{f_{i}}^{\ast}\eta_{i,t_{1}}|^{2}d\nu_{t_{1}}^{i}\apprle|t_{0}|^{2}\int_{M}|\overline{\partial}_{f_{i}}^{\ast}\eta_{i,t_{0}}|^{2}d\nu_{t_{0}}^{i}.
\]
To accomplish this, we show that we can replace $f_{i}$ with the
strong almost-soliton potential functions $h_{i}$ constructed in
\cite{complexsplitting}. We can then prove that $|\overline{\partial}_{h_{i}}^{\ast}\eta_{i,t}|$
is almost a subsolution to the heat equation when our flow is almost-selfsimilar,
which also requires proving new estimates for $h_{i}$ that improve
those known for $f_{i}$. Finally, we take $i\to\infty$ to recover
\[
\left|\int_{\mathcal{R}_{t_{0}}}\langle\eta,\overline{\overline{\partial}v}_{t_{0}}\rangle d\nu_{t_{0}}\right|\leq\left(\int_{\mathcal{R}_{t_{0}}}\frac{1}{\beta}|\overline{\partial}v|^{2}d\nu_{t_{1}}\right)^{\frac{1}{2}}\left(\int_{\mathcal{R}_{t_{0}}}|\overline{\partial}_{f}^{\ast}\eta|^{2}d\nu_{t_{0}}\right)^{\frac{1}{2}},
\]
where $(v_{t})$ solves a backwards drift heat flow on $\mathcal{X}$
with final datum $v$, and then take $t_{0}\nearrow t_{1}$.}}

In Section 2, we review some facts from Bamler's compactness theory
and related constructions which will be needed in later sections,
and we establish notation. In Section 3, we establish evolution equations
satisfied by certain forwards and backwards heat flows appearing in
the proof of Theorem \ref{thm2}. We combine the ideas of the previous
sections in Section 4, proving the parabolic Hörmander's estimate
and establishing Theorem \ref{thm2}. In Section 5, we use Theorem
\ref{thm2} to prove that tangent cones are normal analytic varieties.
In Section 6, we use the parabolic frequency function and ideas from
\cite{szek1} to show that the notions of metric singular set and
complex analytic singular set agree. We also use the torus action
generated by the Reeb vector field of the cone to show that the tangent
cone is actually algebraic, completing the proof of Theorem \ref{thm1}.
In the Appendix, we establish properties of the heat flow on a Ricci
flat cone, and relate it to that of the corresponding metric flow,
which is used in Sections 4 and 5. \\

\noindent \textbf{Acknowledgments: }I thank Gábor Székelyhidi for
essential conversations which made this paper possible. I am also
grateful to Jian Song for useful suggestions and feedback. This material
is based upon work supported by the National Science Foundation under
grant DMS-2202980.

\section{Preliminaries and Notation}

Given parameters $\delta_{1},...,\delta_{n}$ and $\tau_{1},...,\tau_{m}$,
we will let 
\[
\Psi(\delta_{1},...,\delta_{n}|\tau_{1},...,\tau_{m})
\]
denote a quantity which can be made arbitrarily small for a fixed
$\tau_{1},...,\tau_{m}$ by taking $\delta_{1},...,\delta_{n}$ all
sufficiently small. If we say that a property depending on $\delta,\tau_{1},...,\tau_{m}$
holds assuming $\delta\leq\overline{\delta}(\tau_{1},...,\tau_{m})$,
it means that the property holds if $\delta$ is sufficiently small,
and the smallness only depends on the choice of $\tau_{1},...,\tau_{m}$. 

Given a Riemannian manifold $(M,g)$, we let $dg$ denote the Riemannian
volume measure. For a metric space $(X,d)$, and $x\in X$, $r>0$,
we let $B(x,r)$ be the open ball. We let $\nabla^{\mathbb{R}}$ denote
the usual Levi-Civita connection. If $(M,g,J)$ is Kähler, then we
will let $\omega=\omega_{g}\in\mathcal{A}^{1,1}(M)$ denote the corresponding
Kähler form, let $\nabla=\nabla^{1,0}$ and $\overline{\nabla}=\nabla^{0,1}$
denote the composition of the Levi-Civita connection with the projection
onto $(T^{1,0}M)^{\ast}$ and $(T^{0,1}M)^{\ast}$, respectively,
and we let $\Delta$ denote half the Riemannian Laplacian. 

Given $(p,q)$-tensors $S,T$, we set
\[
\langle S,\overline{T}\rangle:=g^{\overline{j}_{1}i_{1}}\cdots g^{\overline{j}_{p}i_{p}}g^{\overline{\ell}_{1}k_{1}}\cdots g^{\overline{\ell}_{q}k_{q}}S_{i_{1}\cdots i_{p}\overline{\ell}_{1}\cdots\overline{\ell}_{q}}\overline{T_{j_{1}\cdots j_{p}\overline{k}_{1}\cdots\overline{k}_{q}}},
\]
in a local holomorphic frame, and $|T|^{2}:=\langle T,\overline{T}\rangle$.
With this convention, we have $|\nabla^{\mathbb{R}}\phi|^{2}=2|\nabla\phi|^{2}$
for any $\phi\in C^{\infty}(M)$, as well as
\[
|\nabla^{\mathbb{R}}\nabla^{\mathbb{R}}\phi|^{2}=2|\nabla\overline{\nabla}\phi|^{2}+2|\nabla\nabla\phi|^{2}.
\]
The Dolbeault Laplacian $\Delta_{\overline{\partial}}=\overline{\partial}\overline{\partial}^{\ast}+\overline{\partial}^{\ast}\overline{\partial}$
and rough Laplacians $\overline{\nabla}^{\ast}\overline{\nabla},\nabla^{\ast}\nabla$
are nonnegative operators, so that on functions, we have $\Delta_{\overline{\partial}}=-\Delta$.
Given a function $f\in C^{\infty}(M)$, we let $\Delta_{\overline{\partial},f}$
denote the Dolbeault Laplacian $\Delta_{\overline{\partial},f}=\overline{\partial}\overline{\partial}_{f}^{\ast}+\overline{\partial}_{f}^{\ast}\overline{\partial}$
corresponding to the trivial bundle with metric $e^{-f}$. We also
let $R=R_{g}$ be the complex Laplacian, given by $R=g^{\overline{j}i}R_{i\overline{j}}$
in holomorphic coordinates, so that $R$ is half the Riemannian Laplacian.
Viewing $Rc$ as a $(1,1)$-tensor, $|Rc|^{2}$ is half of the Euclidean
norm squared of the Ricci curvature by our convention.

Suppose $(M^{n},(g_{t})_{t\in I},J)$ is a Kähler-Ricci flow on a
closed manifold, where $n$ is the complex dimension. We let $d_{g_{t}}$
be the corresponding Riemannian length metric. For $(x,t)\in M\times I$
and $r>0$, we define
\[
B(x,t,r):=B_{g}(x,t,r):=\{y\in M;d_{g_{t}}(x,y)<r\}.
\]
The curvature scale is defined by 
\[
r_{Rm}(x,t):=r_{Rm}^{g}(x,t):=\sup\{r>0;|Rm|(x',t')\leq r^{-2}\text{ for all }(x',t')\in B(x,t,r)\times([t-r^{2},t+r^{2}]\cap I)\}.
\]
We let $K(x,t;\cdot,\cdot)$ denote the fundamental solution of the
conjugate heat equation based at $(x,t)\in M\times I$, so that 
\[
(-\partial_{s}-\Delta_{g_{s}}+R)K(x,t;\cdot,\cdot)=0,
\]
and for any $h\in C^{\infty}(M)$, we have 
\[
\lim_{s\nearrow t}\int_{M}K(x,t;y,s)h(y)dg_{s}(y)=h(x).
\]
We write $d\nu_{x,t;s}=K(x,t;y,s)dg_{s}=(2\pi\tau)^{-n}e^{-f}dg_{s}$,
and write 
\[
\mathcal{N}_{x,t}(\tau):=\int_{M}fd\nu_{x,t;t-\tau}-n
\]
for the corresponding pointed Nash entropy. 

We recall the definition of almost-selfsimilarity from \cite{bamlergen3},
which describes when a region of Ricci flow is close to a singular
shrinking soliton in a certain sense. 
\begin{defn}
(Definition 5.1 of \cite{bamlergen3}) \label{SSdef} A point $(x_{0},t_{0})\in M\times I$
is $(\epsilon,r)$-selfsimilar if $[t_{0}-\epsilon^{-1}r^{2},t_{0}]\subseteq I$
and the following holds for $W:=\mathcal{N}_{x_{0},t_{0}}(r^{2})$:

$(i)$ $\int_{t_{0}-\epsilon^{-1}r^{2}}^{t_{0}-\epsilon r^{2}}\int_{M}\tau\left(\left|Rc+\nabla\overline{\nabla}f-\frac{1}{\tau}g\right|^{2}+|\nabla\nabla f|^{2}\right)d\nu_{x_{0},t_{0};t}dt\leq\epsilon,$

$(ii)$$\int_{M}\left|\tau(R+2\Delta f-|\nabla f|^{2})+f-2n-W\right|d\nu_{x_{0},t_{0};t}\leq\epsilon$
for all $t\in[t_{0}-\epsilon^{-1}r^{2},t_{0}-\epsilon r^{2}]$,

$(iii)$ $r^{2}R\geq-\epsilon$ on $M\times[t_{0}-\epsilon^{-1}r^{2},t_{0}-\epsilon r^{2}]$.
\end{defn}

\begin{rem}
This definition slightly differs from that of \cite{bamlergen3} due
to notational differences as well as the usual convention of defining
Kähler-Ricci flow without a factor of two attached to $Rc$. 
\end{rem}

In \cite{complexsplitting}, a stronger notion of almost selfsimilarity
is given.
\begin{defn}
(Definition 3.1 of \cite{complexsplitting}) \label{SSSdef} A strong
$(\epsilon,r)$-soliton potential based at $(x_{0},t_{0})\in M\times I$
is a function $h\in C^{\infty}(M\times[t_{0}-\epsilon^{-1}r^{2},t_{0}-\epsilon r^{2}])$
such that the following holds for $W:=\mathcal{N}_{x_{0},t_{0}}(r^{2})$:

$(i)$ $(\partial_{t}-\Delta)\left(\tau(h-W-n)\right)=0$,

$(ii)$ $r^{-2}\int_{t_{0}-\epsilon^{-1}r^{2}}^{t_{0}-\epsilon r^{2}}\int_{M}\left|\tau(R+|\nabla h|^{2})-(h-W)\right|d\nu_{x_{0},t_{0};t}dt\leq\epsilon$,

$(iii)$ $\int_{t_{0}-\epsilon^{-1}r^{2}}^{t_{0}-\epsilon r^{2}}\int_{M}\left(|Rc+\nabla\overline{\nabla}h-\frac{1}{\tau}g|^{2}+|\nabla\nabla h|^{2}\right)d\nu_{x_{0},t_{0};t}dt\leq\epsilon$,

$(iv)$ $\sup_{t\in[t_{0}-\epsilon^{-1}r^{2},t_{0}-\epsilon r^{2}]}\int_{M}|\tau(R+2\Delta h-|\nabla h|^{2})+h-2n-W|d\nu_{x_{0},t_{0};t}\leq\epsilon,$

$(v)$ $\int_{M}(h-n-W)d\nu_{x_{0},t_{0};t}=0$ for all $t\in[t_{0}-\epsilon^{-1}r^{2},t_{0}-\epsilon r^{2}]$,

$(vi)$ $r^{2}R\ge-\epsilon$ on $M\times[t_{0}-\epsilon^{-1}r^{2},t_{0}-\epsilon r^{2}]$. 

\noindent We call $(x_{0},t_{0})$ strongly $(\epsilon,r)$-selfsimilar
if it is $(\epsilon,r)$-selfsimilar and if such a function $h$ exists. 
\end{defn}

The usefulness of this definition is a consequence of the following
proposition, which roughly states that almost selfsimilarity implies
strong almost-selfsimilarity.
\begin{prop}
\label{strongselfsimilar} (Proposition 3.2 of \cite{complexsplitting})
For any $\epsilon>0$, $Y<\infty$, the following holds whenever $\delta\leq\overline{\delta}(\epsilon,Y)$.
Suppose $(M^{n},(g_{t})_{t\in I})$ is a closed Ricci flow, assume
$(x_{0},t_{0})\in M\times I$ is $(\delta,r)$-selfsimilar with $W:=\mathcal{N}_{x_{0},t_{0}}(r^{2})\geq-Y$,
and write $d\nu_{x_{0},t_{0};t}=(2\pi\tau)^{-\frac{n}{2}}e^{-f}dg_{t}$.
Then there exists a strong $(\epsilon,r)$-soliton potential $h\in C^{\infty}(M\times[t_{0}-\epsilon^{-1}r^{2},t_{0}-\epsilon r^{2}])$
based at $(x_{0},t_{0})$, such that
\[
\int_{t_{0}-\epsilon^{-1}r^{2}}^{t_{0}-\epsilon r^{2}}\int_{M}|\nabla(f-h)|^{2}d\nu_{x_{0},t_{0};t}dt+\sup_{t\in[t_{0}-\epsilon^{-1}r^{2},t_{0}-\epsilon r^{2}]}\int_{M}|f-h|^{2}d\nu_{x_{0},t_{0};t}\leq\epsilon.
\]
\end{prop}

Next, we establish an improvement of Proposition \ref{strongselfsimilar}
which will be used in the proof of Lemma \ref{Hormander}. 
\begin{lem}
\label{gradclose} Suppose $(M^{n},(g_{t})_{t\in I})$ is a closed
Kähler-Ricci flow and $(x_{0},0)$ is strongly $(\delta,1)$-selfsimilar,
with $h\in C^{\infty}(M\times[-\delta^{-1},-\delta])$ a strong almost-soliton
potential function as in Proposition \ref{strongselfsimilar}. Writing
$d\nu_{t}=d\nu_{x_{0},0;t}$, we then have
\[
\sup_{t\in[-\epsilon^{-1},-\epsilon]}\int_{M}|\nabla(f-h)|^{2}d\nu_{t}\leq\Psi(\delta|\epsilon).
\]
\end{lem}

\begin{proof}
Integration by parts gives
\begin{align*}
\int_{M}|\nabla(f-h)|^{2}d\nu_{t}= & 2\text{Re}\int_{M}\langle\nabla f,\overline{\nabla}(f-h)\rangle d\nu_{t}+\int_{M}\left(|\nabla h|^{2}-|\nabla f|^{2}\right)d\nu_{t}\\
= & 2\int_{M}\Delta(f-h)d\nu_{t}+\int_{M}\left(|\nabla h|^{2}-|\nabla f|^{2}\right)d\nu_{t}\\
= & \int_{M}(R+2\Delta f-|\nabla f|^{2})d\nu_{t}-\int_{M}(R+2\Delta h-|\nabla h|^{2})d\nu_{t}.
\end{align*}
On the other hand, Definitions \ref{SSdef}, \ref{SSSdef} give
\[
\sup_{t\in[-\delta^{-1},-\delta]}\int_{M}\left|\tau(R+2\Delta f-|\nabla f|^{2})+f-n-W\right|d\nu_{t}<\delta,
\]
\[
\sup_{t\in[-\delta^{-1},-\delta]}\int_{M}\left|\tau(R+2\Delta h-|\nabla h|^{2})+h-n-W\right|d\nu_{t}<\delta,
\]
so combining expressions gives
\begin{align*}
\int_{M}|\nabla(f-h)|^{2}d\nu_{t}\leq & \frac{1}{\tau}\int_{M}|f-h|d\nu_{t}+2\delta\leq\Psi(\delta|\epsilon),
\end{align*}
where the last estimate follows from Proposition \ref{strongselfsimilar}. 
\end{proof}
We now recall some additional notions and notations from \cite{bamlergen2,bamlergen3},
giving specific references when the details are omitted. Let $I\subseteq\mathbb{R}$.
A metric flow (Definition 3.2 in \cite{bamlergen2}) over $I$ is
a tuple $(\mathcal{X},\mathfrak{t},(d_{t})_{t\in I},(\nu_{x;s})_{x\in\mathcal{X},s\in(-\infty,\mathfrak{t}(x)]\cap I})$,
where $\mathcal{X}$ is a set, $\mathfrak{t}:\mathcal{X}\to I$ is
a (time) function, $d_{t}$ are complete separable metrics on $\mathcal{X}_{t}:=\mathfrak{t}^{-1}(t)$,
and $\nu_{x;s}$ are Borel probability measures on $\mathcal{X}_{s}$,
which satisfy a time-inhomogeneous version of the semigroup property
and a gradient estimate. The regular part of $\mathcal{X}$ is a dense
open subset $\mathcal{R}$ with the structure of a Ricci flow spacetime
$(\mathcal{R},g,\mathfrak{t},\partial_{\mathfrak{t}})$ in the sense
of Definition 1.2 of \cite{klott1}. Any smooth Ricci flow $(M,(g_{t})_{t\in I})$
admits the structure of a metric flow, by letting $\mathcal{X}:=M\times I$,
taking $\mathfrak{t}:\mathcal{X}\to\mathbb{R}$ to be projection onto
the second factor, and setting $d_{t}:=d_{g_{t}}$ and $\nu_{(x,t);s}:=\nu_{x,t;s}=K(x,t;\cdot,s)dg_{s}$.
A metric flow pair is roughly a pair $(\mathcal{X},(\mu_{t})_{t\in I})$
consisting of a metric flow with a reference conjugate heat kernel. 

The appropriate notion of convergence of metric flows is that of $\mathbb{F}$-convergence,
defined in Definition 5.6 of \cite{bamlergen2}. The notion of a correspondence
of metric flows will be called upon in later sections, so we recall
its full definition to establish notation. 
\begin{defn}
(Definition 5.4 of \cite{bamlergen2}) Given a sequence of Ricci flows
with a choice of reference conjugate heat flow $(M_{i},(g_{i,t})_{t\in[-T_{i},0]},(\mu_{t}^{i})_{t\in[-T_{i},0]})$,
and a metric flow pair $(\mathcal{X},(\mu_{t}^{\infty})_{t\in[-T_{\infty},0]})$
where $T_{\infty}:=\lim_{i\to\infty}T_{i}$, a correspondence is a
pair of the form 
\[
\mathfrak{C}=\left((Z_{t},d_{t})_{t\in I'},(\varphi_{t}^{i})_{t\in[T_{i},0],i\in\mathbb{N}\cup\{\infty\}}\right),
\]
where $(Z_{t},d_{t}^{Z})$ are metric spaces, $I'\subseteq\mathbb{R}$
is an interval containing $[-T_{\infty},0]$, and $\varphi_{t}^{i}:(M_{i},d_{g_{i,t}})\to(Z_{t},d_{t}^{Z})$,
$\varphi_{t}^{\infty}:(\mathcal{X}_{t},d_{t})\to(Z_{t},d_{t})$ are
isometric embeddings. 
\end{defn}

\begin{rem}
\label{convergenotation} The notation
\[
(M_{i},(g_{i,t})_{t\in[-T_{i},0]},(\mu_{t}^{i})_{t\in[-T_{i},0]})\xrightarrow[i\to\infty]{\mathbb{F},\mathfrak{C}}(\mathcal{X},(\mu_{t}^{\infty})_{t\in[-T_{\infty},0]})
\]
means that $\mathbb{F}$-convergence occurs within the correspondence
$\mathfrak{C}$ (see Definition 6.2 in \cite{bamlergen2}). By Theorem
9.1 of \cite{bamlergen2} and Theorem 2.5 of \cite{bamlergen3}, we
know that then the convergence is smooth on the regular part. This
means that there are time-preserving diffeomorphisms $\psi_{i}:U_{i}\to V_{i}\subseteq M_{i}\times[-T_{i},0]$,
where $(U_{i})_{i\in\mathbb{N}}$ is a precompact exhaustion of the
regular set $\mathcal{R}$ of $\mathcal{X}$ and $V_{i}\subseteq M_{i}\times[-T_{0},0]$
are open, such that the following hold:

$(i)$ $\psi_{i}^{\ast}g_{i}\to g$ in $C_{\text{loc}}^{\infty}(\mathcal{R}$),

$(ii)$ $(\psi_{i}^{-1})_{\ast}\partial_{t}\to\partial_{\mathfrak{t}}$
in $C_{\text{loc}}^{\infty}(\mathcal{R})$,

$(iii)$ $\psi_{i}^{\ast}v_{i}\to v_{\infty}$ in $C_{\text{loc}}^{\infty}(\mathcal{R})$
, where $d\mu_{t}^{i}=v_{i,t}=dg_{i,t}$, and $d\mu_{t}^{\infty}=v_{\infty,t}dg_{t}$
on $\mathcal{R}$,

$(iv)$ $(\psi_{i},\psi_{i})^{\ast}K^{i}(\cdot;\cdot)\to K$ in $C_{\text{loc}}^{\infty}(\{(x,y)\in\mathcal{R}\times\mathcal{R};\mathfrak{t}(x)>\mathfrak{t}(y)\})$,
where $K^{i}$ are the conjugate heat kernels of $(M_{i},(g_{i,t}))$,
and $d\nu_{x;s}=K(x;\cdot)dg_{s}$ for $x\in\mathcal{X}_{t}$, $s\in(-T_{\infty},0]$,

$(v)$ For any compact subset $L\subseteq\mathcal{R}_{[t_{0},t_{1}]}$,
we have 
\[
\lim_{i\to\infty}\sup_{t\in[t_{0},t_{1}]}\sup_{x\in L\cap\mathcal{R}_{t}}d_{t}^{Z}\left((\varphi_{t}^{i}\circ\psi_{i,t})(x),\varphi_{t}^{\infty}(x)\right)=0,
\]
where $\psi_{i,t}$ is the restriction of $\psi_{i}$ to $U_{i}\cap\mathcal{R}_{t}$,

$(vi)$ Given $(y_{i},t_{i})\in M_{i}$ and $y\in\mathcal{R}_{t_{\infty}}$,
we have $(y_{i},t_{i})\to y$ within $\mathfrak{C}$ (see Definition
6.18 of \cite{bamlergen2}) if and only if $y_{i}\in V_{i}$ for sufficiently
large $i\in\mathbb{N}$ and $\psi_{i}^{-1}(y_{i})\to y$ in $\mathcal{R}$. 

Note that $(v),(vi)$ tell us in some sense that the smooth convergence
is compatible with the correspondence.
\end{rem}

If $(M_{i},(g_{i,t}))$ are Kähler-Ricci flows, then the regular part
$\mathcal{R}$ of the limit $\mathcal{X}$ also be equipped with a
complex structure on its time slices. For this reason, we make the
following definition.
\begin{defn}
(c.f. Definition 1.2 of \cite{klott1}) A Kähler-Ricci flow spacetime
is a tuple $(\mathcal{R},\mathfrak{t},\partial_{\mathfrak{t}},g,J)$,
where $\mathcal{R}$ is a smooth manifold, $\mathfrak{t}\in C^{\infty}(\mathcal{R})$
is a submersion, $\partial_{\mathfrak{t}}\in\mathfrak{X}(\mathcal{R})$
satisfies $\partial_{\mathfrak{t}}\mathfrak{t}=1$, $g$ is a smooth
bundle metric on $\ker(d\mathfrak{t})$, $J\in\text{End}(\ker(d\mathfrak{t}))$
restricts to an integrable almost-complex structure on each fiber
$\mathfrak{t}^{-1}(t)$, and 
\[
\mathcal{L}_{\partial_{\mathfrak{t}}}g=-Rc(g|_{\mathfrak{t}^{-1}(t)}),\,\,\,\,\,\,\,\,\,\,\,\,\,\,\,\,\,\,\,\,\mathcal{L}_{\partial_{\mathfrak{t}}}J=0.
\]
\end{defn}

\section{Estimates for Heat Flows of $(0,1)$-Forms}

In this section, we establish estimates and evolution equations for
certain forwards and backwards heat flows coupled to the Kähler-Ricci
flow. We choose to study the Dolbeault Laplacian heat flow due to
the fact it preserves $\overline{\partial}$-closed forms, and the
norm of solutions are subsolutions to the usual heat equation. We
also consider a backwards heat flow of the drift Laplacian, which
is chosen so that Proposition \ref{slowchange} holds. 

Throughout this section, suppose $(M^{n},(g_{t})_{t\in I},J)$ is
a closed Kähler-Ricci flow and that $\eta_{t}\in\mathcal{A}^{0,1}(M)$,
$v_{t}\in C^{\infty}(M,\mathbb{C})$, $\varphi_{t}\in C^{\infty}(M)$
satisfy 
\begin{equation}
\partial_{t}\eta_{t}=-\overline{\nabla}^{\ast}\overline{\nabla}\eta_{t}-Rc_{g_{t}}(\eta_{t})=-\Delta_{\overline{\partial}}\eta_{t},\label{forward}
\end{equation}
\begin{equation}
-\partial_{t}v_{t}=\Delta v_{t}-\langle\nabla^{\mathbb{R}}v_{t},\nabla^{\mathbb{R}}(f_{t}+\varphi_{t})\rangle=-\Delta_{\overline{\partial},f_{t}+\varphi_{t}}v_{t}-\nabla_{\nabla^{1,0}(f_{t}+\varphi_{t})}v_{t},\label{backward}
\end{equation}
where $\Delta_{\overline{\partial}}\eta=(\overline{\partial}\overline{\partial}^{\ast}+\overline{\partial}^{\ast}\overline{\partial})\eta=\overline{\nabla}^{\ast}\overline{\nabla}\eta+Rc(\eta)$,
so that in local holomorphic coordinates, we can write
\[
\partial_{t}\eta_{\overline{j}}=g^{\overline{q}p}\left(\nabla_{p}\nabla_{\overline{q}}\eta_{\overline{j}}-R_{p\overline{j}}\eta_{\overline{q}}\right),
\]
\[
\partial_{t}v=-g^{\overline{q}p}(\nabla_{p}\nabla_{\overline{q}}v-\nabla_{p}v\nabla_{\overline{q}}(f+\varphi)-\nabla_{p}(f+\varphi)\nabla_{\overline{q}}v).
\]
Also fix $h\in C^{\infty}(M\times I)$, and let $\overline{\partial}_{h}^{\ast}\eta$
denote the formal adjoint of $\overline{\partial}$ with respect to
the weight $e^{-h}dg_{t}$:
\[
\overline{\partial}_{h}^{\ast}\eta:=-g^{\overline{j}i}\left(\nabla_{i}\eta_{\overline{j}}-\eta\nabla_{i}h\right).
\]
For ease of notation, we set
\[
w:=\tau\left(\tau(R+2\Delta h-|\nabla h|^{2})+h-2n-W\right),
\]
where $W\in(-\infty,0].$ 

We now compute some evolution equations that will be essential for
proving Hörmander's $L^{2}$ estimate.
\begin{lem}
\label{Computations} $(i)$ We have
\[
(\partial_{t}-\Delta)|\eta|^{2}=-|\nabla\eta|^{2}-|\overline{\nabla}\eta|^{2}.
\]

\noindent $(ii)$ We have 
\begin{align*}
(\partial_{t}-\Delta)(\tau\overline{\partial}_{h}^{\ast}\eta)= & -\tau\left\langle Rc+\nabla\overline{\nabla}h-\frac{1}{\tau}g,\overline{\nabla\eta}\right\rangle -\tau\langle\nabla\nabla h,\overline{\nabla}\eta\rangle+\tau\left\langle \nabla\left((\partial_{t}-\Delta)h-\frac{1}{\tau}\left(h-n-W\right)\right),\eta\right\rangle .
\end{align*}
In particular, if $h$ is a strong almost-soliton potential function
based at $(x_{0},0)$, and $W:=\mathcal{W}_{x_{0},0}(1)$, then 
\[
(\partial_{t}-\Delta)|\tau\overline{\partial}_{h}^{\ast}\eta|\leq\tau|Rc+\nabla\overline{\nabla}h-\frac{1}{\tau}g|\cdot|\nabla\eta|+\tau|\nabla\nabla h|\cdot|\overline{\nabla}\eta|.
\]

\noindent $(iii)$ If $h$ is a strong almost-soliton potential function
based at $(x_{0},0)$, and $W:=\mathcal{W}_{x_{0},0}(1)$, then 
\begin{align*}
(\partial_{t}-\Delta)w= & \tau^{2}\left|Rc+\nabla\overline{\nabla}h-\frac{1}{\tau}g\right|^{2}+\tau^{2}|\nabla\nabla h|^{2}.
\end{align*}
\end{lem}

\begin{proof}
$(i)$ We have
\begin{align*}
(\partial_{t}-\Delta)|\eta|^{2}= & (\partial_{t}-g^{\overline{q}p}\nabla_{p}\nabla_{\overline{q}})\left(g^{\overline{j}i}\eta_{\overline{j}}\overline{\eta_{\overline{i}}}\right)\\
= & R^{\overline{j}i}\eta_{\overline{j}}\overline{\eta_{\overline{i}}}+2\text{Re}\langle\partial_{t}\eta,\overline{\eta}\rangle-g^{\overline{q}p}g^{\overline{j}i}\left(\overline{\eta_{\overline{i}}}\nabla_{p}\nabla_{\overline{q}}\eta_{\overline{j}}+\eta_{\overline{j}}\overline{\nabla_{\overline{p}}\nabla_{q}\eta_{\overline{i}}}+\nabla_{p}\eta_{\overline{j}}\overline{\nabla_{q}\eta_{\overline{i}}}+\nabla_{\overline{q}}\eta_{\overline{j}}\overline{\nabla_{\overline{p}}\eta_{\overline{i}}}\right)\\
= & 2\text{Re}\langle\partial_{t}\eta,\overline{\eta}\rangle+g^{\overline{j}i}(\Delta_{\overline{\partial}}\eta)_{\overline{j}}\overline{\eta_{\overline{i}}}-g^{\overline{q}p}g^{\overline{j}i}\eta_{\overline{j}}\overline{\left(\nabla_{q}\nabla_{\overline{p}}\eta_{\overline{i}}-R_{q\overline{p}k\overline{i}}g^{\overline{\ell}k}\eta_{\overline{\ell}}\right)}-|\nabla\eta|^{2}-|\overline{\nabla}\eta|^{2}\\
= & 2\text{Re}\langle\partial_{t}\eta+\Delta_{\overline{\partial}}\eta,\overline{\eta}\rangle-|\nabla\eta|^{2}-|\overline{\nabla}\eta|^{2}\\
= & -|\nabla\eta|^{2}-|\overline{\nabla}\eta|^{2}.
\end{align*}
$(ii)$ We compute
\begin{align*}
\partial_{t}(\overline{\partial}_{h}^{\ast}\eta)= & \partial_{t}\left(-g^{\overline{j}i}(\nabla_{i}\eta_{\overline{j}}-\eta_{\overline{j}}\nabla_{i}h)\right)\\
= & -\langle Rc,\overline{\nabla\eta-\nabla h\otimes\eta}\rangle+\overline{\partial}_{h}^{\ast}(\partial_{t}\eta)+\langle\nabla\partial_{t}h,\eta\rangle,
\end{align*}
\begin{align*}
-\Delta(\overline{\partial}_{h}^{\ast}\eta)= & g^{\overline{q}p}g^{\overline{j}i}\left(\nabla_{p}\nabla_{\overline{q}}\nabla_{i}\eta_{\overline{j}}-\eta_{\overline{j}}\nabla_{p}\nabla_{\overline{q}}\nabla_{i}h-\nabla_{p}\nabla_{i}h\nabla_{\overline{q}}\eta_{\overline{j}}-\nabla_{\overline{q}}\nabla_{i}h\nabla_{p}\eta_{\overline{j}}-\nabla_{i}h\nabla_{p}\nabla_{\overline{q}}\eta_{\overline{j}}\right)\\
= & g^{\overline{q}p}g^{\overline{j}i}\nabla_{p}\left(\nabla_{i}\nabla_{\overline{q}}\eta_{\overline{j}}-R_{\overline{q}i\overline{j}k}g^{\overline{\ell}k}\eta_{\overline{\ell}}\right)-\langle\nabla\Delta h,\eta\rangle-\langle\nabla\nabla h,\overline{\nabla}\eta\rangle-\langle\nabla\overline{\nabla}h,\overline{\nabla\eta}\rangle\\
 & -g^{\overline{j}i}\nabla_{i}h\left(g^{\overline{q}p}(\nabla_{p}\nabla_{\overline{q}}\eta_{\overline{j}}-R_{p\overline{j}}\eta_{\overline{q}})\right)-g^{\overline{j}i}R_{p\overline{j}}\eta_{\overline{q}}\nabla_{i}h\\
= & g^{\overline{j}i}\nabla_{i}\left(g^{\overline{q}p}\nabla_{p}\nabla_{\overline{q}}\eta_{\overline{j}}\right)-g^{\overline{q}p}\nabla_{p}\left(R_{k\overline{q}}g^{\overline{\ell}k}\eta_{\overline{\ell}}\right)-\langle\nabla\Delta h,\eta\rangle-\langle\nabla\nabla h,\overline{\nabla}\eta\rangle-\langle\nabla\overline{\nabla}h,\overline{\nabla\eta}\rangle\\
 & -\langle\nabla h,\Delta_{\overline{\partial}}\eta\rangle-\langle Rc,\overline{\nabla h\otimes\eta}\rangle.
\end{align*}
By the Kodaira-Nakano formula, 
\[
g^{\overline{j}i}\nabla_{i}\left(g^{\overline{q}p}\nabla_{p}\nabla_{\overline{q}}\eta_{\overline{j}}\right)-g^{\overline{j}i}\nabla_{i}\left(R_{p\overline{j}}g^{\overline{q}p}\eta_{\overline{q}}\right)=g^{\overline{j}i}\nabla_{i}(\Delta_{\overline{\partial}}\eta)_{\overline{j}}=-\overline{\partial}_{h}^{\ast}(\Delta_{\overline{\partial}}\eta)+\langle\nabla h,\Delta_{\overline{\partial}}\eta\rangle.
\]
Combining expressions gives
\begin{align*}
(\partial_{t}-\Delta)(\overline{\partial}_{h}^{\ast}\eta)= & -\langle Rc,\overline{\nabla\eta-\nabla h\otimes\eta}\rangle+\langle\nabla(\partial_{t}-\Delta)h,\eta\rangle+\overline{\partial}_{h}^{\ast}(\partial_{t}\eta+\Delta_{\overline{\partial}}\eta)\\
 & -\langle\nabla\nabla h,\overline{\nabla}\eta\rangle-\langle\nabla\overline{\nabla}h,\overline{\nabla\eta}\rangle-\langle Rc,\overline{\nabla h\otimes\eta}\rangle\\
= & \langle\nabla(\partial_{t}-\Delta)h,\eta\rangle-\left\langle Rc+\nabla\overline{\nabla}h,\overline{\nabla\eta}\right\rangle -\langle\nabla\nabla h,\overline{\nabla}\eta\rangle\\
= & \frac{1}{\tau}\overline{\partial}_{h}^{\ast}\eta+\langle\nabla(\partial_{t}-\Delta)h,\eta\rangle-\frac{1}{\tau}\langle\nabla h,\eta\rangle-\left\langle Rc+\nabla\overline{\nabla}h-\frac{1}{\tau}g,\overline{\nabla\eta}\right\rangle -\langle\nabla\nabla h,\overline{\nabla}\eta\rangle.
\end{align*}
$(iii)$ This is a straightforward computation (see Proposition 3.2
of \cite{complexsplitting}).
\end{proof}
\begin{cor}
\label{Essential} For any $a>0$, we have 
\[
(\partial_{t}-\Delta)\left(-aw+\frac{1}{a}|\eta|^{2}+2|\tau\overline{\partial}_{h}^{\ast}\eta|\right)\leq0.
\]
\end{cor}

\begin{proof}
Combine $(i)-(iii)$ of the previous lemma. 
\end{proof}
\begin{rem}
Though it will not be needed, we observe that $\overline{\partial}\eta$
also satisfies a nice evolution equation, which does not rely on almost-selfsimilarity:
\[
(\partial_{t}-\Delta)|\overline{\partial}\eta|^{2}\leq-|\overline{\nabla}\overline{\partial}\eta|^{2}-|\nabla\overline{\partial}\eta|^{2}.
\]
\end{rem}

Next, we derive a formula which will be used to show that if $(\eta_{t})$
solves (\ref{forward}) and $(v_{t})$ solves (\ref{backward}), then
the heat kernel weighted $L^{2}$-inner product of $\eta$ and $\overline{\partial}v$
is almost constant in time. 
\begin{prop}
\label{slowchange} If $h$ is a strong almost-soliton potential function,
then 
\begin{align*}
\frac{d}{dt}\int_{M}(\tau\overline{\partial}_{h+\varphi}^{\ast}\eta)\overline{v}e^{-\varphi}d\nu_{t}= & -\int_{M}\tau\left\langle Rc+\nabla\overline{\nabla}h-\frac{1}{\tau}g,\overline{\nabla\eta}\right\rangle \overline{v}e^{-\varphi}d\nu_{t}-\int_{M}\tau\langle\nabla\nabla h,\overline{\nabla}\eta\rangle\overline{v}e^{-\varphi}d\nu_{t}\\
 & -\int_{M}\tau\left(\langle\nabla\overline{\nabla}\varphi,\overline{\nabla\eta}\rangle+\langle\nabla\nabla\varphi,\overline{\nabla}\eta\rangle\right)\overline{v}e^{-\varphi}d\nu_{t}\\
 & +\int_{M}\tau\left\langle \nabla\left((\partial_{t}-\Delta)\varphi-\frac{1}{\tau}\varphi\right),\eta\right\rangle \overline{v}e^{-\varphi}d\nu_{t}.
\end{align*}
\end{prop}

\begin{proof}
This follows by combining Lemma \ref{Computations} with 
\[
(-\partial_{t}-\Delta+R)\left((2\pi\tau)^{-n}e^{-(\varphi+f)}v\right)=0.
\]
\end{proof}
We will also record the following elementary estimates for $\eta$,
the first of which will allow us to make use of Proposition \ref{slowchange}.
\begin{lem}
\label{veryrough} Suppose $(\eta_{t})_{t\in[t_{0},t_{1}]}$ solves
(\ref{forward}). 

$(i)$ 
\[
\sup_{t\in[t_{0},t_{1}]}\frac{1}{\tau^{2}}\int_{M}|\eta_{t}|^{2}d\nu_{t}+\int_{t_{0}}^{t_{1}}\frac{1}{\tau^{2}}\int_{M}(|\nabla\eta_{t}|^{2}+|\overline{\nabla}\eta_{t}|^{2})d\nu_{t}dt\leq\frac{1}{t_{0}^{2}}\int_{M}|\eta_{t_{0}}|^{2}d\nu_{t_{0}}.
\]

$(ii)$ For any $t\in(t_{0},t_{1}]$ and $p\in[2,\infty)$ such that
$|t_{0}|\geq(p-1)|t|$, we have
\[
\left(\int_{M}|\eta_{t}|^{p}d\nu_{t}\right)^{\frac{1}{p}}\leq\left(\int_{M}|\eta_{t_{0}}|^{2}d\nu_{t_{0}}\right)^{\frac{1}{2}}.
\]
\end{lem}

\begin{proof}
These follows from formula $(i)$ in Lemma \ref{Computations}, the
fact that $|\eta_{t}|$ is a subsolution of the heat equation, and
Theorem 12.1 of \cite{bamlergen3}.
\end{proof}

\section{A Parabolic Version of Hörmander's Estimates}

Suppose $(M_{i},(g_{i,t})_{t\in[-\epsilon_{i}^{-1},0]},J_{i},(\nu_{x_{i},0;t}^{i})_{t\in[-\epsilon_{i}^{-1},0]})$
is a sequence of closed Kähler-Ricci flows, where $(x_{i},0)$ are
strongly $(\epsilon_{i},1)$-selfsimilar, with $h_{i}\in C^{\infty}(M_{i}\times[-\epsilon_{i}^{-1},-\epsilon_{i}])$
strong $(\epsilon_{i},1)$-almost soliton potential as in Proposition
\ref{strongselfsimilar}, and $\epsilon_{i}\searrow0$. Also assume
there exists $Y<\infty$ such that $\mathcal{N}_{x_{0},0}^{g_{i}}(\epsilon_{i}^{-1})\geq-Y$
for all $i\in\mathbb{N}$ (this is satisfied given the assumptions
of Theorem \ref{thm1}). After passing to a subsequence, we may assume
that there is a correspondence $\mathfrak{C}=\left((Z_{t},d_{t})_{t\in(-\infty,0]},(\varphi_{t}^{i})_{t\in[-\epsilon_{i}^{-1},0]}\right)$
such that 
\begin{equation}
(M_{i},(g_{i,t})_{t\in[-\epsilon_{i}^{-1},0]},J_{i},(\nu_{x_{i},0;t}^{i})_{t\in[-\epsilon_{i}^{-1},0]})\xrightarrow[i\to\infty]{\mathbb{\mathbb{F}},\mathfrak{C}}(\mathcal{X},(\nu_{t})_{t\in(-\infty,0]}),\label{F}
\end{equation}
where $(\mathcal{X},(\nu_{t})_{t\in(-\infty,0]})$ is an $H_{2n}$-concentrated
metric flow pair, and the convergence holds uniformly on compact time
intervals of $(-\infty,0)$. Moreover, $(\mathcal{X},(\nu_{t})_{t\in(-\infty,0)})$
is a metric soliton, and there is a decomposition $\mathcal{X}=\mathcal{R}\sqcup\mathcal{S}$,
where $\dim_{\mathcal{H}}(\mathcal{S}_{t})\leq2n-4$ for all $t\in(-\infty,0)$,
and $\mathcal{R}$ admits the structure of a Kähler-Ricci flow spacetime
$(\mathcal{R},g,\mathfrak{t},\partial_{\mathfrak{t}},J)$. If we define
$f\in C^{\infty}(\mathcal{R})$ by $d\nu_{t}=(2\pi\tau)^{-n}e^{-f}dg_{t}$,
then we have 
\[
Rc_{g}+\nabla^{\mathbb{R}}\nabla^{\mathbb{R}}f=\frac{1}{\tau}g,\,\,\,\,\,\,\,\,\,\,\,\,\,\,\,\,\,\,\,\mathcal{L}_{\nabla f}J=0.
\]
Moreover, there is a precompact open exhaustion $(U_{i})$ of $\mathcal{R}$
along with time-preserving open embeddings $\psi_{i}:U_{i}\to M_{i}\times[-\epsilon_{i}^{-1},-\epsilon_{i}]$
such that 
\[
\psi_{i}^{\ast}h_{i},\psi_{i}^{\ast}f_{i}\to f,\,\,\,\,\,\,\,\,\,\psi_{i}^{\ast}g_{i}\to g,\,\,\,\,\,\,\,\,\,\psi_{i}^{\ast}J_{i}\to J,\,\,\,\,\,\,\,\,\,(\psi_{i}^{-1})_{\ast}\partial_{t}=\partial_{\mathfrak{t}}
\]
in $C_{\text{loc}}^{\infty}(\mathcal{R})$, where we define $f_{i}\in C^{\infty}(M_{i}\times[-\epsilon_{i}^{-1},0))$
by $d\nu_{x_{i},0;t}^{i}=(2\pi\tau)^{-n}e^{-f_{i}}dg_{i,t}$. 

The following estimate is an application of Proposition 8.1 of \cite{bamlergen3}.
It will be used to show that an $L^{2}$ solution $(\eta_{t})$ to
(\ref{forward}) is locally bounded after some amount of time has
passed. This regularization feature of the flow (\ref{forward}) is
needed to integrate the Kodaira-Nakano formula by parts. 
\begin{lem}
\label{nearby} Given $-\infty<t_{0}<t_{1}<0$, $\alpha\in(0,\frac{1}{2})$,
and any compact subset $K\subseteq\mathcal{X}_{t_{1}}$, there exists
$C=C(\alpha,K,t_{0},t_{1})<\infty$ such that for all $\sigma>0$,
if $i=i(\alpha,K,t_{0},t_{1},\sigma)\in\mathbb{N}$ is sufficiently
large, then 
\[
d\nu_{x,t_{1};t_{0}}^{i}\leq Ce^{\alpha f_{i}}d\nu_{t_{0}}^{i}
\]
for all $x\in\psi_{i,t_{1}}(K\cap\{r_{Rm}\geq\sigma\})$.
\end{lem}

\begin{rem}
For the precise definition of $r_{Rm}$ on $\mathbb{F}$-limits as
in (\ref{F}), see Lemma 15.16 and Corollary 15.47 in \cite{bamlergen3}. 
\end{rem}

\begin{proof}
\textbf{Claim: }There exists $D=D(K,t_{1})<\infty$ such that for
any $t^{\ast}<t_{1}$, if $i=i(K,t_{1},t^{\ast},\sigma)\in\mathbb{N}$
sufficiently large, we have 
\[
\sup_{x\in\psi_{i,t}(K\cap\{r_{Rm}\geq\sigma\})}d_{W_{1}}^{g_{i,t^{\ast}}}(\nu_{x,t_{1};t^{\ast}}^{i},\nu_{x_{i},0;t^{\ast}}^{i})\leq D.
\]
If $y_{0}\in\mathcal{X}_{t_{1}}$ is an $H_{2n}$-center of $(\nu_{t})_{t\in(-\infty,0)}$,
then for all $y\in K$, we have 
\[
d_{W_{1}}^{\mathcal{X}_{t_{1}}}(\delta_{y},\nu_{t_{1}})\leq d_{t_{1}}(y,y_{0})+\sqrt{H_{2n}|t_{1}|}\leq C(K,t_{1}).
\]
Let $\mathfrak{C}$ realize the $\mathbb{F}$-convergence as in (\ref{F}).
Fix $\overline{t}\in(t^{\ast},t_{1})$ such that the convergence is
timewise at time $\overline{t}$. Then for all $y\in K\cap\{r_{Rm}\geq\sigma\}$,
and $i=i(K,t_{1},\sigma)\in\mathbb{N}$ sufficiently large, 
\begin{align*}
d_{W_{1}}^{g_{i,t^{\ast}}}(\nu_{\psi_{i,t_{1}}(y),t_{1};t^{\ast}}^{i},\nu_{x_{i},0;t^{\ast}}^{i})\leq & d_{W_{1}}^{g_{i,\overline{t}}}(\nu_{\psi_{i,t_{1}}(y),t_{1};\overline{t}}^{i},\nu_{x_{i},0;\overline{t}}^{i})\\
\leq & d_{W_{1}}^{Z_{\overline{t}}}\left((\varphi_{\overline{t}}^{i})_{\ast}\nu_{\psi_{i,t_{1}}(y),t_{1};t^{\ast}}^{i},(\varphi_{\overline{t}}^{\infty})_{\ast}\nu_{y;t^{\ast}}\right)+d_{W_{1}}^{\mathcal{X}_{\overline{t}}}(\nu_{y;\overline{t}},\nu_{\overline{t}})\\
 & +d_{W_{1}}^{Z_{\overline{t}}}\left((\varphi_{\overline{t}}^{\infty})_{\ast}\nu_{\overline{t}},(\varphi_{\overline{t}}^{i})_{\ast}\nu_{x_{i},0;\overline{t}}^{i}\right)\\
\leq & d_{W_{1}}^{\mathcal{X}_{t_{1}}}(\delta_{y},\nu_{t_{1}})+\Psi(i^{-1}|K,t_{1},t^{\ast},\sigma)
\end{align*}
where we used Theorem 9.31(c) of \cite{bamlergen2}. $\square$

Next, we apply Proposition 8.1 of \cite{bamlergen3}, replacing $t_{0}$
with $0$, $s$ with $t_{0}$, and choosing $t^{\ast}\in[t_{0},t_{1}]$
sufficiently close to $t_{1}$ so that (where $\theta=\theta(1)$
is a universal constant)
\[
t_{1}-t^{\ast}\leq\theta\alpha(t^{\ast}-t_{0}),\,\,\,\,\,\,\,\,\,\,\,\,\,\,\,t^{\ast}-t_{0}\geq\frac{1}{2}(t_{1}-t_{0}).
\]
We replace $x_{0}$ with $x_{i}$, and $x_{1}$ with $\psi_{i,t_{1}}(y)$
for $y\in K\cap\{r_{Rm}\geq\sigma\}$, so that the Claim ensures that
the hypotheses of that proposition are satisfied with $D=D(K,t_{1},t_{0})$
when $i=i(K,t_{1},t_{0},\sigma,\alpha)\in\mathbb{N}$ is sufficiently
large, hence 
\[
d\nu_{\psi_{i,t}(y),t_{1};t_{0}}^{i}\leq C(\alpha,K,t_{0},t_{1})e^{\alpha f_{i}}d\nu_{x_{i},0;t_{0}}.
\]
\end{proof}
Next, we prove the Hörmander $L^{2}$ estimate in the case $\varphi=0$. 
\begin{defn}
\label{correcteq} We say $v\in C^{\infty}(\mathcal{R}_{[t_{0},t_{1}]})$
is the backward drift Laplacian flow with final data $v_{t_{1}}\in C_{c}^{\infty}(\mathcal{R}_{t_{1}})$
if $v$ is given by 
\[
v_{t}(y)=(2\pi|t|)^{n}e^{f_{t}(y)}\int_{\mathcal{R}_{t_{1}}}(2\pi|t_{1}|)^{-n}v_{t_{1}}(x)e^{-f_{t_{1}}(x)}K(x;y)dg_{t_{1}}(x)
\]
for $t\in[t_{0},t_{1}]$ and $y\in\mathcal{R}_{t}$. Equivalently,
$(v_{t}d\nu_{t})_{t\in[t_{0},t_{1}]}$ is the unique conjugate heat
flow on $\mathcal{X}$ (in the sense of Definition 3.13 of \cite{bamlergen2})
with final datum $v_{t_{1}}d\nu_{t_{1}}$. 
\end{defn}

\begin{rem}
This implies that $\partial_{t}v=\Delta_{\overline{\partial},f}v+\mathcal{L}_{\nabla^{1,0}f}v$.
After pulling back by the flow of $\partial_{\mathfrak{t}}$, this
is equivalent to $v$ being the backwards heat flow of the Laplacian
of $(\mathcal{R}_{-1},g_{-1})$ (see the Appendix). 
\end{rem}

\begin{lem}
\label{Hormander} Suppose $-\infty<t_{0}\leq t_{1}<0$, and that
$v\in C^{\infty}(\mathcal{R}_{[t_{0},t_{1}]})$ is the backward drift
Laplacian flow with final datum $v_{t_{1}}\in C_{c}^{\infty}(\mathcal{R}_{t_{1}})$.
Then for any $\eta\in\mathcal{A}_{c}^{0,1}(\mathcal{R}_{t_{0}})$,
we have 
\[
\left|\int_{\mathcal{R}_{t_{0}}}\langle\eta,\overline{\overline{\partial}v_{t_{0}}}\rangle d\nu_{t_{0}}\right|\leq\left(|t_{1}|\int_{\mathcal{R}_{t_{1}}}|\overline{\partial}v_{t_{1}}|^{2}d\nu_{t_{1}}\right)^{\frac{1}{2}}\left(\int_{\mathcal{R}_{t_{0}}}|\overline{\partial}_{f_{t_{0}}}^{\ast}\eta|^{2}d\nu_{t_{0}}\right)^{\frac{1}{2}}.
\]
\end{lem}

\begin{proof}
It will suffice to establish the claim for $t_{0}<t_{1}$, since if
$\theta_{t_{0}}:\mathcal{R}_{t_{0}}\to\mathcal{R}_{t_{1}}$ is obtained
from flowing along the (complete) vector field $\tau(\partial_{\mathfrak{t}}-\nabla f)\in\mathfrak{X}(\mathcal{R})$,
then we can apply the estimate to $\theta_{t_{0}}^{\ast}\eta$, and
take $t_{0}\nearrow t_{1}$. In fact, we know from Lemma \ref{fundsol}$(ii)$
that $(\theta_{t_{0}}^{-1})^{\ast}v_{t_{0}}\to v_{t_{1}}$ pointwise
in $C_{\text{loc}}^{\infty}(\mathcal{R}_{t_{1}})$ (Lemma \ref{fundsol}
is stated only in the case $\mathcal{X}$ is a static cone, which
is what we need for Theorem \ref{thm1}, but the proof can be modified
to deal with the case of a general shrinking soliton by instead considering
the drift Laplacian). Because $(v_{t})$ solves a backwards parabolic
equation with smooth coefficients and because $v_{t_{1}}\in C_{c}^{\infty}(\mathcal{R}_{t_{1}})$,
local parabolic regularity gives higher derivative estimates on $(\theta_{t_{0}}^{-1})^{\ast}v_{t_{0}}$,
uniformly for $t_{0}$ close to $t_{1}$. Thus, the convergence $(\theta_{t_{0}}^{-1})^{\ast}v_{t_{0}}\to v_{t_{1}}$
actually occurs in $C_{\text{loc}}^{\infty}(\mathcal{R}_{t_{1}})$
as $t_{0}\nearrow t_{1}$. 

For sufficiently large $i\in\mathbb{N}$, we can choose $\eta_{i,t_{0}}\in\mathcal{A}_{c}^{0,1}(\psi_{i,t_{0}}(U_{i}\cap\mathcal{R}_{t_{0}}))$
and $v_{i,t_{1}}\in C_{c}^{\infty}(\psi_{i,t_{1}}(U_{i}\cap\mathcal{R}_{t_{1}}))$
such that $\psi_{i,t_{0}}^{\ast}\eta_{i,t_{0}}\to\eta$ in $C_{\text{loc}}^{\infty}(\mathcal{R}_{-1})$
and $\psi_{i,t_{1}}^{\ast}v_{i,t_{1}}\to v_{t_{1}}$ in $C_{\text{loc}}^{\infty}(\mathcal{R}_{t_{1}})$.
Write $\nu_{t}^{i}:=\nu_{x_{i},0;t}^{i}$ for $t\in[-\epsilon_{i}^{-1},0]$.
Let $(\eta_{i,t})_{t\in[t_{0},0]}$ solve 
\[
\partial_{t}\eta_{i,t}=-\Delta_{\overline{\partial}}\eta_{i,t}=\Delta_{g_{i,t}}\eta_{i,t}-Rc(\eta_{i,t})
\]
with initial data $\eta_{i,t_{0}}$, and let $(v_{i,t})_{t\in[t_{0},t_{1}]}$
solve 
\begin{align*}
-\partial_{t}v_{i,t}= & \Delta_{g_{i,t}}v_{i,t}-\langle\nabla v_{i,t},\overline{\nabla}f_{i,t}\rangle-\langle\nabla f_{i,t},\overline{\nabla}v_{i,t}\rangle\\
= & -\Delta_{\overline{\partial},f}v_{i,t}-\nabla_{\nabla^{1,0}f}v_{i,t}
\end{align*}
with final data $v_{i,t_{1}}$. The maximum principle applied to $(\partial_{t}-\Delta)|\eta_{i,t}|\leq0$
and $\partial_{t}v_{i,t}=\Delta_{\overline{\partial},f}v_{i,t}+\nabla_{\nabla^{1,0}f}v_{i,t}$
gives
\[
\sup_{M_{i}}|\eta_{i,t}|\leq\sup_{M_{i}}|\eta_{i,t_{0}}|,
\]
\[
\sup_{M_{i}}|v_{i,t}|\leq\sup_{M_{i}}|v_{i,t_{0}}|
\]
for all $i\in\mathbb{N}$, hence
\[
\limsup_{i\to\infty}\sup_{M_{i}\times[t_{0},t_{1}]}|\eta_{i}|\leq\sup_{\mathcal{R}_{t_{0}}}|\eta|,
\]
\[
\limsup_{i\to\infty}\sup_{M\times[t_{0},t_{1}]}|v_{i}|\leq\sup_{\mathcal{R}_{t_{1}}}|v|.
\]
Next, we write $\eta_{i,t_{0}}=\eta_{i,t_{0}}^{(1)}+\eta_{i,t_{0}}^{(2)}$,
where $\overline{\partial}\eta_{i,t_{0}}^{(1)}=0$ and 
\[
\int_{M_{i}}\langle\eta_{i,t_{0}}^{(2)},\overline{\alpha}\rangle d\nu_{t_{0}}^{i}=0
\]
for all $\alpha\in\ker(\overline{\partial})\subseteq L^{2}(M_{i},d\nu_{t_{0}}^{i})$.
Let $(\eta_{i,t}^{(j)})_{t\in[t_{0},t_{1}]}$ solve
\[
\partial_{t}\eta_{i}^{(j)}=-\Delta_{\overline{\partial}}\eta_{i}^{(j)}
\]
with initial data $\eta_{i,t_{0}}^{(j)}$, so that $\eta_{i,t}=\eta_{i,t}^{(1)}+\eta_{i,t}^{(2)}$
for all $t\in[t_{0},t_{1}]$. Observe that because
\[
\partial_{t}(\overline{\partial}\eta_{i}^{(1)})=-\Delta_{\overline{\partial}}(\overline{\partial}\eta_{i}^{(1)})
\]
and $\overline{\partial}\eta_{i,t_{0}}^{(1)}\equiv0$, the maximum
principle implies $\overline{\partial}\eta_{i,t}^{(1)}\equiv0$ for
all $t\in[t_{0},t_{1}]$. Lemma \ref{veryrough} gives
\[
\int_{M_{i}}|\eta_{i,t}^{(j)}|^{2}d\nu_{t}^{i}\leq\int_{M_{i}}|\eta_{i,t_{0}}^{(j)}|^{2}d\nu_{t_{0}}^{i}\leq\int_{M_{i}}|\eta_{i,t_{0}}|^{2}d\nu_{t_{0}}^{i}\leq C
\]
for all $t\in[t_{0},t_{1}]$, where the last inequality follows from
\[
\lim_{i\to\infty}\int_{M_{i}}|\eta_{i,t_{0}}|^{2}d\nu_{t_{0}}^{i}=\int_{\mathcal{R}_{t_{0}}}|\eta|^{2}d\nu_{t_{0}}<\infty.
\]
Combining this fact with Lemma \ref{gradclose} yields

\begin{align*}
\left|\int_{M_{i}}\langle\eta_{i,t_{0}}^{(1)},\overline{\overline{\partial}v_{i,t_{0}}}\rangle d\nu_{t_{0}}^{i}-\int_{M_{i}}(\overline{\partial}_{h_{i}}^{\ast}\eta_{i,t_{0}}^{(1)})\overline{v_{i,t_{0}}}d\nu_{t_{0}}^{i}\right|\leq & \left|\int_{M_{i}}(\overline{\partial}_{f_{i}}^{\ast}\eta_{i,t_{0}}^{(1)}-\overline{\partial}_{h_{i}}^{\ast}\eta_{i,t_{0}}^{(1)})\overline{v_{i,t_{0}}}d\nu_{t_{0}}^{i}\right|\\
\leq & \int_{M_{i}}|\nabla(f_{i}-h_{i})|\cdot|\eta_{i,t_{0}}^{(1)}|\cdot|v_{i,t_{0}}|d\nu_{t_{0}}^{i}\\
\leq & \left(\sup_{M_{i}}|v_{i,t_{1}}|\right)\left(\int_{M_{i}}|\nabla(f_{i}-h_{i})|^{2}d\nu_{t_{0}}^{i}\right)^{\frac{1}{2}}\left(\int_{M_{i}}|\eta_{i,t_{0}}^{(1)}|^{2}d\nu_{t_{0}}^{i}\right)^{\frac{1}{2}}\\
\leq & \Psi(i^{-1}),
\end{align*}
and similarly
\[
\left|\int_{M_{i}}\langle\eta_{i,t_{1}}^{(1)},\overline{\overline{\partial}v_{i,t_{1}}}\rangle d\nu_{t_{1}}^{i}-\int_{M_{i}}(\overline{\partial}_{h_{i}}^{\ast}\eta_{i,t_{1}}^{(1)})\overline{v_{i,t_{1}}}d\nu_{t_{1}}^{i}\right|\leq\Psi(i^{-1}).
\]
Recall from Lemma \ref{slowchange} that
\begin{align*}
\frac{d}{dt}\int_{M_{i}}\tau(\overline{\partial}_{h_{i}}^{\ast}\eta_{i,t}^{(1)})\overline{v_{i,t}}e^{-f_{i,t}}dg_{i,t}= & \int_{M_{i}}\tau\left\langle Rc_{g_{i,t}}+\nabla\overline{\nabla}h_{g_{i,t}}-\frac{1}{\tau}g_{i,t},\overline{\nabla}\eta_{i,t}^{(1)}\right\rangle \overline{v_{i,t}}d\nu_{t}^{i}+\int_{M_{i}}\tau\langle\nabla\nabla h_{i,t},\overline{\nabla\eta_{i,t}^{(1)}}\rangle\overline{v_{i,t}}d\nu_{t}^{i},
\end{align*}
which along with Lemma \ref{veryrough} allows us to estimate
\begin{align*}
\left||t_{0}|\int_{M_{i}}\langle\eta_{i,t_{0}}^{(1)},\overline{\overline{\partial}v_{i,t_{0}}}\rangle d\nu_{t_{0}}^{i}-|t_{1}|\int_{M_{i}}\langle\eta_{i,t_{1}}^{(1)},\overline{\overline{\partial}v_{i,t_{1}}}\rangle d\nu_{t_{1}}^{i}\right|\leq & \Psi(i^{-1})\left(\int_{t_{0}}^{t_{1}}\int_{M_{i}}|\nabla^{\mathbb{R}}\eta_{i,t}^{(1)}|^{2}d\nu_{t}^{i}dt\right)^{\frac{1}{2}}=\Psi(i^{-1}).
\end{align*}
Combining expressions yields
\[
\left||t_{0}|\int_{M_{i}}\langle\eta_{i,t_{0}},\overline{\overline{\partial}v_{i,t_{0}}}\rangle d\nu_{t_{0}}^{i}-|t_{1}|\int_{M_{i}}\langle\eta_{i,t_{1}}^{(1)},\overline{\overline{\partial}v_{i,t_{1}}}\rangle d\nu_{t_{1}}^{i}\right|\leq\Psi(i^{-1}).
\]
Now fix $x_{0}\in\mathcal{X}_{t_{1}}$, and define $\rho_{r}(x):=\rho(r^{-1}d_{t_{1}}(x,x_{0}))$
for $x\in\mathcal{X}_{t_{1}}$, so that $|\nabla\rho_{r}|\leq1$ on
$\mathcal{R}_{t_{1}}$. Because $\mathcal{S}_{t_{1}}\subseteq\mathcal{X}_{t}$
has singularities of Minkowski codimension $4$, we can find $(\chi_{\epsilon})_{\epsilon\in\mathbb{N}}$
in $C_{c}^{\infty}(\mathcal{R}_{t_{1}})$ such that 
\[
\chi_{\epsilon}|_{\{r_{Rm}>2\epsilon\}}\equiv1,\,\,\,\,\,\,\,\,\,\,\,\,\,\,\,\,\,\,\,\,\text{supp}(\chi_{\epsilon})\subseteq\{r_{Rm}>\epsilon\},
\]
and also
\[
\lim_{\epsilon\searrow0}\int_{\mathcal{R}_{t_{1}}\cap\text{supp}(\rho_{r})}|\nabla\chi_{\epsilon}|^{3}dg_{t_{1}}=0
\]
for any $r\in(0,\infty)$. For each $r<\infty$ and $\epsilon>0$,
we define $\chi_{r,\epsilon,i}:=(\rho_{r}\chi_{\epsilon})\circ\psi_{i,t_{1}}^{-1}$
for $i=i(\epsilon,r)\in\mathbb{N}$ sufficiently large. The Kodaira-Nakano
formula then yields
\[
(\overline{\partial}\overline{\partial}_{f_{i}}^{\ast}+\overline{\partial}_{f_{i}}^{\ast}\overline{\partial})(\chi_{r,\epsilon,i}\eta_{i,t_{1}}^{(1)})=(\overline{\nabla}^{f_{i}})^{\ast}\overline{\nabla}(\chi_{r,\epsilon,i}\eta_{i,t_{1}}^{(1)})+\left(Rc_{g_{i,t_{1}}}+\nabla\overline{\nabla}f_{i,t_{1}}\right)(\cdot,\overline{\chi_{r,\epsilon,i}\eta_{i,t_{1}}^{(1)}}),
\]
which we can integrate against $\chi_{r,\epsilon,i}\overline{\eta_{i,t_{1}}^{(1)}}$
to obtain
\begin{align*}
\int_{M_{i}}|\overline{\partial}_{f_{i}}^{\ast}(\chi_{r,\epsilon,i}\eta_{i,t_{1}}^{(1)})|^{2}d\nu_{t_{1}}^{i} & +\int_{M_{i}}|\overline{\partial}(\chi_{r,\epsilon,i}\eta_{i,t_{1}}^{(1)})|^{2}d\nu_{t_{1}}^{i}\\
= & \int_{M_{i}}|\overline{\nabla}(\chi_{r,\epsilon,i}\eta_{i,t_{1}}^{(1)})|^{2}d\nu_{t_{1}}^{i}\\
 & +\int_{M_{i}}\left(Rc_{g_{i,t_{1}}}+\nabla\overline{\nabla}f_{i,t_{1}}\right)(\chi_{r,\epsilon,i}\eta_{i,t_{1}}^{(1)},\overline{\chi_{r,\epsilon,i}\eta_{i,t_{1}}^{(1)}})d\nu_{t_{1}}^{i}.
\end{align*}
We next apply Hölder's inequality and Lemma \ref{veryrough} to obtain
(for some $p=p(t_{0},t_{1})>2$)
\begin{align*}
\int_{M_{i}}\left(Rc_{g_{i,t_{1}}}+\nabla\overline{\nabla}f_{i,t_{1}}\right) & (\chi_{r,\epsilon,i}\eta_{i,t_{1}}^{(1)},\overline{\chi_{r,\epsilon,i}\eta_{i,t_{1}}^{(1)}})d\nu_{t_{1}}^{i}\\
\geq & \frac{1}{|t_{1}|}\int_{M_{i}}|\eta_{i,t_{1}}^{(1)}|^{2}d\nu_{t_{1}}^{i}-\frac{1}{|t_{1}|}\int_{M_{i}}(1-\chi_{r,\epsilon,i}^{2})|\eta_{i,t_{1}}^{(1)}|^{2}d\nu_{t_{1}}^{i}\\
 & -\left(\sup_{\text{supp}(\chi_{r,\epsilon,i})}\left|Rc_{g_{i,t_{1}}}+\nabla\overline{\nabla}f_{i,t_{1}}-\frac{1}{|t_{1}|}g_{i,t_{1}}\right|\right)\int_{M_{i}}|\eta_{i,t_{1}}^{(1)}|^{2}d\nu_{t_{1}}^{i}\\
\geq & \left(\frac{1}{|t_{1}|}-\Psi(i^{-1}|\epsilon,r)\right)\int_{M_{i}}|\eta_{i,t_{1}}^{(1)}|^{2}d\nu_{t_{1}}^{i}\\
 & -C(t_{0},t_{1})\left(\int_{M}|\eta_{i,t_{1}}^{(1)}|^{p}d\nu_{t_{1}}\right)^{\frac{2}{p}}\left(\nu_{t_{1}}^{i}(\text{supp}(1-\chi_{r,\epsilon,i}))\right)^{\frac{2(p-1)}{p}}\\
\geq & \left(\frac{1}{|t_{1}|}-\Psi(i^{-1}|\epsilon,r)-\Psi(\epsilon|r)-\Psi(r^{-1})\right)\int_{M_{i}}|\eta_{i,t_{1}}^{(1)}|^{2}d\nu_{t_{1}}^{i}.
\end{align*}
Combining expressions yields (for $i=i(r,\epsilon)\in\mathbb{N}$
sufficiently large)
\begin{align*}
\left|\int_{M_{i}}\langle\eta_{i,t_{1}},\overline{\overline{\partial}v_{i,t_{1}}}\rangle d\nu_{t_{1}}^{i}\right|\leq & \left(\int_{M_{i}}|\eta_{i,t_{1}}^{(1)}|^{2}d\nu_{t_{1}}^{i}\right)^{\frac{1}{2}}\left(\int_{M_{i}}|\overline{\partial}v_{i,t_{1}}|^{2}d\nu_{t_{1}}^{i}\right)^{\frac{1}{2}}\\
\leq & \left(|t_{1}|\int_{M_{i}}(Rc_{g_{i,t_{1}}}+\nabla\overline{\nabla}f_{i,t_{1}})(\overline{\chi_{r,\epsilon,i}\eta_{i,t_{1}}^{(1)}},\chi_{r,\epsilon,i}\eta_{i,t_{1}}^{(1)})d\nu_{t_{1}}^{i}\right)^{\frac{1}{2}}\\
 & \times\left(\int_{M_{i}}|\overline{\partial}v_{i,t_{1}}|^{2}d\nu_{t_{1}}^{i}\right)^{\frac{1}{2}}+\Psi(i^{-1}|\epsilon,r)+\Psi(\epsilon|r)+\Psi(r^{-1})\\
\leq & \left(\int_{M_{i}}|\overline{\partial}\chi_{r,\epsilon,i}\wedge\eta_{i,t_{1}}^{(1)}|^{2}d\nu_{t_{1}}^{i}+\int_{M_{i}}|\overline{\partial}_{f_{i}}^{\ast}(\chi_{r,\epsilon,i}\eta_{i,t_{1}}^{(1)})|^{2}d\nu_{t_{1}}^{i}\right)^{\frac{1}{2}}\\
 & \times\left(|t_{1}|\int_{M_{i}}|\overline{\partial}v_{i,t_{1}}|^{2}d\nu_{t_{1}}^{i}\right)^{\frac{1}{2}}++\Psi(i^{-1}|\epsilon,r)+\Psi(\epsilon|r)+\Psi(r^{-1}).
\end{align*}
For any $x\in\text{supp}(\nabla\chi_{r,\epsilon,i})$, We can apply
Lemma \ref{nearby} with $\alpha=\frac{1}{4}$, $K$ replaced by $B(x_{0},2r)$,
and $\sigma=\epsilon$ to obtain
\[
d\nu_{x,t_{1};t_{0}}^{i}\leq C(r,t_{1},t_{0})e^{\frac{1}{4}f_{i}}d\nu_{t_{0}}^{i}
\]
for all $x\in\text{supp}(\nabla\chi_{r,\epsilon,i})$ when $i=i(r,\epsilon,t_{1},t_{0})\in\mathbb{N}$
is sufficiently large. Because $|\eta_{i,t_{1}}|$ is a subsolution
of the heat equation, we therefore have
\begin{align*}
|\eta_{i,t_{1}}^{(1)}|(x)\leq & \int_{M_{i}}|\eta_{i,t_{0}}^{(1)}|d\nu_{x,t_{1};t_{0}}^{i}\leq C(r,t_{1},t_{0})\int_{M_{i}}|\eta_{i,t_{0}}^{(1)}|e^{\frac{1}{4}f_{i}}d\nu_{t_{0}}^{i}\\
\leq & C(r,t_{1},t_{0})\left(\int_{M_{i}}|\eta_{i,t_{0}}^{(1)}|^{2}d\nu_{t_{0}}^{i}\right)^{\frac{1}{2}}\left(\int_{M_{i}}e^{\frac{1}{2}f_{i}}d\nu_{t_{0}}^{i}\right)^{\frac{1}{2}}\\
\leq & C(r,t_{1},t_{0})
\end{align*}
for all $x\in\text{supp}(\nabla\chi_{r,\epsilon,i})$ when $i=i(r,\epsilon,t_{1},t_{0})\in\mathbb{N}$
is sufficiently large. We can therefore estimate (suppressing the
dependence on $t_{0},t_{1}$)
\begin{align*}
\int_{M_{i}}|\overline{\partial}\chi_{r,\epsilon,i}\wedge\eta_{i,t_{1}}^{(1)}|^{2}d\nu_{t_{1}}^{i}\leq & C(r)\int_{M_{i}\cap\psi_{i,t_{1}}(\text{supp}(\rho_{r}))}|\nabla\chi_{\epsilon,i}|^{2}d\nu_{t_{1}}^{i}+Cr^{-2}\int_{M_{i}}|\eta_{i,t_{1}}^{(1)}|^{2}d\nu_{t_{1}}^{i}\\
\leq & \Psi(i^{-1}|\epsilon,r)+\Psi(\epsilon|r)+\Psi(r^{-1}),
\end{align*}
and similarly 
\begin{align*}
\int_{M_{i}}|\overline{\partial}_{f_{i}}^{\ast}(\chi_{r,\epsilon,i}\eta_{i,t_{1}}^{(1)})|^{2}d\nu_{t_{1}}^{i}\leq & (1+\Psi(i^{-1}|\epsilon,r)+\Psi(\epsilon|r)+\Psi(r^{-1}))\int_{\text{supp}(\chi_{r,\epsilon,i})}|\overline{\partial}_{f_{i}}^{\ast}\eta_{i,t_{1}}^{(1)}|^{2}d\nu_{t_{1}}^{i}\\
 & +\Psi(i^{-1}|\epsilon,r)+\Psi(\epsilon|r)+\Psi(r^{-1})
\end{align*}
so that (since $\psi_{i,t_{1}}^{\ast}(h_{i,t_{1}}-f_{i,t_{1}})\to0$
in $C_{\text{loc}}^{\infty}(\mathcal{R}_{t_{1}})$)
\begin{align*}
\left|\int_{M_{i}}\langle\eta_{i,t_{1}},\overline{\overline{\partial}v_{i,t_{1}}}\rangle d\nu_{t_{1}}^{i}\right|\leq & \left(|t_{1}|\int_{M_{i}}|\overline{\partial}v_{i,t_{1}}|^{2}d\nu_{t_{1}}^{i}\right)^{\frac{1}{2}}\left(\int_{\text{supp}(\chi_{r,\epsilon,i})}|\overline{\partial}_{h_{i}}^{\ast}\eta_{i,t_{1}}^{(1)}|^{2}d\nu_{t_{1}}^{i}\right)^{\frac{1}{2}}\\
 & +\Psi(i^{-1}|\epsilon,r)+\Psi(\epsilon|r)+\Psi(r^{-1}).
\end{align*}
By Corollary \ref{Essential}, for any $a>0$ we have
\begin{equation}
(\partial_{t}-\Delta)\left(-aw_{i}+\frac{1}{a}|\eta_{i}^{(2)}|^{2}+\tau|\overline{\partial}_{h_{i}}^{\ast}\eta_{i}^{(2)}|\right)\leq0,\label{nice}
\end{equation}
so that 
\[
\tau|\overline{\partial}_{h_{i}}^{\ast}\eta_{i,t}^{(2)}(x)|\leq aw_{i,t}(x)+\left|\int_{M_{i}}\left(-aw_{i,t_{0}}+\frac{1}{a}|\eta_{i,t_{0}}^{(2)}|^{2}+\tau|\overline{\partial}_{h_{i}}^{\ast}\eta_{i,t_{0}}^{(2)}|\right)d\nu_{x,t;t_{0}}^{i}\right|
\]
hence integration against $d\nu_{t}^{i}$ yields
\begin{align*}
\tau\int_{M_{i}}|\overline{\partial}_{h_{i}}^{\ast}\eta_{i,t}^{(2)}|d\nu_{t}^{i}\leq & \frac{2}{a}\int_{M_{i}}|\eta_{i,t_{0}}^{(2)}|^{2}d\nu_{x,t;t_{0}}^{i}+\Psi(i^{-1}|a)\leq\frac{2}{a}\int_{M_{i}}|\eta_{i,t_{0}}|^{2}d\nu_{x,t;t_{0}}^{i}+\Psi(i^{-1}|a)\\
\leq & \Psi(a)+\Psi(i^{-1}|a)
\end{align*}
for all $t\in[t_{0},t_{1}]$. Choosing $a>0$ small gives 
\[
\int_{M_{i}}|\overline{\partial}_{h_{i}}^{\ast}\eta_{i,t}^{(2)}|d\nu_{t_{1}}^{i}\leq\Psi(i^{-1}).
\]
Moreover, because $h_{i},\eta_{i}^{(2)}$ satisfy heat-type equations,
and are uniformly (in $i\in\mathbb{N}$) bounded on $\text{supp}(\chi_{2r,\frac{\epsilon}{2},i})$
(for $r<\infty$, $\epsilon>0$ fixed), local parabolic regularity
gives 
\[
\sup_{\text{supp}(\chi_{r,\epsilon,i})}|\overline{\partial}_{h_{i}}^{\ast}\eta_{i,t_{1}}^{(2)}|\leq C(\epsilon,r).
\]
Combining these facts gives
\[
\int_{\text{supp}(\chi_{r,\epsilon,i})}|\overline{\partial}_{h_{i}}^{\ast}\eta_{i,t}^{(2)}|^{2}d\nu_{t_{1}}^{i}\leq\Psi(i^{-1}|\epsilon,r).
\]
Because $\eta_{i,t}=\eta_{i,t}^{(1)}+\eta_{i,t}^{(2)}$, 
\[
\left||t_{0}|\int_{M_{i}}\langle\eta_{i,t_{0}},\overline{\overline{\partial}v_{i,t_{0}}}\rangle d\nu_{t_{0}}^{i}\right|\leq|t_{1}|\left(\int_{\text{supp}(\chi_{r,\epsilon,i})}|\overline{\partial}_{h_{i}}^{\ast}\eta_{i,t_{1}}|^{2}d\nu_{t_{1}}^{i}\right)^{\frac{1}{2}}\left(|t_{1}|\int_{M_{i}}|\overline{\partial}v_{i,t_{1}}|^{2}d\nu_{t_{1}}^{i}\right)^{\frac{1}{2}}+\Psi(i^{-1}|\epsilon,r)+\Psi(\epsilon|r)+\Psi(r^{-1})
\]
For $\sigma\in(0,1)$ to be determined, (\ref{nice}) yields
\begin{align*}
\tau^{2}|\overline{\partial}_{h_{i}}^{\ast}\eta_{i,t_{1}}(x)|^{2}\leq & C\frac{a^{2}}{\sigma}w_{i,t_{1}}^{2}(x)+(1+\sigma)\left(\int_{M_{i}}\tau|\overline{\partial}_{h_{i}}^{\ast}\eta_{i,t_{0}}|d\nu_{x,t_{1};t_{0}}^{i}\right)^{2}\\
 & +C\frac{a^{2}}{\sigma}\left(\int_{M_{i}}|w_{i,t_{0}}|d\nu_{x,t_{1};t_{0}}^{i}\right)^{2}+\frac{C(n)}{a^{2}\sigma}\left(\int_{M_{i}}|\eta_{i,t_{0}}|^{2}d\nu_{x,t_{1};t_{0}}^{i}\right)^{2},
\end{align*}
so we can integrate against $d\nu_{t_{1}}^{i}$ to obtain

\begin{align*}
|t_{1}|^{2}\int_{\text{supp}(\chi_{r,\epsilon,i})}|\overline{\partial}_{h_{i}}^{\ast}\eta_{i,t_{1}}|^{2}d\nu_{t_{1}}^{i}\leq & C\frac{a^{2}}{\sigma}\sup_{\text{supp}(\chi_{r,\epsilon,i})}w_{i,t_{1}}^{2}+C\frac{a^{2}}{\sigma}\sup_{x\in\text{supp}(\chi_{r,\epsilon,i})}\left(\int_{M_{i}}|w_{i,t_{0}}|d\nu_{x,t_{1};t_{0}}^{i}\right)^{2}\\
 & +(1+\sigma)\int_{M_{i}}\left(\int_{M_{i}}|\tau\overline{\partial}_{h_{i}}^{\ast}\eta_{i,t_{0}}|^{2}d\nu_{x,t_{1};t_{0}}^{i}\right)d\nu_{t_{1}}^{i}(x)+\frac{C(n)}{a^{2}\sigma}\sup_{M_{i}}|\eta_{i,t_{0}}|^{4}\\
\leq & \frac{a^{2}}{\sigma}\Psi(i^{-1}|\epsilon,r)+\frac{C(n)}{a^{2}\sigma}\sup_{M_{i}}|\eta_{i,t_{0}}|^{4}+(1+\sigma)|t_{0}|^{2}\int_{M_{i}}|\overline{\partial}_{h_{i}}^{\ast}\eta_{i,t_{0}}|^{2}d\nu_{t_{0}}^{i}.
\end{align*}
Taking $\sigma>0$ small, then $a=a(\sigma)>0$ large gives
\[
|t_{1}|^{2}\int_{\text{supp}(\chi_{r,\epsilon,i})}|\overline{\partial}_{h_{i}}^{\ast}\eta_{i,t_{1}}|^{2}d\nu_{t_{1}}^{i}\leq\Psi(i^{-1}|\epsilon,r)+(1+\Psi(i^{-1}|\epsilon,r))|t_{0}|^{2}\int_{M_{i}}|\overline{\partial}_{h_{i}}^{\ast}\eta_{i,t_{0}}|^{2}d\nu_{t_{0}}^{i}.
\]
Again combining estimates, we obtain
\begin{align*}
\left||t_{0}|\int_{M_{i}}(\overline{\partial}_{h_{i}}^{\ast}\eta_{i,t_{0}})\overline{v}_{i,t_{0}}d\nu_{t_{0}}^{i}\right|= & \left|\int_{M_{i}}|t_{0}|\langle\eta_{i,t_{0}},\overline{\overline{\partial}v_{i,t_{0}}}\rangle d\nu_{t_{0}}^{i}\right|+\Psi(i^{-1})\\
\leq & \left(\Psi(i^{-1}|\epsilon,r)+(1+\Psi(i^{-1}|\epsilon,r))|t_{0}|^{2}\int_{M_{i}}|\overline{\partial}_{h_{i}}^{\ast}\eta_{i,t_{0}}|^{2}d\nu_{t_{0}}^{i}\right)^{\frac{1}{2}}\left(|t_{1}|\int_{M_{i}}|\overline{\partial}v_{i,t_{1}}|^{2}d\nu_{t_{1}}^{i}\right)^{\frac{1}{2}}\\
 & +\Psi(i^{-1}|\epsilon,r)+\Psi(\epsilon|r)+\Psi(r^{-1}).
\end{align*}
Choosing $r<\infty$ large, then $\epsilon>0$ small, we obtain
\begin{align*}
\left|\int_{M_{i}}(\overline{\partial}_{h_{i}}^{\ast}\eta_{i,t_{0}})\overline{v}_{i,t_{0}}d\nu_{t_{0}}^{i}\right|\leq & (1+\Psi(i^{-1}))\left(\int_{M_{i}}|\overline{\partial}_{h_{i}}^{\ast}\eta_{i,t_{0}}|^{2}d\nu_{t_{0}}^{i}\right)^{\frac{1}{2}}\left(|t_{1}|\int_{M_{i}}|\overline{\partial}v_{i,t_{1}}|^{2}d\nu_{t_{1}}^{i}\right)^{\frac{1}{2}}+\Psi(i^{-1}).
\end{align*}
Taking $i\to\infty$ then gives 
\[
\left|\int_{\mathcal{R}_{t_{0}}}\langle\eta,\overline{\overline{\partial}v_{t_{0}}}\rangle d\nu_{t_{0}}\right|=\left|\int_{\mathcal{R}_{t_{0}}}(\overline{\partial}_{h}^{\ast}\eta)\overline{v}_{t_{0}}d\nu_{t_{0}}\right|\leq\left(\int_{\mathcal{R}_{t_{0}}}|\overline{\partial}_{f}^{\ast}\eta|^{2}d\nu_{t_{0}}\right)^{\frac{1}{2}}\left(|t_{1}|\int_{\mathcal{R}_{t_{1}}}|\overline{\partial}v|^{2}d\nu_{t_{1}}\right)^{\frac{1}{2}},
\]
where the convergence $\psi_{i,t_{0}}^{\ast}v_{i,t_{0}}\to v_{t_{0}}$
in $C_{\text{loc}}^{0}(\mathcal{R}_{t_{0}})$ can be justified using
Remark \ref{convergenotation}, Definition \ref{correcteq}, and the
fact that

\[
\psi_{i,t_{1}}^{\ast}\left((2\pi|t_{1}|)^{-n}v_{i,t_{1}}e^{-f_{i,t_{1}}}\right)\to(2\pi|t_{1}|)^{-n}ve^{-f_{t_{1}}}
\]
in $C_{c}^{\infty}(\mathcal{R}_{t_{1}})$. 
\end{proof}
Lemma \ref{Hormander} immediately implies the following existence
theorem for solutions of the $\overline{\partial}$-equation. 
\begin{cor}
\label{hormandersection} With the same notation as the previous lemma,
there exists $u_{t_{0}}\in L^{2}(\mathcal{R}_{t_{0}},d\nu_{t_{0}})$
satisfying $\overline{\partial}u_{t_{0}}=\overline{\partial}v_{t_{0}}$
and 
\[
\int_{\mathcal{R}_{t_{0}}}|u|^{2}d\nu_{t_{0}}\leq|t_{1}|\int_{\mathcal{R}_{t_{1}}}|\overline{\partial}v_{t_{1}}|^{2}d\nu_{t_{1}}.
\]
\end{cor}

\begin{proof}
The linear functional $\overline{\partial}_{f_{t_{0}}}^{\ast}\left(\mathcal{A}_{c}^{0,1}(\mathcal{R}_{t_{0}})\right)\to\mathbb{C}$
given by
\[
\overline{\partial}_{f_{t_{0}}}^{\ast}\eta\mapsto\int_{\mathcal{R}_{t_{0}}}\langle\eta,\overline{\overline{\partial}v_{t_{0}}}\rangle d\nu_{t_{0}}
\]
is well-defined and bounded by the previous lemma, so the Hahn-Banach
theorem and Riesz isomorphism theorem produce $u\in L^{2}(\mathcal{R}_{t_{0}},d\nu_{t_{0}})$
satisfying
\[
\int_{\mathcal{R}_{t_{0}}}\langle\eta,\overline{\overline{\partial}u}\rangle d\nu_{t_{0}}=\int_{\mathcal{R}_{t_{0}}}(\overline{\partial}_{f_{t_{0}}}^{\ast}\eta)\overline{u}d\nu_{t_{0}}=\int_{\mathcal{R}_{t_{0}}}\langle\eta,\overline{\overline{\partial}v_{t_{0}}}\rangle d\nu_{t_{0}}
\]
for all $\eta\in\mathcal{A}_{c}^{0,1}(\mathcal{R}_{t_{0}})$ and 
\[
\int_{\mathcal{R}_{t_{0}}}|u|^{2}d\nu_{t_{0}}\leq|t_{1}|\int_{\mathcal{R}_{t_{1}}}|\overline{\partial}v_{t_{1}}|^{2}d\nu_{t_{1}}.
\]
In particular, we have $\overline{\partial}u=\overline{\partial}v_{t_{0}}$
in $\mathcal{R}_{t_{0}}$. 
\end{proof}
We now prove Hörmander-type estimates where we allow the weight to
be altered in a compact region of the regular set. This will be used
to prove that we can produce holomorphic functions on $\mathcal{X}$
which separate tangents at any point in $\mathcal{R}$. The proof
is very similar to the case $\varphi\equiv0$, so we mostly indicate
the necessary changes. 
\begin{lem}
\label{HormanderwWeights} Suppose $t_{1}\in(-\infty,0)$ and $v\in C_{c}^{\infty}(\mathcal{R}_{t_{1}},\mathbb{C})$,
$\varphi\in C_{c}^{\infty}(\mathcal{R}_{t_{1}})$ satisfy
\[
\sqrt{-1}\partial\overline{\partial}(f_{t_{1}}+\varphi)\geq\beta\omega_{t_{1}}
\]
for some $\beta>0$. Then for any $\eta\in\mathcal{A}_{c}^{0,1}(\mathcal{R}_{t_{1}})$,
we have 
\[
\left|\int_{\mathcal{R}_{t_{1}}}\langle\eta,\overline{\overline{\partial}v_{t_{1}}}\rangle e^{-\varphi}d\nu_{t_{1}}\right|\leq\left(\frac{1}{\beta}\int_{\mathcal{R}_{t_{1}}}|\overline{\partial}v_{t_{1}}|^{2}e^{-\varphi}d\nu_{t_{1}}\right)^{\frac{1}{2}}\left(\int_{\mathcal{R}_{t_{1}}}|\overline{\partial}_{f_{t_{1}}+\varphi}^{\ast}\eta|^{2}e^{-\varphi}d\nu_{t_{1}}\right)^{\frac{1}{2}}.
\]
\end{lem}

\begin{rem}
Note that in this case, we can only solve the $\overline{\partial}$-equation
at time $t_{1}$, unlike Lemma \ref{Hormander}.
\end{rem}

\begin{proof}
Let $\theta_{t_{0}}:\mathcal{R}_{t_{0}}\to\mathcal{R}_{t_{1}}$ be
as in Lemma \ref{Hormander}, and set $\eta_{t_{0}}:=\theta_{t_{0}}^{\ast}\eta$,
$\varphi_{t}:=\theta_{t}^{\ast}\varphi$ for $t\in(-\infty,t_{0})$.
We correspondingly let $\eta_{i,t_{0}}\in\mathcal{A}_{c}^{0,1}(M_{i})$
and $v_{i,t_{1}}\in C_{c}^{\infty}(M_{i},\mathbb{C})$ be as in Lemma
\ref{Hormander}, and set $\varphi_{i,t}:=(\psi_{i,t}^{-1})^{\ast}\varphi_{t}$
(extended to zero outside the domain of $\psi_{i,t}^{-1}$). We let
$(v_{i,t})_{t\in[t_{0},t_{1}]}$ be the solution of
\[
\partial_{t}v_{i,t}=\Delta_{\overline{\partial},f_{i}+\varphi_{i}}v_{i,t}+\mathcal{L}_{\nabla^{1,0}(f_{i}+\varphi_{i})}v_{i,t}
\]
with final value $v_{i,t_{1}}$, and we write $\eta_{i,t_{0}}=\eta_{i,t_{0}}^{(1)}+\eta_{i,t_{0}}^{(2)}$,
where $\overline{\partial}\eta_{i,t_{0}}^{(1)}\equiv0$ and 
\[
\int_{M_{i}}\langle\eta_{i,t_{0}}^{(2)},\overline{\alpha}\rangle e^{-\varphi_{i,t_{0}}}d\nu_{t_{0}}^{i}=0
\]
for all $\alpha\in\ker(\overline{\partial})\subseteq L^{2}(M_{i},e^{-\varphi_{i,t_{0}}}d\nu_{t_{0}}^{i})$.
We then let $(\eta_{i,t}^{(j)})_{t\in[t_{0},t_{1}]}$ solve 
\[
\partial_{t}\eta_{i,t}^{(j)}=-\Delta_{\overline{\partial}}\eta_{i,t}^{(j)}
\]
with initial condition $\eta_{i,t_{0}}^{(j)}$. Arguing as in Lemma
\ref{Hormander}, we can estimate
\[
\left|\int_{M_{i}}(\overline{\partial}_{h_{i,t_{k}}+\varphi_{i,t_{k}}}^{\ast}\eta_{i,t_{k}}^{(1)})\overline{v_{i,t_{k}}}e^{-\varphi_{i,t_{k}}}d\nu_{t_{k}}^{i}-\int_{M_{i}}\langle\eta_{i,t_{k}}^{(1)},\overline{\overline{\partial}v_{i,t_{k}}}\rangle e^{-\varphi_{i,t_{k}}}d\nu_{t_{k}}^{i}\right|\leq\Psi(i^{-1})
\]
for $k=0,1$. Integrating Lemma \ref{slowchange} in time then leads
to
\begin{align*}
\left||t_{0}|\int_{M}(\overline{\partial}_{h_{t_{0}}+\varphi_{t_{0}}}^{\ast}\eta_{i,t_{0}}^{(1)})\overline{v_{t_{0}}}\right. & \left.e^{-\varphi_{i,t_{0}}}d\nu_{t_{0}}^{i}-|t_{1}|\int_{M}(\overline{\partial}_{h_{t_{1}}+\varphi_{t_{1}}}^{\ast}\eta_{i,t_{1}}^{(1)})\overline{v_{t_{1}}}e^{-\varphi_{i,t_{1}}}d\nu_{t_{1}}^{i}\right|\\
\leq & \left(\sup_{M_{i}\times[t_{0},t_{1}]}|v_{i}|e^{-\varphi_{i}}\right)\left(\int_{t_{0}}^{t_{1}}\int_{M_{i}}|Rc+\nabla\overline{\nabla}h_{i}-\frac{1}{\tau}g_{i}|^{2}d\nu_{t}^{i}dt\right)^{\frac{1}{2}}\left(\int_{t_{0}}^{t_{1}}\int_{M_{i}}|\overline{\nabla}\eta_{i,t}^{(1)}|^{2}d\nu_{t}^{i}dt\right)^{\frac{1}{2}}\\
 & +\left(\sup_{M_{i}\times[t_{0},t_{1}]}|v_{i}|e^{-\varphi_{i}}\right)\left(\int_{t_{0}}^{t_{1}}\int_{M_{i}}|\nabla\nabla h_{i}|^{2}d\nu_{t}^{i}dt\right)^{\frac{1}{2}}\left(\int_{t_{0}}^{t_{1}}\int_{M_{i}}|\nabla\eta_{i,t}^{(1)}|^{2}d\nu_{t}^{i}dt\right)^{\frac{1}{2}}\\
 & +\left(\sup_{M_{i}\times[t_{0},t_{1}]}(|\nabla\overline{\nabla}\varphi_{i}|+|\nabla\nabla\varphi_{i}|)|v_{i}|e^{-\varphi}\right)\int_{t_{0}}^{t_{1}}\int_{M_{i}}(|\nabla\eta_{i,t}^{(1)}|+|\overline{\nabla}\eta_{i,t}^{(1)}|)d\nu_{t}^{i}dt\\
 & +\left(\sup_{M_{i}\times[t_{0},t_{1}]}\left|\nabla(\partial_{t}-\Delta)\varphi_{i}-\frac{1}{\tau}\varphi_{i}\right||v_{i}|e^{-\varphi_{i}}\right)\int_{t_{0}}^{t_{1}}\int_{M_{i}}|\eta_{i,t}^{(1)}|d\nu_{t}^{i}dt\\
\leq & \Psi(i^{-1})+C\sqrt{|t_{1}-t_{0}|}\left(\int_{t_{0}}^{t_{1}}\int_{M_{i}}(|\nabla\eta_{i,t}^{(1)}|^{2}+|\overline{\nabla}\eta_{i,t}^{(1)}|^{2}+|\eta_{i,t}^{(1)}|^{2})d\nu_{t}^{i}dt\right)^{\frac{1}{2}}\\
\leq & \Psi(i^{-1}|t_{1},t_{0})+\Psi(|t_{1}-t_{0}|).
\end{align*}
Combining estimates leads to 
\[
\left||t_{1}|\int_{M_{i}}\langle\eta_{i,t_{1}}^{(1)},\overline{\overline{\partial}v_{i}\rangle}e^{-\varphi_{i,t_{1}}}d\nu_{t_{1}}^{i}-|t_{0}|\int_{M_{i}}\langle\eta_{i},\overline{\overline{\partial}v_{i,t_{0}}}\rangle e^{-\varphi_{i,t_{0}}}d\nu_{t_{0}}^{i}\right|\leq\Psi(i^{-1}|t_{1},t_{0})+\Psi(|t_{1}-t_{0}|).
\]
Next, we estimate
\[
\left|\int_{M_{i}}\langle\eta_{i,t_{1}}^{(1)},\overline{\overline{\partial}v_{i}\rangle}e^{-\varphi_{i,t_{1}}}d\nu_{t_{1}}^{i}\right|\leq\left(\int_{M_{i}}\beta|\eta_{i,t_{1}}^{(1)}|^{2}e^{-\varphi_{i,t_{1}}}d\nu_{t_{1}}^{i}\right)^{\frac{1}{2}}\left(\int_{M_{i}}\frac{1}{\beta}|\overline{\partial}v_{i}|^{2}e^{-\varphi_{i,t_{1}}}d\nu_{t_{1}}^{i}\right)^{\frac{1}{2}}.
\]
Let $\chi_{r,\epsilon,i}\in C_{c}^{\infty}(M_{i})$ be as in Lemma
\ref{Hormander}, so that 
\begin{align*}
\int_{M_{i}}\left(Rc_{g_{i,t_{1}}}+\nabla\overline{\nabla}(f_{i,t_{1}}+\varphi_{i,t_{1}})\right) & (\chi_{r,\epsilon,i}\eta_{i,t_{1}}^{(1)},\overline{\chi_{r,j,i}\eta_{i,t_{1}}^{(1)}})e^{-\varphi_{i,t_{1}}}d\nu_{t_{1}}^{i}\\
\geq & \int_{M_{i}}\beta|\eta_{i,t_{1}}^{(1)}|^{2}e^{-\varphi_{i,t_{1}}}d\nu_{t_{1}}^{i}-\beta\int_{M_{i}}(1-\chi_{r,\epsilon,i}^{2})|\eta_{i,t_{1}}^{(1)}|^{2}e^{-\varphi_{i,t_{1}}}d\nu_{t_{1}}^{i}\\
 & -\left(\sup_{\text{supp}(\chi_{r,\epsilon,i})}\left(Rc_{g_{i,t_{1}}}+\nabla\overline{\nabla}(f_{i,t_{1}}+\varphi_{i,t_{1}})-\beta\omega_{i,t_{1}}\right)_{-}\right)\int_{M_{i}}\chi_{r,j,i}^{2}|\eta_{i,t_{1}}^{(1)}|^{2}e^{-\varphi_{i,t_{1}}}d\nu_{t_{1}}^{i}\\
\geq & \left(\beta-\Psi(i^{-1}|r,\epsilon,t_{1})-\Psi(\epsilon|r,t_{1})-\Psi(r^{-1}|t_{1})\right)\int_{M_{i}}|\eta_{i,t_{1}}^{(1)}|^{2}e^{-\varphi_{i,t_{1}}}d\nu_{t_{1}}^{i}.
\end{align*}
Analogously to Lemma \ref{Hormander}, this implies
\begin{align*}
\int_{M_{i}}|\overline{\partial}_{f_{i}+\varphi_{i}}^{\ast}(\chi_{r,\epsilon,i}\eta_{i,t_{1}}^{(1)})|^{2}e^{-\varphi_{i,t_{1}}}d\nu_{t_{1}}^{i} & +\int_{M_{i}}|\overline{\partial}\chi_{r,\epsilon,i}\wedge\eta_{i,t_{1}}^{(1)}|^{2}e^{-\varphi_{i,t_{1}}}d\nu_{t_{1}}^{i}\\
= & \int_{M_{i}}|\overline{\nabla}(\chi_{r,\epsilon,i}\eta_{i,t_{1}}^{(1)})|^{2}e^{-\varphi_{i,t_{1}}}d\nu_{t_{1}}^{i}\\
 & +\int_{M_{i}}\left(Rc_{g_{i,t_{1}}}+\nabla\overline{\nabla}(f_{i,t_{1}}+\varphi_{i,t_{1}})\right)(\chi_{r,\epsilon,i}\eta_{i,t_{1}}^{(1)},\overline{\chi_{r,\epsilon,i}\eta_{i,t_{1}}^{(1)}})e^{-\varphi_{i,t_{1}}}d\nu_{t_{1}}^{i}\\
\geq & \int_{M_{i}}\beta|\eta_{i,t_{1}}^{(1)}|^{2}e^{-\varphi_{i,t_{1}}}d\nu_{t_{1}}^{i}-\Psi(i^{-1}|r,\epsilon,t_{0})-\Psi(\epsilon|r,t_{0})-\Psi(r^{-1},t_{0}).
\end{align*}
Almost exactly as in the $\varphi_{i}\equiv0$ case, we also have
\begin{align*}
\int_{M_{i}}\left(|\overline{\partial}_{f_{i}+\varphi_{i}}^{\ast}(\chi_{r,\epsilon,i}\eta_{i,t_{1}}^{(1)})|^{2}+|\overline{\partial}(\chi_{r,\epsilon,i}\eta_{i,t_{1}}^{(1)})|^{2}\right)e^{-\varphi_{i,t_{1}}}d\nu_{t_{1}}^{i}\leq & \int_{\text{supp}(\chi_{r,\epsilon,i})}|\overline{\partial}_{f_{i}+\varphi_{i}}^{\ast}\eta_{i,t_{1}}^{(1)}|^{2}e^{-\varphi_{i,t_{1}}}d\nu_{t_{1}}^{i}\\
 & +\Psi(i^{-1}|\epsilon,r,t_{0})+\Psi(\epsilon|r,t_{0})+\Psi(r^{-1}|t_{0}).
\end{align*}
Next, we recall that 
\begin{align*}
(\partial_{t}-\Delta)\left(-aw_{i}+\frac{|\eta_{i}^{(2)}|^{2}}{a}+\tau|\overline{\partial}_{h_{i}+\varphi_{i}}^{\ast}\eta_{i}^{(2)}|\right)\leq & \tau|\nabla\overline{\nabla}\varphi_{i}|\cdot|\overline{\nabla}\eta_{i}^{(2)}|+\tau|\nabla\nabla\varphi_{i}|\cdot|\nabla\eta_{i}^{(2)}|\\
 & +\tau\left|\nabla\left((\partial_{t}-\Delta)\varphi_{i}-\frac{1}{\tau}\varphi_{i}\right)\right|\cdot|\eta_{i}^{(2)}|.\\
\leq & C\left(|\overline{\nabla}\eta_{i}^{(2)}|+|\nabla\eta_{i}^{(2)}|+|\eta_{i}^{(2)}|\right),
\end{align*}
where $C$ is independent of $i\in\mathbb{N}$. We also have 
\begin{align*}
\left|(\partial_{t}+\Delta-R)\left((4\pi\tau)^{-n}e^{-(\varphi_{i}+f_{i})}\right)\right|= & \left|\left(-(\partial_{t}+\Delta)\varphi_{i}+|\nabla\varphi_{i}|^{2}+2\text{Re}\langle\nabla\varphi_{i},\overline{\nabla}f_{i}\rangle\right)\right|(4\pi\tau)^{-n}e^{-(\varphi+f_{i})}\\
\leq & C(4\pi\tau)^{-n}e^{-(\varphi+f_{i})}
\end{align*}
Combining expressions allows us to estimate (since we have uniform
estimates for $|\nabla f_{i}|$ on $\text{supp}(\nabla\varphi_{i})$)
\begin{align*}
\frac{d}{dt}e^{-At}\int_{M_{i}}\left(-aw_{i}+\frac{|\eta_{i}^{(2)}|^{2}}{a}+\tau|\overline{\partial}_{h_{i}+\varphi_{i}}^{\ast}\eta_{i}^{(2)}|\right)e^{-\varphi_{i}}d\nu_{t}^{i}\leq & Ce^{-At}\int_{M_{i}}\left(|\overline{\nabla}\eta_{i}^{(2)}|+|\nabla\eta_{i}^{(2)}|+|\eta_{i}^{(2)}|\right)e^{-\varphi_{i}}d\nu_{t}^{i}\\
 & +a(C+A)e^{-At}\int_{M_{i}}|w_{i}|e^{-\varphi_{i}}d\nu_{t}^{i}
\end{align*}
if we choose $A<\infty$ sufficiently large (again independently of
$i\in\mathbb{N}$). Integrating in time therefore leads to 
\begin{align*}
\int_{M_{i}}\left(-aw_{i}+\frac{|\eta_{i,t_{1}}^{(2)}|^{2}}{a}+\tau|\overline{\partial}_{h_{i}+\varphi_{i}}^{\ast}\eta_{i,t_{1}}^{(2)}|\right)e^{-\varphi_{i}}d\nu_{t_{1}}^{i}\leq & e^{A(t_{1}-t_{0})}\int_{M_{i}}\left(-aw_{i}+\frac{|\eta_{i}^{(2)}|^{2}}{a}+\tau|\overline{\partial}_{h_{i}+\varphi_{i}}^{\ast}\eta_{i}^{(2)}|\right)e^{-\varphi_{i}}d\nu_{t_{0}}^{i}\\
 & +\Psi(i^{-1}|t_{1},t_{0})+\Psi(|t_{1}-t_{0}|),
\end{align*}
so that (choosing $a>0$ small as in Lemma \ref{Hormander})
\begin{align*}
\int_{M_{i}}|\overline{\partial}_{h_{i}+\varphi_{i}}^{\ast}\eta_{i,t_{1}}^{(2)}|e^{-\varphi_{i}}d\nu_{t_{1}}^{i}\leq & e^{A(t_{1}-t_{0})}\int_{M_{i}}|\overline{\partial}_{h_{i}+\varphi_{i}}^{\ast}\eta_{i}^{(2)}|e^{-\varphi_{i}}d\nu_{t_{0}}^{i}+\Psi(i^{-1}|t_{1},t_{0})+\Psi(|t_{1}-t_{0}|)\\
\leq & \Psi(i^{-1}|t_{1},t_{0})+\Psi(|t_{1}-t_{0}|).
\end{align*}
Arguing as in the $\varphi_{i}\equiv0$ case leads to
\[
\int_{\text{supp}(\chi_{r,\epsilon,i})}|\overline{\partial}_{h_{i}+\varphi_{i}}^{\ast}\eta_{i,t}^{(2)}|^{2}e^{-\varphi_{i}}d\nu_{t_{1}}^{i}\leq\Psi(i^{-1}|\epsilon,r,t_{0},t_{1})+\Psi(\epsilon|r,t_{0},t_{1})+\Psi(r^{-1}|t_{0},t_{1}),
\]
so that combining estimates gives 
\begin{align*}
\int_{M_{i}}\beta|\eta_{i,t_{1}}^{(1)}|^{2}e^{-\varphi_{i,t_{1}}}d\nu_{t_{1}}^{i}\leq & \int_{\text{supp}(\chi_{r,\epsilon,i})}|\overline{\partial}_{h_{i,t_{1}}+\varphi_{i,t_{1}}}^{\ast}\eta_{i,t_{1}}|^{2}e^{-\varphi_{i,t_{1}}}d\nu_{t_{1}}^{i}+\Psi(i^{-1}|t_{1},t_{0},\epsilon,r)+\Psi(|t_{1}-t_{0}||\epsilon,r)\\
 & +\Psi(\epsilon|r)+\Psi(r^{-1}).
\end{align*}
Next, observe that if $A<\infty$ is chosen sufficiently large, then
we have
\[
(\partial_{t}-\Delta)\left(e^{-At}\left|-aw_{i}+\frac{|\eta_{i}|^{2}}{a}+\tau|\overline{\partial}_{h_{i}+\varphi_{i}}^{\ast}\eta_{i}|\right|e^{-\varphi_{i}}\right)\leq e^{-At}\left(|\overline{\nabla}\eta_{i}|+|\nabla\eta_{i}|+|\eta_{i}|\right)\chi_{\text{supp}(\varphi_{i})},
\]
so that for all $x\in\text{supp}(\chi_{r,\epsilon,i})$, we have 
\begin{align*}
\tau^{2}|\overline{\partial}_{h_{i}}^{\ast}\eta_{i,t_{1}}(x)|^{2}\leq & e^{A(t_{1}-t_{0})}\left(C\sigma^{-1}a^{2}w_{i,t_{1}}^{2}(x)+(1+\sigma)\left(\int_{M_{i}}\tau|\overline{\partial}_{h_{i}}^{\ast}\eta_{i,t_{0}}|d\nu_{x,t_{1};t_{0}}^{i}\right)^{2}+\frac{C(n)}{a^{2}\sigma}\left(\int_{M_{i}}|\eta_{i,t_{0}}|^{2}d\nu_{x,t_{1};t_{0}}^{i}\right)^{2}\right)\\
 & +Ce^{-At}\left(|\overline{\nabla}\eta_{i}|^{2}+|\nabla\eta_{i}|^{2}+|\eta_{i}|^{2}\right)\chi_{\text{supp}(\varphi_{i})}+C\frac{a^{2}}{\sigma}\left(\int_{M_{i}}|w_{i,t_{0}}|d\nu_{x,t_{1};t_{0}}^{i}\right)^{2},
\end{align*}
so we can argue as in Lemma \ref{Hormander} to get
\begin{align*}
|t_{1}|^{2}\int_{\text{supp}(\chi_{r,\epsilon,i})}|\overline{\partial}_{h_{i}}^{\ast}\eta_{i,t_{1}}|^{2}d\nu_{t_{1}}^{i}\leq & a^{2}\sigma^{-1}e^{A(t_{1}-t_{0})}\Psi(i^{-1}|\epsilon,r,t_{0},t_{1})+(1+\sigma)|t_{0}|^{2}e^{A(t_{1}-t_{0})}\int_{M_{i}}|\overline{\partial}_{h_{i}}^{\ast}\eta_{i,t_{0}}|^{2}d\nu_{t_{0}}^{i}\\
 & +e^{A(t_{1}-t_{0})}\frac{C(n)}{a^{2}\sigma}\sup_{M_{i}}|\eta_{i,t_{0}}|^{4}+C|t_{1}-t_{0}|e^{-A(t_{1}-t_{0})}\sup_{\text{supp}(\varphi_{i})}\left(|\overline{\nabla}\eta_{i}|^{2}+|\nabla\eta_{i}|^{2}+|\eta_{i}|^{2}\right)\\
\leq & \frac{a^{2}}{\sigma}\Psi(i^{-1}|\epsilon,r)+\frac{C}{a^{2}\sigma}+C(1+\sigma)|t_{0}|^{2}e^{A(t_{1}-t_{0})}\int_{M_{i}}|\overline{\partial}_{h_{i}}^{\ast}\eta_{i,t_{0}}|^{2}d\nu_{t_{0}}^{i}.
\end{align*}
This is because $\eta_{i,t_{0}}$ satisfy uniform $C_{\text{loc}}^{\infty}$
estimates, and because the geometry is controlled near $\text{supp}(\varphi_{i})$,
so we can apply local parabolic estimates to get bounds for $|\nabla\eta_{i}|,|\overline{\nabla}\eta_{i}|$
there (uniformly in $i\in\mathbb{N}$). Choosing $\sigma>0$ small
and then $a<\infty$ large gives
\begin{align*}
\int_{\text{supp}(\chi_{r,\epsilon,i})}|\overline{\partial}_{h_{i}+\varphi_{i}}^{\ast}\eta_{i}|^{2}e^{-\varphi_{i}}d\nu_{t_{1}}^{i}\leq & (1+\Psi(|t_{1}-t_{0}||\epsilon,r^{-1}))\int_{M_{i}}|\overline{\partial}_{h_{i}+\varphi_{i}}^{\ast}\eta_{i}|^{2}e^{-\varphi_{i}}d\nu_{t_{0}}^{i}\\
 & +\Psi(i^{-1}|t_{1},t_{0},\epsilon,r^{-1})+\Psi(|t_{1}-t_{0}||\epsilon,r^{-1})+\Psi(\epsilon|r^{-1})+\Psi(r^{-1}).
\end{align*}
Again choosing $r<\infty$ large and then $\epsilon>0$ small, and
combining estimates, we get
\begin{align*}
\left|\int_{M_{i}}\langle\eta_{i,t_{0}},\overline{\overline{\partial}v_{i,t_{0}}}\rangle e^{-\varphi_{t_{0}}}d\nu_{t_{0}}^{i}\right|\leq & (1+\Psi(|t_{1}-t_{0}|)\left(\int_{M_{i}}|\overline{\partial}_{h_{i}+\varphi_{i}}^{\ast}\eta_{i}|^{2}e^{-\varphi_{i}}d\nu_{t_{0}}^{i}\right)^{\frac{1}{2}}\\
 & \times\left(\int_{M_{i}}\frac{1}{\beta}|\overline{\partial}v_{i,t_{1}}|^{2}e^{-\varphi_{i,t_{1}}}d\nu_{t_{1}}^{i}\right)^{\frac{1}{2}}+\Psi(|t_{1}-t_{0}|)+\Psi(i^{-1}|t_{0},t_{1}).
\end{align*}
We take $i\to\infty$ to obtain
\begin{align*}
\left|\int_{\mathcal{R}_{t_{0}}}\langle\eta_{t_{0}},\overline{\partial}v_{t_{0}}\rangle e^{-\varphi_{t_{0}}}d\nu_{t_{0}}\right|\leq & (1+\Psi(|t_{1}-t_{0}|)\left(\int_{\mathcal{R}_{t_{0}}}|\overline{\partial}_{f_{t_{0}}+\varphi_{t_{0}}}^{\ast}\eta_{t_{0}}|^{2}e^{-\varphi_{t_{0}}}d\nu_{t_{0}}\right)^{\frac{1}{2}}\\
 & \times\left(\frac{1}{\beta}\int_{\mathcal{R}_{t_{1}}}|\overline{\partial}v|^{2}e^{-\varphi_{t_{1}}}d\nu_{t_{1}}\right)^{\frac{1}{2}}+\Psi(|t_{1}-t_{0}|).
\end{align*}
Finally, taking $t_{0}\nearrow t_{1}$ gives the claim (arguing as
in the beginning of Lemma \ref{Hormander}). 
\end{proof}
As before, the Hörmander estimate implies an $L^{2}$ existence theorem.
\begin{cor}
\label{lemmaforseparation} With the same notation as in Lemma \ref{HormanderwWeights},
there exists $u\in L^{2}(\mathcal{R}_{t_{1}},e^{-\varphi}d\nu_{t_{1}})$
satisfying $\overline{\partial}u=\overline{\partial}v$ and 
\[
\int_{\mathcal{R}_{t_{1}}}|u|^{2}d\nu_{t_{1}}\leq\int_{\mathcal{R}_{t_{1}}}\frac{1}{\beta}|\overline{\partial}v|^{2}d\nu_{t_{1}}.
\]
\end{cor}

\begin{rem}
The usual approximation argument tells us that $\varphi$ need not
be smooth, and we can in fact take $\varphi$ to have logarithmic
singularities (see the proof of Proposition \ref{separate}). 
\end{rem}

\begin{proof}[Proof of Theorem 2]
 This is just the statement of Corollary \ref{lemmaforseparation}. 
\end{proof}

\section{Constructing Local Holomorphic Embeddings}

We now indicate how Theorem \ref{thm2} allows us to show that any
tangent cone of a singular soliton is homeomorphic to a normal analytic
space. Much of the arguments are taken from \cite{donaldsun1,donaldsun2,szek1,szek2},
so certain details will be omitted, and the reader will be directed
to those references when necessary. 

Our first step will be constructing peak holomorphic functions, using
the partial regularity theory developed in \cite{bamlergen3} and
the arguments of \cite{donaldsun1,CDSII}. We first include some definitions
to clarify the notion of convergence we will need.

We first show that the holomorphic functions on $\mathcal{R}$ constructed
by $L^{2}$ methods are locally Lipschitz on $\mathcal{X}$, and in
particular extend to continuous functions along the singular set of
$\mathcal{X}_{t_{0}}$. The proof is similar to that of Proposition
10 in \cite{szek2}.

Throughout this section, we will assume $\mathcal{X}$ is a singular
soliton obtained as an $\mathbb{F}$-limit as in (\ref{F}), with
the same hypotheses as in Section 4, but now also assuming $(\mathcal{R},g)$
is Ricci-flat, hence $(\mathcal{X}_{-1},d_{-1})$ corresponds to a
metric cone with vertex $x_{\ast}\in\mathcal{X}_{-1}$.
\begin{lem}
\label{LocLipschitz} Assume $\mathcal{X}$ is also a static cone,
and suppose $u:\mathcal{R}_{t_{0}}\to\mathbb{C}$ is a holomorphic
function satisfying
\[
\int_{\mathcal{R}_{t_{0}}}|u|^{2}d\nu_{t_{0}}<\infty.
\]

$(i)$ For all $x\in\mathcal{R}_{t_{0}}$ and $t\in(t_{0},0)$, we
have $u\in L^{1}(\mathcal{R}_{t_{0}},\nu_{x(t);t_{0}})$ and
\[
u(x)=\int_{\mathcal{R}_{t_{0}}}ud\nu_{x(t);t_{0}}=\int_{\mathcal{R}_{t_{0}}}K(x(t);y)u(y)dg_{t_{0}}(y).
\]

$(ii)$ For any compact subset $L\subseteq\mathcal{X}_{t_{0}}$, there
exists $C=C(L)<\infty$ such that 
\[
\sup_{L\cap\mathcal{R}_{t_{0}}}(|u|+|\nabla u|)\leq C\int_{\mathcal{R}_{t_{0}}}|u|^{2}d\nu_{t_{0}}.
\]
\end{lem}

\begin{rem}
\label{holoflow} Another way of restating this lemma is that if $(u_{t})_{t\in(-\infty,0)}$
is defined by $\partial_{\mathfrak{t}}u_{t}\equiv0$ and $u_{t_{0}}=u$,
then $(\partial_{\mathfrak{t}}-\Delta)u_{t}\equiv0$. In fact, the
lemma gives
\[
u_{t}(x)=u(x(t))=\int_{\mathcal{R}_{t_{0}}}u_{t_{0}}d\nu_{x;t_{0}},
\]
which exactly means that $(u_{t})$ is the heat flow starting from
$u_{t_{0}}$. 
\end{rem}

\begin{proof}
$(i)$ By \cite{bamlergen3}, we can find a sequence $(\eta_{\epsilon})$
in $C_{c}^{\infty}(\mathcal{R}_{t_{0}})$ such that $\eta_{\epsilon}|_{\{r_{Rm}\geq2\epsilon\}}\equiv1$,
$\text{supp}(\eta_{\epsilon})\subseteq\{r_{Rm}>\epsilon\}$, and
\[
\lim_{\epsilon\to0}\int_{\mathcal{R}_{t_{0}}\cap B(x_{0},\rho)}|\nabla\eta_{\epsilon}|^{\frac{7}{2}}dgdt=0
\]
for all $\rho\in[1,\infty)$. Let $\chi\in C^{\infty}(\mathbb{R})$
satisfy $\chi|_{[2,\infty)}\equiv0$ and $\chi|_{(-\infty,1]}\equiv1$,
and set $\chi_{\rho}(x):=\chi(\rho^{-1}r(x))$, where $r:=d_{t_{0}}(x_{0},\cdot)$,
and $x_{0}$ is the vertex of $\mathcal{X}_{t_{0}}$. For fixed $\epsilon>0$
and $\rho<\infty$, we thus have $\text{supp}(\chi_{\rho}\eta_{\epsilon}u)\subset\subset\mathcal{R}_{t_{0}},$where
$x_{0}\in\mathcal{X}_{t_{0}}$ is fixed and $r:=d_{t_{0}}(x_{0},\cdot)$.
For fixed $x\in\mathcal{R}_{t_{0}}$, $t\in(t_{0},0)$, we can estimate
(arguing as in Lemma \ref{goodgrad})
\begin{align*}
\int_{\mathcal{R}_{t_{0}}}|u|d\nu_{x(t);t_{0}}\leq & C(x,t)\int_{\mathcal{R}_{t_{0}}}|u|e^{\frac{1}{4}f}d\nu_{t_{0}}\leq C(x,t)\left(\int_{\mathcal{R}_{t_{0}}}|u|^{2}d\nu_{t_{0}}\right)^{\frac{1}{2}}\left(\int_{\mathcal{R}_{t_{0}}}e^{\frac{1}{2}f}d\nu_{t_{0}}\right)^{\frac{1}{2}}<\infty,
\end{align*}
so that the dominated convergence theorem ensures
\begin{equation}
\int_{\mathcal{R}_{t_{0}}}ud\nu_{x(t);t_{0}}=\lim_{r\to\infty}\lim_{\epsilon\searrow0}\int_{\mathcal{R}_{t_{0}}}u\chi_{\rho}\eta_{\epsilon}d\nu_{x(t);t_{0}}.\label{domconverge}
\end{equation}

By Theorem 15.60 of \cite{bamlergen3}, 
\[
\partial_{\mathfrak{t}}|_{x}K(\cdot;y)+\partial_{\mathfrak{t}}|_{y}K(x;\cdot)=0
\]
for any $y\in\mathcal{R}_{t_{0}}$ and $x\in\mathcal{R}_{>t_{0}}$.
Using this and $\partial_{\mathfrak{t}}|_{y}K(x;\cdot)=-\Delta_{y}K(x;y)$,
we may compute (for fixed $x\in\mathcal{R}_{t_{0}}$)
\begin{align*}
\frac{d}{dt}\int_{\mathcal{R}_{t_{0}}}K(x(t);y)\chi_{\rho}(y)\eta_{\epsilon}(y)u(y)dg_{t_{0}}(y)= & \int_{\mathcal{R}_{t_{0}}}\partial_{\mathfrak{t}}|_{x(t)}K(\cdot;y)\chi_{\rho}(y)\eta_{\epsilon}(y)u(y)dg_{t_{0}}(y)\\
= & \int_{\mathcal{R}_{t_{0}}}\Delta_{y}K(x(t);y)\chi_{\rho}(y)\eta_{\epsilon}(y)u(y)dg_{t_{0}}(y)\\
= & -\int_{\mathcal{R}_{t_{0}}}\langle\nabla_{y}\log K(x(t);y),\nabla\chi_{\rho}(y)\rangle\eta_{\epsilon}(y)u(y)d\nu_{x(t);t_{0}}(y)\\
 & -\int_{\mathcal{R}}\langle\nabla_{y}\log K(x(t);y),\nabla\eta_{\epsilon}(y)\rangle\chi_{\rho}(y)u(y)d\nu_{x(t);t_{0}}(y),
\end{align*}
where we used $\overline{\partial}u\equiv0$. Next, we use Hölder's
inequality to estimate
\begin{align*}
\biggl|\int_{\mathcal{R}_{t_{0}}}\langle\nabla_{y}\log K(x(t);y),\nabla\chi_{\rho}(y)\rangle & \eta_{\epsilon}(y)u(y)d\nu_{x(t);t_{0}}(y)\biggr|\\
\leq & \frac{C}{\rho}\int_{\mathcal{R}_{t_{0}}}|\nabla_{y}\log K(x(t);y)|\cdot|u(y)|d\nu_{x(t);t_{0}}(y)\\
\leq & \frac{C}{\rho}\left(\int_{\mathcal{R}_{t_{0}}}|\nabla_{y}\log K(x(t);y)|^{3}d\nu_{x(t);t_{0}}(y)\right)^{\frac{1}{3}}\left(\int_{\mathcal{R}_{t_{0}}}|u|^{\frac{3}{2}}d\nu_{x(t);t_{0}}(y)\right)^{\frac{2}{3}},
\end{align*}
and 
\begin{align*}
\biggl|\int_{\mathcal{R}_{t_{0}}}\langle\nabla_{y}\log K(x(t);y),\nabla\eta_{\epsilon}(y)\rangle & \chi_{\rho}(y)u(y)d\nu_{x(t);t_{0}}(y)\biggr|\\
\leq & \int_{B(x_{0},2\rho)}|\nabla_{y}\log K(x(t);y)|\cdot|\nabla\eta_{\epsilon}(y)|\cdot|u(y)|d\nu_{x(t);t_{0}}(y)\\
\leq & \left(\int_{B(x_{0},2\rho)}|\nabla\eta_{\epsilon}(y)|^{\frac{7}{2}}d\nu_{x(t);t_{0}}(y)\right)^{\frac{2}{7}}\left(\int_{\mathcal{R}_{t_{0}}}|\nabla_{y}\log K(x(t);y)|^{21}d\nu_{x(t);t_{0}}(y)\right)^{\frac{1}{21}}\\
 & \times\left(\int_{\mathcal{R}_{t_{0}}}|u|^{\frac{3}{2}}d\nu_{x(t);t_{0}}(y)\right)^{\frac{2}{3}}\\
\leq & \Psi(\epsilon|\rho)\left(\int_{\mathcal{R}_{t_{0}}}|\nabla_{y}\log K(x(t);y)|^{21}d\nu_{x(t);t_{0}}(y)\right)^{\frac{1}{21}}\left(\int_{\mathcal{R}_{t_{0}}}|u|^{\frac{3}{2}}d\nu_{x(t);t_{0}}(y)\right)^{\frac{2}{3}}
\end{align*}
where we also used the on-diagonal heat kernel upper bound (Lemma
15.9 of \cite{bamlergen3}). Moreover, Hölder's inequality and conjugate
heat kernel comparison (see Proposition 8.1 of \cite{bamlergen3}
and Lemma \ref{goodgrad}) give
\begin{align*}
\int_{\mathcal{R}_{t_{0}}}|u|^{\frac{3}{2}}d\nu_{x(t);t_{0}}\leq & C(x,t,t_{0})\int_{\mathcal{R}_{t_{0}}}|u|^{\frac{3}{2}}e^{\frac{1}{8}f}d\nu_{t_{0}}\leq C(x,t,t_{0},\alpha)\left(\int_{\mathcal{R}_{t_{0}}}|u|^{2}d\nu_{t_{0}}\right)^{\frac{3}{4}}\left(\int_{\mathcal{R}_{t_{0}}}e^{\frac{1}{2}f}d\nu_{t_{0}}\right)^{\frac{1}{4}}\\
\leq & C(x,t_{0})\left(\int_{\mathcal{R}_{t_{0}}}|u|^{2}d\nu_{t_{0}}\right)^{\frac{3}{4}},
\end{align*}
From Proposition 4.2 of \cite{bamlergen1} and Lemma \ref{conjheatgradient},
we also have
\begin{align*}
\int_{\mathcal{R}_{t_{0}}}|\nabla_{y}\log K(x(t);y)|^{p}d\nu_{x(t);t_{0}}(y)\leq & C(n,p)(t-t_{0})^{\frac{p}{2}}.
\end{align*}
Combining estimates gives
\begin{align*}
\left|\frac{d}{dt}\int_{\mathcal{R}_{t_{0}}}K(x(t);y)\chi_{\rho}(y)u(y)\eta_{\epsilon}(y)dg_{t_{0}}(y)\right|\leq & \Psi(\epsilon|\rho,t_{0})+\Psi(\rho^{-1}|t_{0}).
\end{align*}
Moreover, Lemma \ref{fundsol} implies that for $\rho\geq\underline{\rho}(x)$
and $\epsilon\leq\overline{\epsilon}(x)$, we have
\[
\lim_{t\searrow t_{0}}\int_{\mathcal{R}_{t_{0}}}K(x(t);y)\chi_{\rho}(y)u(y)\eta_{\epsilon}(y)dg_{t_{0}}(y)=u(x).
\]
The claim follows by choosing $t-t_{0}$ small, and then integrating
in time and appealing to (\ref{domconverge}). 

\noindent $(ii)$ This follows from heat kernel comparison (Proposition
8.1 of \cite{bamlergen3}, see also Lemma \ref{goodgrad}) once we
show that 
\[
\nabla u_{t}(x)=\int_{\mathcal{R}_{t_{0}}}u_{t_{0}}(y)\nabla_{x}\log K(x;y)d\nu_{x;t_{0}}(y)
\]
for all $x\in\mathcal{R}_{t}$. We just need to verify that 
\[
\int_{\mathcal{R}_{t_{0}}}|u_{t_{0}}|(y)|\nabla_{x}\log K(x;y)|d\nu_{x;t_{0}}(y)<\infty,
\]
but this follows from Proposition 4.2 of \cite{bamlergen1} (as in
Lemma \ref{goodgrad}).
\end{proof}
\begin{defn}
If $u_{t_{0}}\in L^{2}(\mathcal{R}_{t_{0}},d\nu_{t_{0}})$ is holomorphic
for some $t_{0}<0$, then we call (by abuse of notation) the unique
extension to $u\in C(\mathcal{X},\mathbb{C})\cap C^{\infty}(\mathcal{R},\mathbb{C})$
satisfying $\partial_{\mathfrak{t}}u\equiv0$ holomorphic. 
\end{defn}

\begin{rem}
If $u$ is a holomorphic function, then $\overline{\partial}u_{t}=0$
for each $t\in(-\infty,0)$ since $\mathcal{L}_{\partial_{\mathfrak{t}}}J\equiv0$. 
\end{rem}

\begin{defn}
\label{polarizeddef} A (pointed) polarized Kähler singular space
is a tuple $(X,d,p,\mathcal{R},g,J,L,h)$ where $(X,d,p)$ is a pointed
time slice of a metric flow obtained as an $\mathbb{F}$-limit of
closed Kähler-Ricci flows as in (\ref{F}), $(\mathcal{R},g,J)$ is
its regular part, and $(L,h)$ is a holomorphic Hermitian line bundle
over $\mathcal{R}$ whose Chern connection has curvature equal to
the Kähler form of $(g,J)$. We say a sequence of polarized Kähler
singular spaces $(X_{j},d_{j},p_{j},\mathcal{R}_{j},g_{j},J_{j},L_{j},h_{j})$
converges to a polarized Kähler singular space $(X,d,p,\mathcal{R},g,J,L,h)$
if there are diffeomorphisms $\psi_{j}:U_{j}\to V_{j}$, where $(U_{j})$
is a precompact exhaustion of $\mathcal{R}$ and $V_{i}\subseteq\mathcal{R}_{i}$
are open, such that the following hold:

$(i)$ $\psi_{j}^{\ast}g_{j}\to g$, $\psi_{j}^{\ast}J_{j}\to J$
in $C_{\text{loc}}^{\infty}(\mathcal{R})$,

$(ii)$ $\psi_{j}$ realize pointed Gromov-Hausdorff convergence $(X_{j},d_{j},p_{j})\to(X,d,p)$,
in the sense that they can be taken to be restrictions of $\epsilon_{j}$-Gromov-Hausdorff
isometries $(X,d,p)\to(X_{j},d_{j},p_{j})$,

$(iii)$ there are line bundle isomorphisms $\xi_{j}:L|_{U_{j}}\to\psi_{j}^{\ast}(L_{j}|_{V_{j}})$
such that $\xi_{j}^{\ast}h_{j}\to h$ and $\xi_{j}^{\ast}A_{j}\to A$
in $C_{\text{loc}}^{\infty}(\mathcal{R})$, where $A_{j}$ is the
Chern connection of $\omega_{j}$. 
\end{defn}

As a preliminary step, we prove a uniform $C^{1}$ bound for sections
of the Hermitian line bundle $(\mathcal{R}_{-1},e^{-f_{-1}})$, using
the heat kernel estimates proved in the Appendix. 
\begin{lem}
\label{therightway} There exists $C=C(Y)<\infty$ such that any holomorphic
$u\in L^{2}(\mathcal{R}_{-1},d\nu_{-1})$ satisfies the following:

$(i)$ $|u|_{h}^{2}\leq C\int_{\mathcal{R}_{-1}}|u|_{h}^{2}dg$,

$(ii)$ $|\nabla^{h}u|_{h}^{2}\leq C\int_{\mathcal{R}_{-1}}|u|_{h}^{2}dg$,

\noindent where we view $h=(2\pi)^{-n}e^{-f_{-1}}$ as a Hermitian
metric on $\mathcal{O}_{\mathcal{R}_{-1}}$, with corresponding Chern
connection $\nabla^{h}=\partial+\partial f\wedge$. 
\end{lem}

\begin{proof}
Let $\chi_{\rho},\eta_{\epsilon}$ be as in Lemma \ref{LocLipschitz},
but recall that we now denote the vertex of $\mathcal{X}_{-1}$ by
$x_{\ast}$. Set $\mathcal{R}:=\mathcal{R}_{-1}$, and let $K:\mathcal{R}\times\mathcal{R}\times(0,\infty)\to(0,\infty)$
be the heat kernel of $\mathcal{R}$ (as in the Appendix). Then we
can integrate 
\begin{equation}
\Delta|u|_{h}^{2}\geq|\nabla^{h}u|_{h}^{2}-n|u|_{h}^{2}\label{c0elliptic}
\end{equation}
to obtain
\begin{align*}
\int_{\mathcal{R}}\chi_{\rho}^{2}(y)\eta_{\epsilon}^{2}(y)|\nabla^{h}u|_{h}^{2}(y)dg(y)\leq & \int_{\mathcal{R}}\chi_{\rho}^{2}(y)\eta_{\epsilon}^{2}(y)\Delta|u|_{h}^{2}(y)dg(y)+n\int_{\mathcal{R}}\chi_{\rho}^{2}(y)\eta_{\epsilon}^{2}(y)|u|_{h}^{2}(y)dg(y)\\
\leq & 8\int_{\mathcal{R}}(\chi_{\rho}\eta_{\epsilon}|\nabla^{h}u|_{h})(y)|u|_{h}(y)\left(|\nabla\chi_{\rho}|(y)+\chi_{\rho}|\nabla\eta_{\epsilon}|(y)\right)dg(y)\\
 & +n\int_{\mathcal{R}}\chi_{\rho}^{2}(y)\eta_{\epsilon}^{2}(y)|u|_{h}^{2}(y)dg(y)\\
\leq & \delta\int_{\mathcal{R}}\chi_{\rho}^{2}(y)\eta_{\epsilon}^{2}(y)|\nabla^{h}u|_{h}^{2}(y)dg(y)+n\int_{\mathcal{R}}\chi_{\rho}^{2}(y)\eta_{\epsilon}^{2}(y)|u|_{h}^{2}(y)dg(y)\\
 & +C(\delta)\int_{\mathcal{R}}|u|_{h}^{2}(y)\left(\frac{1}{\rho^{2}}+\chi_{\rho}^{2}|\nabla\eta_{\epsilon}|^{2}(y)\right)dg(y).
\end{align*}
Because $|u|_{h}$ is Locally bounded by Lemma \ref{LocLipschitz},
we can take $\epsilon\searrow0$, then $\rho\to\infty$, then $\delta\searrow0$
to obtain
\begin{equation}
\int_{\mathcal{R}}|\nabla^{h}u|_{h}^{2}(y)dg(y)\leq n\int_{\mathcal{R}}|u|_{h}^{2}(y)dg(y).\label{IBPidentity}
\end{equation}
Next, we integrate (\ref{c0elliptic}) against the heat kernel of
$\mathcal{R}$ to obtain
\begin{align*}
\frac{d}{dt}\int_{\mathcal{R}}|u|_{h}^{2}(y) & \chi_{\rho}^{2}(y)\eta_{\epsilon}^{2}(y)K(x,y,1-t)dg(y)\\
= & -\int_{\mathcal{R}}\Delta\left(|u|_{h}^{2}\chi_{\rho}^{2}\eta_{\epsilon}^{2}\right)(y)K(x,y,1-t)dg(y)\\
\leq & \int_{\mathcal{R}}\left(n|u|_{h}^{2}-|\nabla^{h}u|_{h}^{2}\right)(y)\chi_{\rho}^{2}\eta_{\epsilon}^{2}(y)K(x,y,1-t)dg(y)\\
 & -2\text{Re}\int_{\mathcal{R}}\langle\nabla|u|_{h}^{2},\eta_{\epsilon}^{2}\overline{\nabla}\chi_{\rho}^{2}+\chi_{\rho}^{2}\overline{\nabla}\eta_{\epsilon}^{2}\rangle(y)K(x,y,1-t)dg(y)\\
 & -\text{Re}\int_{\mathcal{R}}|u|_{h}^{2}(y)\left(\eta_{\epsilon}^{2}\Delta\chi_{\rho}^{2}+2\langle\nabla\chi_{\rho}^{2},\overline{\nabla}\eta_{\epsilon}^{2}\rangle\right)(y)K(x,y,1-t)dg(y)\\
 & -\int_{\mathcal{R}}|u|_{h}^{2}(y)\chi_{\rho}^{2}(y)\Delta\eta_{\epsilon}^{2}(y)K(x,y,1-t)dg(y)\\
\leq & n\int_{\mathcal{R}}|u|_{h}^{2}\chi_{\rho}^{2}\eta_{\epsilon}^{2}(y)K(x,y,1-t)dg(y)\\
 & +C\int_{\mathcal{R}}|\nabla^{h}u|_{h}(y)(|\nabla\chi_{\rho}|+\chi_{\rho}|\nabla\eta_{\epsilon}|)(y)(|u|_{h}\eta_{\epsilon}\chi_{\rho})(y)K(x,y,1-t)dg(y)\\
 & +C\int_{\mathcal{R}}|u|_{h}^{2}(y)\left(|\Delta\chi_{\rho}^{2}|(y)+\chi_{\rho}|\nabla\eta_{\epsilon}|\right)K(x,y,1-t)dg(y)\\
 & +C\int_{\mathcal{R}}|u|_{h}^{2}(y)\chi_{\rho}^{2}(y)\eta_{\epsilon}(y)|\nabla\eta_{\epsilon}|(y)|\nabla_{y}K(x,y,1-t)|dg(y).
\end{align*}
Recalling that $|u|_{h},|\nabla^{h}u|_{h}$ are locally bounded by
Lemma \ref{LocLipschitz}, and appealing to (\ref{IBPidentity}),
we obtain 
\begin{align*}
\frac{d}{dt}e^{-nt}\int_{\mathcal{R}}|u|_{h}^{2}(y) & \chi_{\rho}^{2}(y)\eta_{\epsilon}^{2}(y)K(x,y,1-t)dg(y)\\
\leq & C(\rho)\int_{B(x_{\ast},2\rho)}|\nabla\eta_{\epsilon}|(y)\cdot|\nabla_{y}K(x,y,1-t)|dg(y)+\Psi(\rho^{-1})+\Psi(\epsilon|\rho)\\
\leq & C(\rho)\int_{B(x_{\ast},2\rho)}|\nabla\eta_{\epsilon}|(y)\frac{C(Y)}{(1-t)^{n+\frac{1}{2}}}\exp\left(-\frac{r_{0}^{2}}{10(1-t)}\right)dg(y)+\Psi(\rho^{-1})+\Psi(\epsilon|\rho)\\
\leq & \Psi(\rho^{-1})+\Psi(\epsilon|\rho)+\Psi(\epsilon|\rho,r_{0}),
\end{align*}
for $\epsilon\leq\overline{\epsilon}(r_{0})$, where $r_{0}>0$ is
such that $B(x,2r_{0})\subset\subset\mathcal{R}$. Integrating in
time, using Lemma \ref{fundsol}, and then taking $\epsilon\searrow0$
and $\rho\to\infty$ gives
\[
|u|_{h}^{2}(x)\leq e^{n}\int_{\mathcal{R}}|u|_{h}^{2}(y)K(x,y,1)dg(y)\leq C(Y)\int_{\mathcal{R}}|u|_{h}^{2}(y)dg(y).
\]
where in the last step, we used the on-diagonal heat kernel bound
of Lemma \ref{heatkernel}. 

Next, we integrate 
\begin{equation}
\Delta|\nabla^{h}u|_{h}^{2}\geq|\nabla^{h}\nabla^{h}u|_{h}^{2}-C(n)|\nabla^{h}u|_{h}^{2}\label{IBPgrad}
\end{equation}
against the heat kernel to obtain
\begin{align*}
\frac{d}{dt}\int_{\mathcal{R}}|\nabla^{h}u|_{h}^{2}(y) & \chi_{\rho}^{2}(y)\eta_{\epsilon}^{2}(y)K(x,y,1-t)dg(y)\\
= & -\int_{\mathcal{R}}\Delta\left(|\nabla^{h}u|_{h}^{2}\chi_{\rho}^{2}\eta_{\epsilon}^{2}\right)(y)K(x,y,1-t)dg(y)\\
\leq & \int_{\mathcal{R}}\left(C(n)|\nabla^{h}u|_{h}^{2}-|\nabla^{h}\nabla^{h}u|_{h}^{2}\right)(y)\chi_{\rho}^{2}\eta_{\epsilon}^{2}(y)K(x,y,1-t)dg(y)\\
 & -2\text{Re}\int_{\mathcal{R}}\langle\nabla|\nabla^{h}u|_{h}^{2},\eta_{\epsilon}^{2}\overline{\nabla}\chi_{\rho}^{2}+\chi_{\rho}^{2}\overline{\nabla}\eta_{\epsilon}^{2}\rangle(y)K(x,y,1-t)dg(y)\\
 & -\text{Re}\int_{\mathcal{R}}|\nabla^{h}u|_{h}^{2}(y)\left(\eta_{\epsilon}^{2}\Delta\chi_{\rho}^{2}+2\langle\nabla\chi_{\rho}^{2},\overline{\nabla}\eta_{\epsilon}^{2}\rangle\right)(y)K(x,y,1-t)dg(y)\\
 & -\int_{\mathcal{R}}|\nabla^{h}u|_{h}^{2}(y)\chi_{\rho}^{2}(y)\Delta\eta_{\epsilon}^{2}(y)K(x,y,1-t)dg(y)\\
\leq & C(n)\int_{\mathcal{R}}|\nabla^{h}u|_{h}^{2}\chi_{\rho}^{2}\eta_{\epsilon}^{2}(y)K(x,y,1-t)dg(y)\\
 & -\int_{\mathcal{R}}|\nabla^{h}\nabla^{h}u|_{h}^{2}(y)\chi_{\rho}^{2}\eta_{\epsilon}^{2}(y)K(x,y,1-t)dg(y)\\
 & +C\int_{\mathcal{R}}|\nabla^{h}\nabla^{h}u|_{h}(y)(|\nabla\chi_{\rho}|+\chi_{\rho}|\nabla\eta_{\epsilon}|)(y)(|\nabla^{h}u|_{h}\eta_{\epsilon}\chi_{\rho})(y)K(x,y,1-t)dg(y)\\
 & +C\int_{\mathcal{R}}|\nabla^{h}u|_{h}^{2}(y)\left(|\Delta\chi_{\rho}^{2}|(y)+\chi_{\rho}|\nabla\eta_{\epsilon}|\right)K(x,y,1-t)dg(y)\\
 & +C\int_{\mathcal{R}}|\nabla^{h}u|_{h}^{2}(y)\chi_{\rho}^{2}(y)\eta_{\epsilon}(y)|\nabla\eta_{\epsilon}|(y)|\nabla_{y}K(x,y,1-t)|dg(y).
\end{align*}
We use Cauchy's inequality to estimate
\begin{align*}
C\int_{\mathcal{R}}|\nabla^{h}\nabla^{h}u|_{h}(y) & (|\nabla\chi_{\rho}|+\chi_{\rho}|\nabla\eta_{\epsilon}|)(y)(|\nabla^{h}u|_{h}\eta_{\epsilon}\chi_{\rho})(y)K(x,y,1-t)dg(y)\\
\leq & \int_{\mathcal{R}}|\nabla^{h}\nabla^{h}u|_{h}^{2}(y)\chi_{\rho}^{2}\eta_{\epsilon}^{2}(y)K(x,y,1-t)dg(y)\\
 & +C^{2}\int_{\mathcal{R}}(|\nabla\chi_{\rho}|^{2}+\chi_{\rho}^{2}|\nabla\eta_{\epsilon}|^{2})(y)|\nabla^{h}u|_{h}^{2}(y)K(x,y,1-t)dg(y).
\end{align*}
The remaining terms can then be estimated as in the proof of $(i)$. 
\end{proof}
We now show how Corollary \ref{hormandersection} and Lemma \ref{c0elliptic}
can be used to perturb almost-holomorphic peak functions to holomorphic
peak functions.
\begin{prop}
\label{perturb} For any $\epsilon>0$, there exists $\zeta=\zeta(\epsilon)>0$
and $D_{0}=D_{0}(\epsilon)<\infty$ such that the following holds.

\noindent Suppose $x_{0}\in\mathcal{R}_{-1}$, and that there exists
$v\in C_{c}^{\infty}(\mathcal{R}_{-1},\mathbb{C})$ and an open subset
$U\subseteq\mathcal{R}_{-1}$ satisfying the following for some $D\in[D_{0},\infty)$:

$(i)$ $\int_{\mathcal{R}_{-1}}|v|^{2}d\nu_{-1}<1+\zeta,$ 

$(ii)$ $\sup_{U}\left|e^{-\frac{1}{2}d^{2}(x_{0},\cdot)}-|v|^{2}e^{-f_{-1}}\right|<\zeta$, 

$(iii)$ $\int_{\mathcal{R}_{-1}}|\overline{\partial}v|^{2}d\nu_{-1}<\zeta$,

$(iv)$ $\sup_{U}|\overline{\partial}v|^{2}e^{-f_{-1}}<\zeta$, 

$(v)$ $B(x_{0},2D)\cap\{r_{Rm}\geq\zeta\}\subseteq U$,

$(vi)$ $\text{supp}(v)\subseteq B(x_{0},D)$. 

\noindent Then there exists a holomorphic function $u\in L^{2}(\mathcal{R}_{-1},d\nu_{-1})$
satisfying the following:

$(a)$ $\int_{\mathcal{R}_{-1}}|u|^{2}d\nu_{-1}<1+\epsilon$,

$(b)$ $\sup_{\mathcal{R}_{-1}}\left|e^{-\frac{1}{2}d^{2}(x_{0},\cdot)}-|u|^{2}e^{-f_{-1}}\right|<\epsilon.$
\end{prop}

\begin{proof}
By Lemma \ref{hormandersection} and assumption $(iii)$, we can find
$\sigma\in L^{2}(\mathcal{R}_{-1},d\nu_{-1})$ satisfying $\overline{\partial}\sigma=\overline{\partial}v$
and 
\begin{equation}
\int_{\mathcal{R}_{-1}}|\sigma|_{h}^{2}dg=\int_{\mathcal{R}_{-1}}|\sigma|^{2}d\nu_{-1}\leq\int_{\mathcal{R}_{-1}}|\overline{\partial}v|^{2}d\nu_{-1}<\zeta.\label{l2normissmall}
\end{equation}
Thus $u:=v-\sigma\in L^{2}(\mathcal{R}_{-1},d\nu_{-1})$ is a holomorphic
function satisfying (using $(i)$)
\[
\int_{\mathcal{R}_{-1}}|u|^{2}d\nu_{-1}\leq2(1+\zeta).
\]
By $(iv)$, we know that for any $\epsilon_{0}\in(0,1]$ (to be determined)
and $x\in V:=\{r_{Rm}\geq2\epsilon_{0}\}$,
\[
\sup_{B(x,\epsilon)}|\overline{\partial}\sigma|_{h}^{2}=\sup_{B(x,\epsilon)}|\overline{\partial}v|_{h}^{2}<\zeta,
\]
so combining with (\ref{l2normissmall}) and local elliptic regularity
gives 
\[
\sup_{V}|\sigma|_{h}\leq\Psi(\zeta|\epsilon_{0}),
\]
and then combining with $(ii),(vi)$ gives 
\[
\sup_{V}\left|e^{-\frac{1}{2}d^{2}(x_{0},\cdot)}-|u|_{h}^{2}\right|<\Psi(\zeta|\epsilon_{0})+\Psi(D^{-1}).
\]
By Lemma \ref{therightway}, we know that $e^{-\frac{1}{2}d^{2}(x_{0},\cdot)},|u|_{h}^{2}$
are $C(Y)$-Lipschitz. By 
\[
\text{Vol}\left(\{r_{Rm}<10\epsilon_{0}\}\cap B(x,1)\right)\leq C(Y)\epsilon^{3},
\]
for any $x\in\mathcal{R}_{-1}$, we know that
\[
\text{Vol}\left(B(x,\epsilon_{0}^{\frac{2}{n}})\setminus(\{r_{Rm}<2\epsilon_{0}\}\cap B(x,1))\right)\geq c(Y)\epsilon_{0}^{2}-C(Y)\epsilon_{0}^{3}>0,
\]
hence by choosing $\epsilon_{0}\leq\overline{\epsilon}_{0}(Y)$, we
find $x'\in B(x,\epsilon_{0}^{\frac{2}{n}})$ satisfying $r_{Rm}(x)\geq2\epsilon_{0}$.
We can therefore integrate along minimizing geodesics to obtain
\[
\sup_{B(x_{0},\frac{3}{2}D)}\left|e^{-\frac{1}{2}d^{2}(x_{0},\cdot)}-|u|_{h}^{2}\right|<\epsilon_{0}^{\frac{2}{n}}C(Y)+\Psi(\zeta|\epsilon_{0})+\Psi(D^{-1}).
\]
Choosing $\epsilon_{0}\leq\overline{\epsilon}_{0}(\epsilon,Y)$, we
obtain the desired inequality when $D_{0}\geq\underline{D}(\epsilon,Y)$
and $\zeta\leq\overline{\zeta}(\epsilon,Y)$. 
\end{proof}
Next, we use the infinitesimal cone structure of $\mathcal{X}$ to
find functions satisfying the assumptions of Proposition \ref{perturb}.
\begin{lem}
\label{peaksections} (See Proposition 8 of \cite{CDSII}, Proposition
3.1 of \cite{szek1}) For any $x_{0}\in\mathcal{X}_{-1}$ and $\zeta>0$,
there exists $k\in\mathbb{N}$ with $k\geq\zeta^{-1}$ and some holomorphic
$u\in L^{2}(\mathcal{R}_{-\frac{1}{k}},d\nu_{-\frac{1}{k}})$ satisfying

$(i)$ $\sup_{\mathcal{R}_{-\frac{1}{k}}}\left||u|^{2}e^{-f_{-\frac{1}{k}}}-e^{-\frac{k^{2}}{2}d^{2}(x_{0}(-\frac{1}{k}),\cdot)}\right|<\zeta$,

$(ii)$ $\int_{\mathcal{R}_{-\frac{1}{k}}}|u|^{2}d\nu_{-\frac{1}{k}}<1+\zeta$.
\end{lem}

\begin{proof}
Suppose there exists $x_{0}\in\mathcal{X}_{-1}$ and $\zeta>0$ such
that no such $k\in\mathbb{N}$ exists. Because any tangent cone of
$(\mathcal{X}_{-1},d_{-1},x_{0})$ is a metric cone, we can find a
sequence $k_{j}\in\mathbb{N}$ such that $(\mathcal{X}_{-1},\sqrt{k_{j}}d_{-1},x_{0})$
converges in the pointed Gromov-Hausdorff sense to a metric cone $(C(Y),d_{C(Y)},y_{\ast})$,
where $y_{\ast}$ is the vertex. Consider the sequence

\[
(X_{j},d_{j},x_{j},\mathcal{R}_{j},\omega_{j},J_{j},L_{j},h_{j}):=\left(\mathcal{X}_{-\frac{1}{k_{j}}},\sqrt{k_{j}}d_{-\frac{1}{k_{j}}},x_{0}(-\frac{1}{k_{j}}),\mathcal{R}_{-\frac{1}{k_{j}}},k_{j}\omega_{-\frac{1}{k_{j}}},J_{-\frac{1}{k_{j}}},\mathcal{O}_{\mathcal{R}_{-\frac{1}{k_{j}}}},(4\pi)^{-n}e^{-f_{-\frac{1}{k_{j}}}}\right)
\]
of polarized Kähler singular spaces. Because $\mathcal{X}$ is static,
we moreover know that $(X_{j},d_{j},x_{j})$ are pointed isometric
to $(\mathcal{X}_{-1},\sqrt{k_{j}}d_{-1},x_{0})$, so that $(X_{j},d_{j},x_{j})\to(C(Y),d_{C(Y)},y_{\ast})$
in the pointed Gromov-Hausdorff sense. Moreover, we have smooth convergence
on the regular part, in the sense that there are diffeomorphisms $\psi_{j}$
as in $(i),(ii)$ of Definition \ref{polarizeddef}, where $(X,d,p)=(C(Y),d_{C(Y)},y_{\ast})$. 

\noindent \textbf{Claim: }After passing to a further subsequence,
there is a holomorphic Hermitian line bundle $(L,h)$ on $\mathcal{R}_{C(Y)}$
such that we have the following convergence of polarized singular
Kähler spaces as $j\to\infty$:
\[
(X_{j},d_{j},x_{j},\mathcal{R}_{j},\omega_{j},J_{j},L_{j},h_{j})\to(C(Y),d_{C(Y)},y_{\ast},\mathcal{R}_{C(Y)},g_{C(Y)},J_{C(Y)},L,h).
\]

It remains to construct $(L,h)$ such that item $(iii)$ of Definition
\ref{polarizeddef} holds. Many of the details are given in Proposition
3.1 of \cite{szek1}; our case is similar (but somewhat easier due
to more regularity) so we mostly summarize their arguments. Let $(\psi_{j})$
be as in $(i),(ii)$ of Definition \ref{polarizeddef}. For any compact
subset $K\subseteq\mathcal{R}_{C(Y)}$, we can choose a finite covering
by normal holomorphic coordinates; let $(B(q,4r),(z^{\alpha})_{\alpha=1}^{n})$
be one such chart, and note that we can (by taking $r\leq\overline{r}(\epsilon,K)$)
assume 
\[
\sup_{B(q,4r)}\left|\omega_{C(Y)}+\frac{1}{2}dJ_{C(Y)}d\left(\sum_{\alpha=1}^{n}|z^{\alpha}|^{2}\right)\right|_{\omega_{C(Y)}}\leq\epsilon,
\]
where $\epsilon>0$ is to be determined. Then a local Hörmander estimate
using the weight $(\sum_{\alpha=1}^{n}|z^{\alpha}|^{2})\circ\psi_{i}^{-1}$
on $B(\psi_{j}(q),2r)$ can be used to perturb the $(\psi_{j}^{-1})^{\ast}J_{C(Y)}$-holomorphic
coordinates to nearby $J_{j}$-holomorphic coordinates $(B_{g_{j}}(\psi_{j}(q),r),(z_{j}^{\alpha})_{\alpha=1}^{n})$
of $\mathcal{R}_{j}$, so that $\psi_{j}^{\ast}z_{j}^{\alpha}\to z^{\alpha}$
in $C_{\text{loc}}^{\infty}(B(q,2r))$ as $j\to\infty$. in particular,
we have
\[
\sup_{B_{g_{j}}(\psi_{j}(q),2r)}\left|\omega_{j}+\frac{1}{2}dJ_{j}d\left(\sum_{\alpha=1}^{n}|z^{\alpha}|^{2}\right)\right|_{\omega_{j}}\leq\Psi(\epsilon)
\]
for $j=j(\epsilon)\in\mathbb{N}$ sufficiently large. Noting that
the proof of the local $\partial\overline{\partial}$-Lemma shows
that a solution $\varphi$ of $\sqrt{-1}\partial\overline{\partial}\varphi=\alpha$
on a holomorphic coordinate ball can be made to have $C^{0}$ norm
bounded by $|\alpha|_{g_{\mathbb{C}^{n}}}$, we conclude the existence
of $s_{j}\in H^{0}(B_{g_{j}}(\psi_{j}(q),\frac{3}{2}r))$ satisfying
\[
\sup_{B_{g_{j}}(\psi_{j}(q),\frac{3}{2}r)}\left||s_{j}|_{h_{j}}^{2}-e^{-\sum_{\alpha=1}^{n}|z_{j}^{\alpha}|^{2}}\right|\leq\Psi(\epsilon).
\]
In particular, we can choose $\epsilon>0$ sufficiently small so that
$\frac{1}{2}\leq|s_{j}|_{h_{j}}^{2}\leq2$ on $B_{g_{j}}(\psi_{j}(q),\frac{3}{2}r)$.
Then the transition functions corresponding to these sections are
uniformly bounded, and bounded away from zero, so we can pass to subsequences
so that they converge to a cocycle defining a holomorphic line bundle
$L$ on a neighborhood of $K$ in $\mathcal{R}_{C(Y)}$. Moreover,
the representatives for $h_{j}$ with respect to these sections are
uniformly bounded as well, so via local elliptic regularity and the
polarization assumption, these also pass to the limit to produce a
Hermitian metric $h$ on $L$. More concretely, the bundle isomorphisms
$\xi_{j}$ can be defined using the sections $s_{j}$, limiting sections
in the line bundle $L$, and a partition of unity. $\square$

\noindent Next, recall that $C(Y)$ has singularities of codimension
4, so there exist cutoff functions $\eta_{\epsilon}\in C_{c}^{\infty}(\mathcal{R}_{C(Y)})$
satisfying $\text{supp}(\eta_{\epsilon})\subseteq\{r_{Rm}\geq\epsilon\}$
and $\eta_{\epsilon}|_{\{r_{Rm}>2\epsilon\}}\equiv1$

\noindent 
\[
\lim_{\epsilon\searrow0}\int_{\mathcal{R}_{C(Y)}\cap B(y_{\ast},r)}|\nabla\eta_{\epsilon}|^{3}\omega_{C(Y)}^{n}=0
\]
for each $r\in(0,\infty)$. Thus we can argue as in Proposition 8
of \cite{CDSII} to obtain (after replacing $k_{j}$ with $t_{j}k_{j}$
for some $t_{j}\in\mathbb{N}^{\times}$) $\widetilde{v}_{j}\in C_{c}^{\infty}(\mathcal{R}_{-\frac{1}{k_{j}}})$
and open subsets $U_{j}\subseteq\mathcal{R}_{-\frac{1}{k_{j}}}$ satisfying
the following for some sequences $\zeta_{j}\nearrow\infty$ and $D_{j}\nearrow\infty$: 

$(i)'$ $\int_{\mathcal{R}_{-\frac{1}{k_{j}}}}|\widetilde{v}_{j}|_{h_{j}}^{2}\omega_{j}^{n}=1+\zeta_{j},$

$(ii)'$ $\sup_{U_{j}}\left|e^{-\frac{1}{2}d_{g_{j}}^{2}(x_{j},\cdot)}-|\widetilde{v}_{j}|_{h_{j}}^{2}\right|<\zeta_{j}$, 

$(iii)'$ $\int_{\mathcal{R}_{-\frac{1}{k_{j}}}}|\overline{\partial}\widetilde{v}_{j}|_{\omega_{j}\otimes h_{j}}^{2}\omega_{j}^{n}<\zeta_{j}$,

$(iv)'$ $\sup_{U_{j}}|\overline{\partial}\widetilde{v}_{j}|_{\omega_{j}\otimes h_{j}}<\zeta_{j}$

$(v)'$ $B_{g_{j}}(x_{j},2D_{j})\cap\{r_{Rm}^{g_{j}}\geq\zeta'\}\subseteq U_{j}$,

$(vi)'$ $\text{supp}(\widetilde{v}_{j})\subseteq B_{g_{j}}(x_{j},D_{j})$.

\noindent Let $\Lambda_{j}:\mathcal{R}_{-\frac{1}{k_{j}}}\to\mathcal{R}_{-1}$
be the holomorphic map induced by the (backwards) flow of $\partial_{\mathfrak{t}}-\nabla f$,
which satisfies $\Lambda_{j}^{\ast}\omega_{-1}=k_{j}\omega_{-\frac{1}{k_{j}}}=\omega_{j}$,
$\Lambda_{j}^{\ast}f_{-1}=f_{-\frac{1}{k_{j}}}$, and extends to a
metric isometry $\Lambda_{j}:\mathcal{X}_{-\frac{1}{k_{j}}}\to\mathcal{X}_{-1}$
satisfying $\Lambda_{j}^{\ast}\left(d_{-1}(\Lambda(x_{j}),\cdot)\right)=k_{j}^{\frac{1}{2}}d_{-\frac{1}{k_{j}}}(x_{j},\cdot)=d_{g_{j}}(x_{j},\cdot)$.
Then $v_{j}:=(\Lambda_{j}^{-1})^{\ast}\widetilde{v}_{j}$ satisfies
the hypotheses of Lemma \ref{perturb} with $\zeta$ replaced by $\zeta_{j}$,
$D_{0}$ replaced by $D_{j}$, and $U$ replaced by $\Lambda_{j}(U_{j})$).
We thus obtain $u_{j}'\in L^{2}(\mathcal{R}_{-1},d\nu_{-1})$ holomorphic
satisfying 
\[
\int_{\mathcal{R}_{-1}}|u_{j}'|^{2}d\nu_{-1}<1+\Psi(j^{-1}),
\]
 
\[
\sup_{\mathcal{R}_{-1}}\left|e^{-\frac{1}{2}d_{-1}^{2}(x_{0},\cdot)}-|u_{j}'|^{2}e^{-f_{-1}}\right|<\Psi(j^{-1}).
\]
Then $u_{j}:=\Lambda_{j}^{\ast}u_{j}'$ is a holomorphic function
satisfying $(i),(ii)$ for sufficiently large $j\in\mathbb{N}$, a
contradiction. 
\end{proof}
\noindent Next, we show how to obtain holomorphic sections which separate
points and tangents. We identify any $h_{t_{0}}\in L^{2}(\mathcal{R}_{t_{0}},\nu_{t_{0}})$
with its corresponding holomorphic function $h:\mathcal{X}\to\mathbb{C}$. 
\begin{prop}
\label{separate} $(i)$ For any $x_{0},x_{1}\in\mathcal{X}_{-1}$,
there exist $t_{0}\in(-\infty,0)$ and $h\in L^{2}(\mathcal{R}_{t_{0}},\nu_{t_{0}})$
holomorphic such that $h(x_{0})\neq h(x_{1})$.

$(ii)$ For any $x_{0}\in\mathcal{R}_{-1}$, there exist $t_{0}\in(-\infty,0)$
and holomorphic $h_{1},...,h_{n}\in L^{2}(\mathcal{R}_{t_{0}},\nu_{t_{0}})$
whose restriction to a neighborhood of $x_{0}$ is an embedding. 
\end{prop}

\begin{proof}
$(i)$ This is a standard argument, given Lemma \ref{peaksections}.
By choosing a sequence $k_{j}\nearrow\infty$ such that both $(\mathcal{X}_{-1},k_{j}^{\frac{1}{2}}d_{-1},x_{0}),(\mathcal{X}_{-1},k_{j}^{\frac{1}{2}}d_{-1},x_{1})$
converge to tangent cones, and passing to further appropriate subsequences,
we can find $k\in\mathbb{N}$ arbitrarily large and holomorphic $u_{\ell}\in L^{2}(\mathcal{R}_{-\frac{1}{k}},d\nu_{-\frac{1}{k}})$
for $\ell=0,1$, such that 
\[
\sup_{\mathcal{R}_{-\frac{1}{k}}}\left||u_{\ell}|^{2}e^{-f_{-\frac{1}{k}}}-e^{-\frac{k^{2}}{2}d^{2}(x_{\ell}(-\frac{1}{k}),\cdot)}\right|<\frac{1}{100}.
\]
In particular, if we set $h:=e^{-f_{-\frac{1}{k}}}$, then 
\[
|u_{1}(x_{0})|_{h}^{2},|u_{0}(x_{1})|_{h}^{2}\leq\frac{1}{100},\qquad|u_{1}|_{h}^{2}(x_{1}),|u_{0}|_{h}^{2}(x_{0})\geq\frac{99}{100},
\]
so that $u:=u_{0}-\frac{u_{0}(x_{1})}{u_{1}(x_{1})}u_{1}$ satisfies
$u(x_{1})=0$ and 
\[
|u|_{h}(x_{0})\geq|u_{0}|_{h}(x_{0})-\frac{|u_{0}(x_{1})|_{h}}{|u_{1}(x_{1})|_{h}}|u_{1}(x_{0})|_{h}\geq\frac{1}{2}.
\]
In particular, $u(x_{0})\neq0$. 

$(ii)$ Choose holomorphic coordinates $(U,(z_{i})_{i=1}^{n})$ in
a neighborhood of $x_{0}$. Let $\chi\in C_{c}^{\infty}(U)$ be a
cutoff function such that $\chi|_{U'}\equiv1$ for some neighborhood
$U'\subset\subset U$ of $x_{0}$. Define
\[
\varphi_{\epsilon}:=B\chi(z)\log(|z|^{2}+\epsilon^{2})
\]
for $\epsilon\in[0,1]$, where $B>2n$. A routine computation gives
\[
\sqrt{-1}\partial\overline{\partial}\varphi_{\epsilon}\geq-CB\omega_{-1}
\]
for some $C<\infty$ independent of $\epsilon>0$. By choosing $\delta=\delta(B)>0$
sufficiently small, we therefore have
\[
\sqrt{-1}\partial\overline{\partial}(f_{-\delta}+\varphi_{\epsilon})=\delta^{-1}\omega_{-\delta}-CB\omega_{-1}=\left(\delta^{-1}-CB\right)\omega_{-\delta}.
\]
Moreover, we know $v_{j}:=\chi(z)z_{j}$ satisfy $\overline{\partial}v_{j}=0$
in $U'$, so that 
\[
\sup_{\epsilon\in(0,1]}\int_{\mathcal{R}_{-1}}|\overline{\partial}v_{j}|^{2}e^{-\varphi_{\epsilon}}d\nu_{-\delta}\leq C<\infty.
\]
For each $\epsilon>0$, Lemma \ref{lemmaforseparation} thus produces
$\sigma_{j,\epsilon}\in L^{2}(\mathcal{R}_{-\delta},e^{-\varphi_{\epsilon}}d\nu_{-\delta})$
such that $\overline{\partial}\sigma_{j,\epsilon}=\overline{\partial}v_{j}$
and 
\[
\int_{\mathcal{R}_{-1}}|\sigma_{j,\epsilon}|^{2}e^{-\varphi_{\epsilon}}d\nu_{-\delta}\leq C<\infty.
\]
Because $e^{-\varphi_{\epsilon}}$ is pointwise decreasing as a function
of $\epsilon$, we know that 
\[
\sup_{\epsilon\in(0,\epsilon_{0}]}\int_{\mathcal{R}_{-1}}|\sigma_{j,\epsilon}|^{2}e^{-\varphi_{\epsilon_{0}}}d\nu_{-\delta}\leq C<\infty
\]
for any fixed $\epsilon_{0}\in(0,1]$, where $C$ is independent of
$\epsilon_{0}$. We can therefore find a subsequence $\epsilon_{k}\searrow0$
such that $(\sigma_{j,\epsilon_{k}})_{k\in\mathbb{N}}$ converges
weakly in $L^{2}(\mathcal{R}_{-\delta},e^{-\varphi_{\epsilon_{0}}}d\nu_{-\delta})$
to some $\sigma_{j}\in L^{2}(\mathcal{R}_{-\delta},e^{-\varphi_{\epsilon_{0}}}d\nu_{-\delta})$
which satisfies $\overline{\partial}\sigma_{j}=\overline{\partial}v_{j}$
and 
\[
\int_{\mathcal{R}_{-1}}|\sigma_{j}|^{2}e^{-\varphi_{\epsilon}}d\nu_{-\delta}\leq C<\infty
\]
for any $\epsilon\in(0,1]$. By taking $\epsilon\searrow0$, the monotone
convergence theorem ensures
\[
\int_{\mathcal{R}_{-1}}|\sigma_{j}|^{2}e^{-\varphi_{0}}d\nu_{-\delta}<\infty.
\]
By the choice of $B$, this tells us that $\sigma_{j}(0)=0$ and $d\sigma_{j}(0)=0$,
so that $u_{j}:=v_{j}-\sigma_{j}$ are holomorphic functions satisfying
$u_{j}(0)=v_{j}(0)$ and $du_{j}(0)=dv_{j}(0)$ for $1\leq j\leq n$.
The claim then follows from the implicit function theorem.
\end{proof}
Let $x_{\ast}\in\mathcal{X}_{-1}$ correspond to the vertex of the
cone. The following is an easy consequence of Lemma \ref{peaksections}. 
\begin{lem}
(See Proposition 2.1 of \cite{donaldsun2}) There are $C$, $N\in\mathbb{N},$$0<\rho_{2}<\rho_{1}<1$
and holomorphic $h_{j}\in L^{2}(\mathcal{R}_{-1},d\nu_{-1})$ such
that $H:=(h_{1},...,h_{N})$ satisfies $\inf_{\partial B(x_{\ast},\rho_{1})}|H|>\frac{1}{2}$,
$\sup_{B(x_{\ast},\rho_{2})}|H|\leq\frac{1}{100}$, and $|dH|\leq C$. 
\end{lem}

\begin{proof}
This follows from Proposition \ref{separate}.
\end{proof}
Now define $\Omega:=H^{-1}(B(0^{N},\frac{1}{4}))\cap B(x_{\ast},\rho_{1})$
so that $H|_{\Omega}$ is proper. Set $F:=(f_{1},...,f_{N+Q}):=(h_{1},...,h_{N},u_{1},...,u_{Q})$,
where $u_{j}\in L^{2}(\mathcal{R}_{-1},d\nu_{-1})$ are holomorphic
functions such that $(u_{1},...,u_{Q})$ restricts to an embedding
of some nonempty open subset $D$ of $\Omega\cap\mathcal{R}_{-1}$
(this is possible by part $(ii)$ of Proposition \ref{separate}).
We may also assume (using part $(i)$ of Proposition \ref{separate},
properness, and compactness) that $F^{-1}(F(D))=D$. Note that $F|_{\Omega}$
is also proper. 
\begin{lem}
\label{techimage} $F(\Omega)$ is an irreducible complex analytic
subvariety of $B(0^{N},\frac{1}{4})\times\mathbb{C}^{Q}$ of dimension
$n$.
\end{lem}

\begin{proof}
Because $F$ is locally Lipschitz, we know that $F(\mathcal{S}(X)\cap\Omega)$
has Hausdorff dimension $\leq2(n-2)$, so that $F(\Omega)\setminus F(\mathcal{S}(X)\cap\Omega)\neq\emptyset$.
Choose $\eta_{\epsilon}\in C_{c}^{\infty}(\mathcal{R}_{-1})$ as in
Lemma \ref{LocLipschitz}. 

\noindent \textbf{Claim 1: }If we set $\langle T,\phi\rangle:=\int_{\Omega\cap\mathcal{R}_{-1}}F^{\ast}\phi$
for $\phi\in\mathcal{A}_{c}^{n,n}(B(0^{N},\frac{1}{4})\times\mathbb{C}^{Q})$,
then $T$ is a well-defined closed positive current.

We can write $\phi=\sum_{|I|=|J|=n}\phi_{I\overline{J}}dz^{I}\wedge d\overline{z}^{J}$
for some $\phi_{I\overline{J}}\in C_{c}^{\infty}(B(0^{N},\frac{1}{4})\times\mathbb{C}^{Q})$.
Because $|F^{\ast}dz^{I}|_{g}=|df^{I}|_{g}\leq C$ for all $I$, we
know that $|F^{\ast}\phi|_{g}\leq C||\phi||_{C^{0}(B(0^{N},c)\times\mathbb{C}^{Q})}$,
hence
\[
\int_{\Omega\cap\mathcal{R}}|F^{\ast}\phi|_{g}dg\leq C||\phi||_{C^{0}(B(0^{N},\frac{1}{4})\times\mathbb{C}^{Q})}\text{Vol}_{g}(\Omega)
\]
for all $\phi\in\mathcal{A}_{c}^{n,n}(B(0^{N},\frac{1}{4})\times\mathbb{C}^{Q})$.
This implies $T$ is a well-defined current, and can moreover be expressed
as
\[
\langle T,\phi\rangle=\lim_{\epsilon\searrow0}\langle T,\eta_{\epsilon}F^{\ast}\phi\rangle=\lim_{\epsilon\searrow0}\langle F_{\ast}(\eta_{\epsilon}T),\phi\rangle,
\]
which makes sense because $F|_{\text{supp}(\chi_{j})}$ is a proper
smooth function for each $j\in\mathbb{N}$. For any strongly positive
$\phi\in\mathcal{A}_{c}^{n,n}(B(0^{N},\frac{1}{4})\times\mathbb{C}^{Q})$,
because $F|_{\mathcal{R}}$ is holomorphic, we know $F^{\ast}\phi$
is also strongly positive, hence is a volume form; in particular,
we have 
\[
\langle T,\phi\rangle=\int_{\Omega}F^{\ast}\phi\geq0
\]
for all such $\phi$, so $T$ is in fact a positive current, and it
suffices to show that $T$ is closed. Suppose $\phi\in\mathcal{A}_{c}^{2n-1}(B(0^{N},\frac{1}{4})\times\mathbb{C}^{n})$,
so that $\eta_{\epsilon}F^{\ast}\phi\in C_{c}^{\infty}(\mathcal{R})$
for all $\epsilon\searrow0.$ Then 
\begin{align*}
|\langle T,d\phi\rangle|= & \left|\lim_{\epsilon\searrow0}\int_{\Omega\cap\mathcal{R}}\eta_{\epsilon}d(F^{\ast}\phi)\right|\leq\limsup_{\epsilon\searrow0}\int_{\Omega\cap\mathcal{R}}|F^{\ast}\phi|_{g}\cdot|\nabla\eta_{\epsilon}|d\mu_{g}\\
\leq & C\limsup_{\epsilon\searrow0}\int_{\Omega\cap\mathcal{R}}|\nabla\eta_{\epsilon}|d\mu_{g}=0.\hfill\hfill\square
\end{align*}

\noindent \textbf{Claim 2: }$T$ is a locally rectifiable current
of bidimension $(n,n)$. 

We recall (see Definition 2.1.4 in \cite{king}) that locally rectifiable
currents are exactly those which can be approximated in mass by Lipschitz
images of polyhedral chains. By taking such an approximation for the
locally rectifiable current (of finite mass) $\omega^{n}$ on $\mathcal{R}\cap\Omega$,
and composing with the Lipschitz function $F$, we obtain the desired
approximations for $T$. $\square$

By Claims 1 and 2, $T$ is a locally rectifiable closed positive current
of dimension $2n$, so Theorem 5.2.1 of \cite{king} states that there
exist $n$-dimensional analytic (irreducible) subvarieties $Z_{1},...,Z_{m}$
of $B(0^{N},\frac{1}{4})\times\mathbb{C}^{Q}$ and some $a_{1},...,a_{m}\in\mathbb{N}^{\times}$
such that 
\[
\int_{\Omega\cap\mathcal{R}}F^{\ast}\alpha=\langle T,\alpha\rangle=\sum_{j=1}^{m}a_{j}\int_{Z_{j}}\alpha
\]
for all $\alpha\in\mathcal{A}_{c}^{n,n}(B(0^{N},c)\times\mathbb{C}^{Q})$.
If $m\neq1$, then there exists $f\in\mathcal{O}(B(0^{N},c)\times\mathbb{C}^{Q})$
such that $f|_{Z_{j}}=0$ but $f|_{Z_{k}}\not\equiv0$ for $k\neq j$.
This means $F^{\ast}f\in\mathcal{O}(\Omega\cap\mathcal{R})$ vanishes
on the nonempty open set $F^{-1}(\Omega\setminus\bigcup_{k\neq j}Z_{k})$;
because $\mathcal{R}$ is connected, this implies $F^{\ast}f\equiv0$,
so that $F(\Omega\cap\mathcal{R})\subseteq Z_{j}$. Thus $a_{k}=0$
for $k\neq j$, a contradiction, hence $T=a[Z]$ for some $a\in\mathbb{N}^{\times}$,
and $Z=\text{supp}(T)=F(\Omega)$. Recalling that there is a subset
$D\subseteq\Omega\cap\mathcal{R}$ such that $F|_{D}$ is an embedding
and $F^{-1}(F(D))=D$, we moreover obtain $a=1$. 
\end{proof}
Next, we show that sections of $\mathcal{O}_{\mathcal{X}_{-1}}$ extend
to locally Lipschitz functions. 
\begin{lem}
\label{bddimplieslip} If $u:\Omega\cap\mathcal{R}_{-1}\to\mathbb{C}$
is holomorphic and locally bounded, then it extends to a locally Lipschitz
function $u:\Omega\to\mathbb{C}$. 
\end{lem}

\begin{proof}
We appeal to Corollary 6.3 of \cite{bam1}, once we verify that the
singular space $\mathcal{X}$ satisfies the ``tameness'' assumptions
used in its proof. Fix $x_{0}\in\Omega$. By Lemma \ref{peaksections},
there are $k_{j}\nearrow\infty$ and some holomorphic $s_{j}\in L^{2}(\mathcal{R}_{-1},e^{-k_{j}f_{-1}}d\nu_{-1})$
such that the following hold:
\[
\sup_{\mathcal{R}_{-1}}\left||s_{j}|^{2}e^{-k_{j}f_{-1}}-e^{-\frac{k_{j}}{2}d_{-1}^{2}(x_{0},\cdot)}\right|<\Psi(j^{-1}).
\]
Fix $j\in\mathbb{N}$ sufficiently large so that $B(x_{0},10\sqrt{k_{j}})\subseteq\Omega$
and 
\[
\sup_{\mathcal{R}_{-1}}\left||s_{j}|^{2}e^{-k_{j}f_{-1}}-e^{-\frac{k_{j}}{2}d_{-1}^{2}(x_{0},\cdot)}\right|<\frac{1}{100}.
\]
and write $s:=s_{j}$, $k:=k_{j}$, $h:=e^{-kf_{-1}}$, viewed as
a Hermitian metric on $\mathcal{O}_{\mathcal{R}_{-1}}$. Let $\chi\in C^{\infty}(\mathbb{R})$
be a cutoff function satisfying
\[
\chi|_{[-\frac{2}{3},\infty)}\equiv1,\qquad\chi|_{(-\infty,-1]}\equiv0,
\]
so that 
\[
\chi(\log|s|_{h}^{2})(x)=1
\]
for all $x\in\mathcal{R}_{-1}$ satisfying $|s|_{h}^{2}\geq e^{-\frac{2}{3}}$,
and 
\[
\chi(\log|s|_{h}^{2})(x)=0
\]
for $x\in\mathcal{R}_{-1}$ satisfying $|s|_{h}^{2}\leq e^{-1}$.
By the choice of $k$, we know that if $x\in B(x_{0},k^{-\frac{1}{2}})$,
then 
\[
|s|_{h}^{2}(x)\geq e^{-\frac{1}{2}}-\frac{1}{100}\geq e^{-\frac{2}{3}}.
\]
If $x\in\mathcal{R}_{-1}\setminus B(x_{0},2k^{-\frac{1}{2}})$, then
\[
|s|_{h}^{2}(x)\leq e^{-2}+\frac{1}{100}\leq e^{-1}.
\]
It follows that if $r_{0}:=k^{-\frac{1}{2}}$, then $\varphi:=\chi(\log|s|_{h}^{2})$
satisfies $\varphi|_{B(x_{0},r_{0})}\equiv1$ and $\text{supp}(\varphi)\subset\subset B(x_{0},2r_{0})$,
where we extend $\varphi$ to be zero on $\mathcal{X}_{-1}\setminus B(x_{0},2r_{0})$. 

\noindent \textbf{Claim 1: }$\varphi$ also satisfies
\[
\sup_{\mathcal{R}_{-1}}\left(r_{0}|\nabla\varphi|+r_{0}^{2}|\Delta\varphi|\right)\leq C(Y).
\]
In a neighborhood of $\text{supp}(\varphi)$, we have $|s|>0$, so
we can compute (in local holomorphic coordinates)
\[
\partial_{\overline{j}}\varphi=\chi'(\log|s|_{h}^{2})s^{-1}\overline{\nabla_{j}^{h}s},
\]
\[
\partial_{i}\partial_{\overline{j}}\varphi=\chi''(\log|s|_{h}^{2})|s|^{-2}\nabla_{i}^{h}s\overline{\nabla_{j}^{h}s}-k\chi'(\log|s|_{h}^{2})\partial_{i}\partial_{\overline{j}}f,
\]
so that (because $|s|\geq e^{-6}$ on $\text{supp}(\chi')$ and by
Lemma \ref{LocLipschitz} $|\nabla s|$ is locally bounded)
\[
|\nabla\varphi|^{2}\leq2(\chi'(\log|s|_{h}^{2}))^{2}\frac{|\nabla^{h}s|_{h}^{2}}{|s|_{h}^{2}}\leq\frac{C(Y)}{r_{0}^{2}},
\]
and similarly
\[
|\Delta\varphi|\leq2|\chi''(\log|s|_{h}^{2})|\frac{|\nabla^{h}s|_{h}^{2}}{|s|_{h}^{2}}+k|\chi'(\log|s|_{h}^{2})|\cdot|\Delta f_{-1}|\leq\frac{C(Y)}{r_{0}^{2}}.\,\,\,\,\,\,\square
\]

\noindent \textbf{Claim 2: }The heat kernel $K(x,y,t):=K(x;y(-1-t))$
for $x,y\in\mathcal{R}_{-1}$ and $t>0$ satisfies tameness assumption
$(5)$ of \cite{bam1}.

This is immediate from Lemmas \ref{heatkernel},\ref{fundsol},\ref{symmetry}.$\square$
\end{proof}
Set $Z:=F(\Omega)$. By adding finitely many locally bounded holomorphic
functions $u_{j}$ as components of $F$ to obtain $F'=(h_{1},...,h_{N},u_{1},...,u_{Q'})$,
Lemma \ref{bddimplieslip} and the proof of Lemma \ref{techimage}
guarantee that $Z:=F'(\Omega)$ is another analytic subvariety. If
$\pi:\mathbb{C}^{N+Q'}\to\mathbb{C}^{N+Q}$ is projection onto the
first $N+Q$ factors, we moreover obtain $\pi\circ F'=F$. 

Because the fibers of $\pi:Z'\to Z$ are compact subvarieties of $\mathbb{C}^{N+Q'}$,
they are singletons, so $\pi:Z'\to Z$ is a finite holomorphic map.
Moreover, $\pi|_{F'(D)}:F'(D)\to F(D)$ is a biholomorphism between
open subsets of $Z',Z$, so $\pi:Z'\to Z$ is generically one-to-one. 

Because $F$ is generically one-to-one, the same is true for $\pi:Z'\to Z$,
so we have an inclusion
\[
\mathcal{O}_{Z}\hookrightarrow\pi_{\ast}\mathcal{O}_{Z'}\hookrightarrow\widehat{\mathcal{O}_{Z}},
\]
where $\widehat{\mathcal{O}_{Z}}$ is the normalization of $\mathcal{O}_{Z}$.
As in \cite{donaldsun2}, the Noether property for coherent analytic
subsheaves of $\widehat{\mathcal{O}_{Z}}$ implies that we can choose
$Z'$ and $\pi$ so that $\pi_{\ast}\mathcal{O}_{Z'}$ is maximal
(after possibly shrinking $\Omega$) in the following sense: if $u$
is another locally bounded holomorphic function on $\Omega\cap\mathcal{R}$
(extended to a Lipschitz map on $\Omega$), and $F''=(F',u)$ is the
map obtained by adjoining $u$ to $F'$, with $Z'':=F''(\Omega)\subseteq\mathbb{C}^{N+Q'+1}$
the corresponding variety and $\pi':Z''\to Z'$ is the projection
map, then $(\pi')_{\ast}\mathcal{O}_{Z''}=\mathcal{O}_{Z'}$. By replacing
$Z'$ with $Z$, we may assume $Z=Z'$ and $\pi=\text{id}$. Let $F^{\ast}:\mathcal{O}(Z)\to\mathcal{O}(\Omega)$
be the ring homomorphism induced by the sheaf homomorphism $F^{\ast}:\mathcal{O}_{Z}\to F_{\ast}(\mathcal{O}|_{\Omega})$. 
\begin{prop}
\label{homeo} $F:\Omega\to Z$ is a homeomorphism mapping $\mathcal{R}\cap\Omega$
into the regular set of $Z$, and $F^{\ast}(\mathcal{O}(Z))=\mathcal{O}(\Omega)$. 
\end{prop}

\begin{proof}
\noindent \textbf{Claim 1: }$F^{\ast}(\mathcal{O}(Z))=\mathcal{O}(\Omega)$. 

Suppose $u\in\mathcal{O}(\Omega)$, and adjoin $u$ as another component
of $F$ to obtain $F':\Omega\to Z'$. The above discussion then shows
that $\pi:Z'\to Z$ satisfies $\pi_{\ast}\mathcal{O}_{Z'}=\mathcal{O}_{Z}$.
Moreover, $u=(F')^{\ast}\pi_{N+1}$, where $\pi_{N+1}$ is the restriction
to $Z'$ of projection onto the last component $\mathbb{C}^{N+1}\to\mathbb{C}$.
Choose $h\in\mathcal{O}(Z)$ such that $\pi^{\ast}h=\pi_{N+1}$. Because
$F=\pi\circ F'$, it follows that 
\[
F^{\ast}h=(F')^{\ast}\pi^{\ast}h=(F')^{\ast}\pi_{N+1}=u.\hfill\square
\]

\noindent \textbf{Claim 2: }$F:\Omega\to Z$ is a homeomorphism.

Because $F$ is proper, continuous, and surjective, it suffices to
show it is injective. Because $Z$ is a subvariety of an open subset
of $\mathbb{C}^{N+Q}$, we know $\mathcal{O}(Z)$ separates the points
of $Z$. Suppose $z_{1},z_{2}\in\Omega$ satisfy $F(z_{1})=F(z_{2})$,
so that every $h'\in\mathcal{O}(Z)$ satisfies $(h'\circ F)(z_{1})=(h'\circ F)(z_{2})$.
This means $h(z_{1})=h(z_{2})$ for every $h\in F^{\ast}(\mathcal{O}(Z))=\mathcal{O}(\Omega)$.
However, Proposition \ref{separate} tells us that $\mathcal{O}(\Omega)$
separates the points of $\Omega$, so that that $z_{1}=z_{2}$, hence
$F$ is injective. $\square$

It remains to show that for any $z\in\Omega\cap\mathcal{R}$, $F(z)$
is in the nonsingular part of $Z$. By Proposition \ref{separate},
we can find $h_{1},...,h_{n}\in\mathcal{O}(\Omega)$ such that $H:=(h_{1},...,h_{n})$
is a biholomorphism from a neighborhood $U$ of $z$ onto an open
subset $\widehat{U}$ of $\mathbb{C}^{n}$. By Claim 1, there are
$\widetilde{h}_{1},...,\widetilde{h}_{n}\in\mathcal{O}(Z)$ such that
$h_{j}=\widetilde{h}_{j}\circ F$ for $1\leq j\leq n$. Moreover,
Claim 2 implies that $\widetilde{U}:=F(U)$ is a neighborhood of $F(z)$
in $W$; we claim that $\widetilde{H}:=(\widetilde{h}_{1},...,\widetilde{h}_{n})$
is a biholomorphism to $\widehat{U}$. Because $F,H$ are homeomorphisms
and $H=\widetilde{H}\circ F$, it follows that $\widetilde{H}$ is
a homeomorphism. Moreover, Claim 1 gives that $H^{\ast}\mathcal{O}(\widehat{U})=F^{\ast}(\widetilde{H}^{\ast}\mathcal{O}(\widehat{U}))$,
so that the surjectivity of $H$ and injectivity of $F^{\ast}$ imply
$\widetilde{H}^{\ast}\mathcal{O}(\widehat{U})=\mathcal{O}(\widetilde{U})$.
Thus $\widetilde{H}^{-1}$ is biholomorphic. 
\end{proof}
To show that $Z$ is a normal variety, we argue as in \cite{donaldsun2},
with some modifications due to the fact that the Hörmander technique
cannot be done on the smooth approximants $M_{i}$.
\begin{prop}
After further shrinking $\Omega$, we can assume $Z$ is normal and
$F^{\ast}:\mathcal{O}_{Z}\to F_{\ast}(\mathcal{O}|_{\Omega})$ is
bijective. 
\end{prop}

\begin{proof}
By the openness of the normal locus, it suffices to show that $Z$
is normal in a neighborhood of the vertex $F(x_{\ast})$. 

\noindent \textbf{Claim: }By possibly adding components to $F$ and
shrinking $\Omega$, we may assume that $F^{\ast}:\mathcal{O}_{Z,F(x_{\ast})}\to\mathcal{O}_{x_{\ast}}$
is surjective.

By adding locally bounded functions near $x_{\ast}$ whose restriction
to $\mathcal{R}$ is holomorphic as components of $F$, we obtain
a map $F'$ such that $Z':=F'(\Omega)$ is an analytic variety, and
a projection $\pi:Z'\to Z$ such that $F=\pi\circ F'$. Arguing as
in Proposition \ref{homeo}, but now simply using the fact that $\mathcal{O}_{Z,F(x_{\ast})}$
is a Noetherian ring, we can find $F'$ such that $\pi_{\ast}(\mathcal{O}_{Z',F'(x_{\ast})})\subseteq\widehat{\mathcal{O}}_{Z,F(x_{\ast})}$
is maximal among all such $(Z',\pi)$. Again, we can replace $Z$
with $Z'$ to assume that $Z$ already has this property. Now argue
as in Proposition \ref{homeo} again. $\square$

Suppose $h\in\mathcal{O}_{Z}(U\cap Z_{\text{reg}})$ is bounded in
a neighborhood $U$ of $F(x_{\ast})$. Then $F^{\ast}h\in\mathcal{O}(F^{-1}(U)\cap\mathcal{R})$
is bounded, and $F^{-1}(U)$ is a neighborhood of $x_{\ast}$. Note
that $F^{\ast}h$ extends uniquely to a continuous function on $F^{-1}(U)$,
so $h$ extends to a continuous function $h:U\to\mathbb{C}$; it remains
to check that $h$ is holomorphic. We have $F^{\ast}h\in\mathcal{O}_{x_{\ast}}=F^{\ast}(\mathcal{O}_{Z,F(x_{\ast})})$,
so the Claim gives $\widetilde{h}\in\mathcal{O}_{Z,F(x_{\ast})}$
such that $F^{\ast}\widetilde{h}=F^{\ast}h$ as elements of $\mathcal{O}_{x_{\ast}}$.
This implies (by the connectedness of $\mathcal{R}$) that there exists
$\widehat{h}\in\mathcal{O}_{Z}(V)$ for some neighborhood $V\subseteq U$
of $F(x_{\ast})$ in $Z$, so that the germ of $\widehat{h}$ at $F(x_{\ast})$
is $\widetilde{h}$, and $F^{\ast}\widehat{h}=F^{\ast}h$ on $F^{-1}(V)\cap\mathcal{R}$.
Because $F$ is a homeomorphism, and because $h,\widehat{h}$ are
continuous on $V$, we conclude that $h=\widehat{h}$, so that $\widehat{h}$
is the desired extension, and $Z$ is normal at $F(x_{\ast})$. 

It remains to show that $F^{\ast}:\mathcal{O}_{Z}\to F_{\ast}(\mathcal{O}|_{\Omega})$
is surjective. Fix $q\in\Omega$, and a holomorphic function $u$
defined on a neighborhood $U$ of $q$. This determines a holomorphic
map $F'=(F,u):U\to\mathbb{C}^{N+1}$, and $\pi:F'(U)\to F(U)$ is
generically one-to-one. Moreover, we may identify $u$ as a function
on $F'(U)$ given by projection onto the last factor. Because $Z$
is normal, $\pi$ must be a holomorphic equivalence, hence $u$ is
holomorphic on $F(U)$.
\end{proof}
\begin{rem}
Note that unlike Proposition 2.4 in \cite{donaldsun2}, we are not
claiming that sections of $\mathcal{O}_{\mathcal{X}_{-1}}$ can be
approximated by holomorphic functions on smooth manifolds. 
\end{rem}

\section{Algebraic Structure and the Regular Set}

In order to prove the second half of Theorem \ref{thm1}, we adapt
ideas of \cite{szek2} to the parabolic setting using the parabolic
frequency function for solutions of the heat equation coupled to Ricci
flow, previously studied in \cite{julius}. The following roughly
states that a heat flow with parabolic frequency close to $1$ must
be an almost-splitting function at a smaller scale.

Throughout this section, we set $N:=2^{10}$. 
\begin{lem}
\label{frequency} Suppose $\delta>0$, $B<\infty$, $(M^{n},(g_{t})_{t\in I})$
is a closed Ricci flow, and $y\in C^{\infty}(M\times[-N,-\delta])$
is a solution of the heat equation such that the following hold for
some $x_{0}\in M$:

$(i)$ $(x_{0},0)$ is $(\delta,1)$-selfsimilar,

$(ii)$ $\int_{M}yd\nu_{x_{0},0;t}=0$,

$(iii)$ $\int_{M}y^{2}d\nu_{x_{0},0;-N}\leq B$,

$(iv)$ $a:=\frac{1}{2}\int_{M}y^{2}d\nu_{x_{0},0;-1}\geq\frac{1}{4}$,

$(v)$ $\int_{M}y^{2}d\nu_{x_{0},0;-2}\leq2^{1+\delta}\int_{M}y^{2}d\nu_{x_{0},0;-1}$.

\noindent Then $\widehat{y}:=a^{-\frac{1}{2}}y$ satisfies 
\[
\int_{-1}^{-\epsilon}\int_{M}|\nabla^{2}\widehat{y}|^{2}d\nu_{t}dt+\int_{-1}^{-\epsilon}\int_{M}\left||\nabla\widehat{y}|^{2}-1\right|d\nu_{t}dt\leq\Psi(\delta|\epsilon,B).
\]
\end{lem}

\begin{proof}
Write $\nu_{t}:=\nu_{x_{0},0;t}$ and let $U(\tau):=\frac{2\tau\int_{M}|\nabla y|^{2}d\nu_{t}}{\int_{M}y^{2}d\nu_{t}}$
denote the parabolic frequency function of $y$ based at $(x_{0},0)$.
We first observe that 
\[
\frac{d}{dt}\log\int_{M}y^{2}d\nu_{t}=-\frac{2\int_{M}|\nabla y|^{2}d\nu_{t}}{\int_{M}y^{2}d\nu_{t}}=-\frac{1}{\tau}U(\tau).
\]
Moreover, we can use the computation in Section 1 of \cite{julius}
to estimate
\begin{align*}
\frac{d}{dt}U(\tau)= & -2\left(\frac{\int_{M}(\Delta_{f}y)^{2}d\nu_{t}\int_{M}y^{2}d\nu_{t}-\left(\int_{M}y\Delta_{f}yd\nu_{t}\right)^{2}}{\int_{M}|\nabla y|^{2}d\nu_{t}\int_{M}y^{2}d\nu_{t}}\right)\\
 & +\frac{2\int_{M}(Rc+\nabla^{2}f-\frac{1}{2}g)(\nabla y,\nabla y)d\nu_{t}}{\int_{M}y^{2}d\nu_{t}}\\
\leq & \frac{2\int_{M}|Rc+\nabla^{2}f-\frac{1}{2}g|\cdot|\nabla y|^{2}d\nu_{t}}{\int_{M}y^{2}d\nu_{t}}.
\end{align*}
For all $t\in[-2,-1]$, we have $\int_{M}|y|^{2}d\nu_{t}\geq\int_{M}|y|^{2}d\nu_{-1}\geq\frac{1}{4}$,
so Theorem 12.1 of \cite{bamlergen1} allows us to estimate
\begin{align*}
\int_{M}|\nabla y|^{4}d\nu_{t}\leq & \int_{M}|\nabla y|^{2}d\nu_{-N+1}\leq\int_{-N}^{-N+1}\int_{M}|\nabla y|^{2}d\nu_{t}dt\\
\leq & \int_{M}y^{2}d\nu_{-N}\leq B.
\end{align*}
It follows that $\frac{d}{dt}U(\tau)\leq4\sqrt{B\delta}$, so for
$\tau\in[1,2]$, we have $U(1)\leq U(\tau)+4\sqrt{B\delta}$. Integrating
from $t=-2$ to $t=-1$ gives 
\begin{align*}
-(1+\delta)\log2\leq & \log\left(\frac{\int_{M}y^{2}d\nu_{-1}}{\int_{M}y^{2}d\nu_{-2}}\right)=-\int_{-2}^{-1}\frac{1}{\tau}U(\tau)dt\\
\leq & -U(1)\int_{-2}^{-1}\frac{1}{\tau}dt+4\sqrt{B\delta}\\
= & -(\log2)U(1)+4\sqrt{B\delta}
\end{align*}
which implies $U(1)\leq1+\delta+8\sqrt{B\delta}$.

\noindent \textbf{Claim: }$U(\tau)\leq1+\Psi(\delta|\epsilon,B)$
for all $\tau\in[\epsilon,1]$. 

From $U(1)\leq1+C(B)\sqrt{\delta}$, we get
\[
2\int_{M}|\nabla y|^{2}d\nu_{-1}\leq\left(1+C(B)\sqrt{\delta}\right)\int_{M}y^{2}d\nu_{-1},
\]
hence 
\[
\frac{d}{dt}\int_{M}y^{2}d\nu_{t}=-2\int_{M}|\nabla y|^{2}d\nu_{t}\geq-2\int_{M}|\nabla y|^{2}d\nu_{-1}\geq-\left(1+C(B)\sqrt{\delta}\right)\int_{M}y^{2}d\nu_{-1}.
\]
Integrating from $t=-1$ to $t=-\epsilon$ yields
\[
\int_{M}y^{2}d\nu_{-\epsilon}-\int_{M}y^{2}d\nu_{-1}\geq-(1-\epsilon)\left(1+C(B)\sqrt{\delta}\right)\int_{M}y^{2}d\nu_{-1},
\]
so that 
\[
\int_{M}y^{2}d\nu_{-\epsilon}\geq(\epsilon-\Psi(\delta|B))\int_{M}y^{2}d\nu_{-1}\geq\frac{1}{2}\epsilon
\]
if $\delta\leq\overline{\delta}(\epsilon,B)$. For any $\tau\in[\epsilon,1]$,
we can thus estimate
\begin{align*}
U(\tau)-U(1)\leq & 2\epsilon^{-1}\left(\int_{-1}^{-\epsilon}\int_{M}|Rc+\nabla^{2}f-\frac{1}{2}g|^{2}d\nu_{t}dt\right)^{\frac{1}{2}}\left(\int_{-1}^{-\epsilon}\int_{M}|\nabla y|^{4}d\nu_{t}dt\right)^{\frac{1}{2}}\\
\leq & \Psi(\delta|\epsilon,B).\,\,\,\,\,\,\,\,\,\,\,\,\,\,\,\,\,\,\,\,\,\,\,\,\,\,\,\,\,\,\,\,\,\,\,\,\,\,\,\,\,\,\,\,\,\,\,\,\,\,\,\,\,\,\,\,\,\,\,\,\,\,\,\,\,\,\,\,\,\,\,\,\,\,\,\,\,\,\,\,\,\,\,\,\,\,\,\,\,\,\,\,\,\,\,\,\,\,\,\,\,\,\,\,\,\,\,\,\,\,\,\,\,\,\,\,\,\,\,\,\,\,\square
\end{align*}

Combining the Claim with the $L^{2}$ Poincare inequality (Theorem
1.10 of \cite{hein}) gives
\[
\int_{M}y^{2}d\nu_{t}\leq2\tau\int_{M}|\nabla y|^{2}d\nu_{t}\leq(1+\Psi(\delta|\epsilon,B))\int_{M}y^{2}d\nu_{t}
\]
for all $\tau\in[\epsilon,1]$. We can therefore integrate
\[
\frac{d}{dt}\left(2\tau\int_{M}|\nabla y|^{2}d\nu_{t}-\int_{M}y^{2}d\nu_{t}\right)=-4\tau\int_{M}|\nabla^{2}y|^{2}d\nu_{t}
\]
from $t=-1$ to $t=-\epsilon$ to obtain
\[
\int_{-1}^{-\epsilon}\int_{M}\tau|\nabla^{2}y|^{2}d\nu_{t}dt\leq\Psi(\delta|\epsilon,B).
\]
Recalling that $a^{-\frac{1}{2}}\leq2$, we thus have$\int_{-1}^{-\epsilon}\int_{M}|\nabla^{2}\widehat{y}|^{2}d\nu_{t}dt\leq\Psi(\delta|\epsilon,B)$.
Since $1\leq U(1)\leq1+\Psi(\delta|B)$ implies
\[
1\leq\int_{M}|\nabla\widehat{y}|^{2}d\nu_{-1}\leq1+\Psi(\delta|B),
\]
we can use $\left.\int_{M}|\nabla\widehat{y}|^{2}d\nu_{t}\right|_{t=t_{1}}^{t=t_{2}}\leq\Psi(\delta|B,\epsilon)$
for all $t_{1},t_{2}\in[\epsilon,1]$ and the $L^{2}$-Poincare inequality
to estimate
\begin{align*}
\int_{-1}^{-\epsilon}\int_{M}\left||\nabla\widehat{y}|^{2}-1\right|d\nu_{t}dt\leq & 2\left(\int_{-1}^{-\epsilon}\int_{M}|\nabla^{2}\widehat{y}|^{2}d\nu_{t}dt\right)^{\frac{1}{2}}\left(\int_{-1}^{-\epsilon}\int_{M}|\nabla\widehat{y}|^{2}d\nu_{t}dt\right)^{\frac{1}{2}}+\int_{-1}^{-\epsilon}\left|\int_{M}(|\nabla\widehat{y}|^{2}-1)d\nu_{t}\right|dt\\
\leq & \Psi(\delta|B,\epsilon).
\end{align*}
\end{proof}
We will now observe some technical results which will be necessary
for proving that regular points of $\mathcal{X}$ as an analytic variety
are also regular in the metric sense. The following lemma is a consequence
of heat kernel estimates proved in \cite{bamlergen1,bamlergen3}.
\begin{lem}
\label{goodgrad} $(i)$ For any $A<\infty$ and $\alpha>0$ there
exists $C=C(A,Y,\alpha)<\infty$ such that the following holds for
any $t_{0}\in(-\infty,0)$, $x_{0}\in\mathcal{X}_{t_{0}}$, and $r>0$.
For any $x\in B(x_{0},Ar)$, we have 
\[
d\nu_{x;t_{0}-r^{2}}\leq Ce^{\alpha f_{x_{0};t_{0}-r^{2}}}d\nu_{x_{0};t_{0}-r^{2}},
\]
where $f_{x_{0};t_{0}-r^{2}}\in C^{\infty}(\mathcal{R}_{t_{0}-r^{2}})\cap C(\mathcal{R}_{t_{0}-r^{2}})$
is defined by $d\nu_{x_{0};t_{0}-r^{2}}=(4\pi r^{2})^{-n}e^{-f_{x_{0};t_{0}-r^{2}}}dg_{t_{0}-r^{2}}.$

$(ii)$ For any $p\in[1,\infty)$, there exists $C=C(n,p)<\infty$
such that for any $x\in\mathcal{X}_{t_{0}}$ and $t\in(-\infty,t_{0})$,
we have
\[
\int_{\mathcal{R}_{t}}|\nabla_{y}\log K(x;y)|^{p}d\nu_{x;t}(y)=\int_{\mathcal{R}_{t}}|\nabla_{x}\log K(x;y)|^{p}d\nu_{x;t}(y)\leq\frac{C}{(t_{0}-t)^{p}}.
\]
\end{lem}

\begin{proof}
$(i)$ By parabolic rescaling, we may assume $r=1$. Suppose $z\in B(x_{0};A)\cap\mathcal{R}_{t_{0}}$.
By Theorem 6.19 of \cite{bamlergen2}, there exist $z_{i}\in M_{i}$
and $x_{0,i}\in M_{i}$ satisfying 
\[
(z_{i},t_{0})\xrightarrow[i\to\infty]{\mathfrak{C}}z,\hfill(x_{0,i},t_{0})\xrightarrow[i\to\infty]{\mathfrak{C}}x_{0}.
\]
Because $z\in P^{\ast}(x_{0};A,t_{0}-\theta,t_{0}+\theta)$ for any
$\theta\in(0,\frac{1}{2})$, Lemma 15.8(b) of \cite{bamlergen3} gives
\[
(z_{i},t_{0})\in P^{\ast}(x_{0,i},t_{0};A,t_{0}-\theta,t_{0}+\theta)
\]
for sufficiently large $i=i(\theta)\in\mathbb{N}$, so that 
\[
d_{W_{1}}^{g_{i,t_{0}-\frac{1}{2}}}(\nu_{z_{i},t_{0};t_{0}-\theta}^{i},\nu_{x_{0,i},t_{0};t_{0}-\theta}^{i})\leq A.
\]
We then apply Proposition 8.1 of \cite{bamlergen3} with $t_{0},t_{1}$
in that proposition replaced by $t_{0}$, with $t^{\ast}$ replaced
by $t_{0}-\theta$, and $s$ replaced by $t_{0}-1$, obtaining
\[
d\nu_{z_{i},t_{0};t_{0}-1}^{i}\leq C(A,Y)e^{\frac{1}{8}f_{x_{0,i},t_{0};t_{0}-1}}d\nu_{x_{0,i},t_{0};t_{0}-1}^{i},
\]
where $f_{x_{0,i},t_{0};t_{0}-1}$ is defined by $d\nu_{x_{0,i},t_{0};t_{0}-1}^{i}=(4\pi)^{\frac{n}{2}}e^{-f_{x_{0,i},t_{0};t_{0}-1}}dg_{i,t_{0}-1}$,
if we choose $\theta\leq\overline{\theta}(A)$. By Remark \ref{convergenotation},
we moreover know that for sufficiently large $i\in\mathbb{N}$, we
have $x_{0,i},z_{i}\in\psi_{i,t_{0}}(U_{i}\cap\mathcal{R}_{t_{0}})$,
and $\psi_{i,t_{0}}^{-1}(x_{0,i})\to x_{0}$, $\psi_{i,t_{0}}^{-1}(z_{i})\to z$
as $i\to\infty$. Again using the convergence in Remark \ref{convergenotation},
we conclude that 
\[
d\nu_{z;t_{0}-1}\leq C(A,Y)e^{\frac{1}{8}f_{x_{0};t_{0}-1}}d\nu_{x_{0};t_{0}-1}.
\]
$(ii)$ This follows from Remark \ref{convergenotation}, Proposition
4.2 of \cite{bamlergen1}, and Proposition \ref{conjheatgradient}.
\end{proof}
Next, we will construct cutoff functions with improved estimates using
the peak holomorphic functions constructed in Section 4.

Suppose $(\mathcal{X},(\mu_{t})_{t\in(-\infty,0)})$ is a static metric
cone which is an $\mathbb{F}$-limit of smooth Kähler-Ricci flows
as in (\ref{F}). Let $\chi_{\rho}\in C_{c}(B(x_{\ast},2\rho))$ be
as in Lemma \ref{LocLipschitz}, so that $\chi_{\rho}|_{B(x_{\ast},\rho)}\equiv1$,
$\partial_{\mathfrak{t}}\chi_{\rho}\equiv0$, and 
\[
|\nabla\chi_{\rho}|^{2}+|\Delta\chi_{\rho}|\leq\frac{C(Y)}{\rho^{2}}
\]
on $\mathcal{R}$. 
\begin{lem}
There exist cutoff functions $\eta_{\epsilon}\in C^{\infty}(\mathcal{R})$
such that the following hold for each $\rho\in(0,\infty)$:

$(i)$ $\text{supp}(\chi_{\rho}\eta_{\epsilon})\subset\subset\mathcal{R}$,

$(ii)$ $\text{supp}(\chi_{\rho}(1-\eta_{\epsilon}))\subseteq\{x\in\mathcal{X};d(x,\mathcal{S})<\epsilon\}$,

$(iii)$ $\partial_{\mathfrak{t}}\eta_{\epsilon}\equiv0$,

$(iv)$ $\lim_{\epsilon\searrow0}\int_{\mathcal{R}\cap B(x_{\ast},2\rho)}|\nabla\eta_{\epsilon}|^{\frac{7}{2}}dg=\lim_{\epsilon\searrow0}\int_{\mathcal{R}\cap B(x_{\ast},2\rho)}|\Delta\eta_{\epsilon}|^{\frac{7}{4}}dg=0$. 
\end{lem}

\begin{proof}
By the proof of Proposition \ref{bddimplieslip} and a Vitali covering
lemma, we can find a cover $\{B(x_{i},r_{i})\}_{i\in\mathcal{I}}$
of $\mathcal{S}_{-1}:=\mathcal{S}\cap\mathcal{X}_{-1}$ with $x_{i}\in\mathcal{S}_{-1}$
such that the following hold:

$(a)$ For any $i\in\mathcal{I}$, there exists $\varphi_{i}\in C_{c}(B(x_{i},2r_{i}))\cap C^{\infty}(\mathcal{R}\cap B(x_{i},2r_{i}))$
such that 
\[
r_{i}|\nabla\varphi_{i}|+r_{i}^{2}|\Delta\varphi_{i}|\leq C(Y),
\]

$(b)$ For distinct $i,j\in\mathcal{I}$, we have $B(x_{i},\frac{r_{i}}{5})\cap B(x_{j},\frac{r_{j}}{5})=\emptyset$.

\noindent For each $m\in\mathbb{N}$ with $2^{-m}<\epsilon$, we define
\[
\mathcal{I}_{m,\rho}:=\{i\in\mathcal{I};2^{-m}<r_{i}<2^{-m+1}\text{ and }x_{i}\in B(x_{\ast},4\rho)\}.
\]
Then for any $p\in(0,4)$, using $(b)$ and the fact that $\mathcal{X}_{-1}$
has singularities of codimension four,
\begin{align*}
|\mathcal{I}_{m,\rho}|\leq & C(Y)2^{2mn}\sum_{i\in\mathcal{I}_{m,\rho}}\text{Vol}\left(B(x_{i},\frac{1}{5}r_{i})\cap\mathcal{R}\right)\\
\leq & C(Y)2^{2mn}\text{Vol}\left(\{x\in B(x_{\ast};8\rho)\cap\mathcal{R};d(x,\mathcal{S}_{-1})<2^{-m}\}\right)\\
\leq & C(Y,\rho,p)2^{m(2n-p)}.
\end{align*}
Combining this with $(a)$ yields
\begin{align*}
\int_{B(x_{\ast},2\rho)\cap\mathcal{R}}|\nabla\eta_{\epsilon}|^{\frac{7}{2}}dg\leq & C(Y)\sum_{m\in\mathbb{N}:2^{-m}\leq\epsilon}\sum_{i\in\mathcal{I}_{m,\rho}}r_{i}^{-\frac{7}{2}}\text{Vol}(B(x_{i},2r_{i})\cap\mathcal{R})\\
= & C(Y,\rho,p)\sum_{m\in\mathbb{N}:2^{-m}\leq\epsilon}2^{m(2n+\frac{7}{2}-p)}2^{-2mn}\\
\leq & C(Y,\rho,p)\epsilon^{\frac{1}{2}(p-\frac{7}{2})}\sum_{m=1}^{\infty}2^{-\frac{m}{2}(p-\frac{7}{2})}\\
\le & C(Y,\rho,p)\epsilon^{\frac{1}{2}(p-\frac{7}{2})}
\end{align*}
if we choose $p\in(\frac{7}{2},4)$. Similarly, we have 
\begin{align*}
\int_{B(x_{\ast},2\rho)\cap\mathcal{R}}|\Delta\eta_{\epsilon}|^{\frac{7}{4}}dg\leq & C(Y)\sum_{m\in\mathbb{N}:2^{-m}\leq\epsilon}\sum_{i\in\mathcal{I}_{m,\rho}}r_{i}^{-\frac{7}{2}}\text{Vol}(B(x_{i},2r_{i})\cap\mathcal{R})\\
\leq & C(Y,\rho,p).
\end{align*}
\end{proof}
The following will be used to estimate the error to solving the heat
equation introduced when multiplying $u$ by appropriate cutoff functions.
\begin{lem}
\label{L2gradient} If $x\in\mathcal{X}_{t_{\ast}}$ and $u\in C^{\infty}(\mathcal{R})\cap C(\mathcal{X})$
is a holomorphic function satisfying
\[
\int_{\mathcal{R}_{t_{0}}}|u|^{2}d\nu_{x;t_{0}}<\infty
\]
for some $t_{0}<t_{\ast}$, then
\[
\int_{\mathcal{R}_{t_{1}}}|u|^{2}d\nu_{x;t_{1}}+\int_{t_{0}}^{t_{1}}\int_{\mathcal{R}_{t}}|\nabla u|^{2}d\nu_{x;t}dt\leq\int_{\mathcal{R}_{t_{0}}}|u|^{2}d\nu_{x;t_{0}}
\]
for $t_{1}\in(t_{0},t_{\ast})$. 
\end{lem}

\begin{proof}
For ease of notation, write $\nu_{t}':=\nu_{x;t}$. Note that $\int_{\mathcal{R}_{t}}|u|^{2}d\nu_{x;t}<\int_{\mathcal{R}_{t_{0}}}|u|^{2}d\nu_{x;t_{0}}$
for all $t\in(t_{0},0)$, which can be seen by integrating
\[
|u|^{2}(y)=\left|\int_{\mathcal{R}_{t_{0}}}ud\nu_{y;t_{0}}\right|^{2}\leq\int_{\mathcal{R}_{t_{0}}}|u|^{2}d\nu_{y;t_{0}}
\]
for $y\in\mathcal{R}_{t}$ against $d\nu_{t}$':
\[
\int_{\mathcal{R}_{t}}|u|^{2}d\nu_{t}'\leq\int_{\mathcal{R}_{t}}\left(\int_{\mathcal{R}_{t_{0}}}|u|^{2}d\nu_{y;t_{0}}\right)d\nu_{t}'(y)=\int_{\mathcal{R}_{t_{0}}}|u|^{2}d\nu_{t_{0}}'.
\]
For $\sigma>0$ to be determined, we estimate
\begin{align*}
\frac{d}{dt}\int_{\mathcal{R}_{t}}\eta_{\epsilon}^{2}\chi_{\rho}^{2}|u|^{2}d\nu_{t}'\leq & -\int_{\mathcal{R}_{t}}\eta_{\epsilon}^{2}\chi_{\rho}^{2}|\nabla u|^{2}d\nu_{t}'+8\int_{\mathcal{R}_{t}}\left(\eta_{\epsilon}|\nabla\chi_{\rho}|+\chi_{\rho}|\nabla\eta_{\epsilon}|\right)|u|\cdot\eta_{\epsilon}\chi_{\rho}|\nabla u|d\nu_{t}'\\
 & +8\int_{\mathcal{R}_{t}}|\nabla\chi_{\rho}|\cdot|\nabla\eta_{\epsilon}|\eta_{\epsilon}\chi_{\rho}|u|^{2}d\nu_{t}'+8\int_{\mathcal{R}_{t}}\left(|\nabla\eta_{\epsilon}|^{2}+\eta_{\epsilon}|\Delta\eta_{\epsilon}|\right)\chi_{\rho}^{2}|u|^{2}d\nu_{t}'\\
 & +8\int_{\mathcal{R}_{t}}\left(|\nabla\chi_{\rho}|^{2}+\chi_{\rho}|\Delta\chi_{\rho}|\right)|u|^{2}d\nu_{t}'\\
\leq & -(1-\sigma)\int_{\mathcal{R}_{t}}\eta_{\epsilon}^{2}\chi_{\rho}^{2}|\nabla u|^{2}d\nu_{t}'+\frac{C(\sigma)}{\rho^{2}}\int_{\mathcal{R}_{t}}|u|^{2}d\nu_{t}'\\
 & +C(\sigma)\left(\sup_{B(x_{\ast},2\rho)}|u|^{2}\right)\int_{B(x_{\ast},2\rho)}\chi_{\rho}^{2}\left(|\nabla\eta_{\epsilon}|^{2}+|\Delta\eta_{\epsilon}|\right)|u|^{2}d\nu_{t}'\\
\leq & -(1-\sigma)\int_{\mathcal{R}_{t}}\eta_{\epsilon}^{2}\chi_{\rho}^{2}|\nabla u|^{2}d\nu_{t}'+\Psi(\epsilon|\rho,\sigma)+\Psi(\rho^{-1}|\sigma).
\end{align*}
Integrating in time gives
\[
\int_{\mathcal{R}_{t_{1}}}\eta_{\epsilon}^{2}\chi_{\rho}^{2}|u|^{2}d\nu_{t_{1}}-\int_{\mathcal{R}_{t_{0}}}\eta_{\epsilon}^{2}\chi_{\rho}^{2}|u|^{2}d\nu_{t_{0}}\leq-(1-\sigma)\int_{t_{0}}^{t_{1}}\int_{\mathcal{R}_{t}}\eta_{\epsilon}^{2}\chi_{\rho}^{2}|\nabla u|^{2}d\nu_{t}dt+\Psi(\rho^{-1}|\sigma)+\Psi(\epsilon|\rho,\sigma).
\]
Taking $\epsilon\searrow0$, the dominated convergence theorem gives
\[
\int_{\mathcal{R}_{t_{1}}}\chi_{\rho}^{2}|u|^{2}d\nu_{t_{1}}+(1-\sigma)\int_{t_{0}}^{t_{1}}\int_{\mathcal{R}_{t}}\chi_{\rho}^{2}|\nabla u|^{2}d\nu_{t}dt\leq\int_{\mathcal{R}_{t_{0}}}\chi_{\rho}^{2}|u|^{2}d\nu_{t_{0}}+\Psi(\rho^{-1}|\sigma).
\]
We then take $\rho\nearrow\infty$, so that the monotone convergence
theorem gives
\[
\int_{\mathcal{R}_{t_{1}}}|u|^{2}d\nu_{t_{1}}(1-\sigma)\int_{t_{0}}^{t_{1}}\int_{\mathcal{R}_{t}}|\nabla u|^{2}d\nu_{t}dt\leq\int_{\mathcal{R}_{t_{0}}}|u|^{2}d\nu_{t_{0}}.
\]
Finally, take $\sigma\searrow0$. 
\end{proof}
The following estimates will be used to show that a heat flow can
be grafted onto a smooth approximating Ricci flow to produce an approximate
solution to the heat equation, with $L^{p}$ estimates. 
\begin{lem}
\label{boringestimates} Suppose $(u_{t})_{t\in(-\infty,0)}\in C^{\infty}(\mathcal{R})\cap C(\mathcal{X})$
is a holomorphic function, and define $u_{\rho,\epsilon,t}:=\chi_{\rho}\eta_{\epsilon}u_{t}$
for $\epsilon\in(0,1]$ and $\rho\geq1$. Then we have the following
for any $-\infty<t_{0}<t_{1}<0$:

$(i)$ $\text{supp}(u_{\rho,\epsilon})\cap\mathcal{X}_{[t_{0},t_{1}]}\subset\subset\mathcal{R}_{[t_{0},t_{1}]}$. 

$(ii)$ For any $p\in[1,\infty)$, if $t_{\ast}\in(-\infty,(p-1)t_{0}]$,
then 

\[
\sup_{t\in[t_{0},t_{1}]}\int_{\mathcal{R}_{t}}|u_{t}|^{p}d\nu_{t}\leq C(p)\int_{\mathcal{R}_{t_{\ast}}}|u_{t_{\ast}}|^{2}d\nu_{t_{\ast}}.
\]

$(iii)$ For any $p\in[1,\infty)$, if $t_{\ast}\in(-\infty,(p-1)t_{0}]$,
then
\[
\sup_{t\in[t_{0},t_{1}]}\int_{\mathcal{R}_{t}}|u_{t}-u_{\rho,\epsilon,t}|^{p}d\nu_{t}\leq\left(\Psi(\epsilon|p,\rho,t_{0},t_{1})+\Psi(\rho^{-1}|p)\right)\int_{\mathcal{R}_{t_{\ast}}}|u_{t_{\ast}}|^{2}d\nu_{t_{\ast}}.
\]

$(iv)$ For any $t_{\ast}\in(-\infty,t_{0})$, we have
\[
\int_{t_{0}}^{t_{1}}\int_{\mathcal{R}_{t}}|(\partial_{t}-\Delta)u_{\rho,\epsilon,t}|^{\frac{3}{2}}d\nu_{t}dt\leq\left(\Psi(\epsilon|\rho,t_{\ast})+\Psi(\rho^{-1}|t_{\ast})\right)\int_{\mathcal{R}_{t_{\ast}}}|u_{t_{\ast}}|^{2}d\nu_{t_{\ast}}.
\]
\end{lem}

\begin{proof}
$(ii)$ Without loss of generality, we have $\int_{\mathcal{R}_{t_{\ast}}}|u|^{2}d\nu_{t_{\ast}}<\infty$,
so that if $(v_{\rho,\epsilon,t})_{t\in[t_{\ast},0)}$ denotes the
heat flow starting at $\chi_{\rho}\eta_{\epsilon}u_{t_{\ast}}$, we
have
\[
u_{t}(x)=\int_{\mathcal{R}_{t_{\ast}}}u_{t_{\ast}}d\nu_{x;t_{\ast}},\,\,\,\,\,\,\,\,\,\,\,\,\,\,v_{\rho,\epsilon,t}(x)=\int_{\mathcal{R}_{t_{\ast}}}\chi_{\rho}\eta_{\epsilon}u_{t_{\ast}}d\nu_{x;t_{\ast}}
\]
for all $x\in\mathcal{X}_{t}$. For sufficiently large $i=i(\epsilon,\rho)\in\mathbb{N}$,
we can then define $v_{i,\rho,\epsilon,t_{\ast}}:=(\psi_{i,t_{\ast}}^{-1})^{\ast}v_{\rho,\epsilon,t_{\ast}}$
, and let $(v_{i,\rho,\epsilon,t})_{t\in[t_{\ast},0]}$ be the corresponding
heat flows. By see Theorem 12.1 of \cite{bamlergen1}, we then have
\[
\left(\int_{M_{i}}|v_{i,\rho,\epsilon,t}|^{p}d\nu_{t}^{i}\right)^{\frac{1}{p}}\leq\left(\int_{M_{i}}|v_{i,\rho,\epsilon,t_{\ast}}|^{2}d\nu_{t_{\ast}}^{i}\right)^{\frac{1}{2}}
\]
if $|t_{\ast}|\geq(p-1)|t|$. For any $t\in[\frac{1}{p-1}t_{\ast},0)$,
Remark \ref{convergenotation}, the uniform boundedness (in $i\in\mathbb{N}$)
of $|v_{i,\rho,\epsilon,t}|$, and the dominated convergence theorem
give the following for any compact subset $L\subseteq\mathcal{R}_{t}$:

\begin{align*}
\int_{L}|v_{\rho,\epsilon,t}|^{p}d\nu_{t}= & \int_{L}\left|\int_{\mathcal{R}_{t_{\ast}}}v_{\rho,\epsilon,t_{\ast}}(y)K(x;y)dg_{t_{\ast}}(y)\right|^{p}d\nu_{t}(x)\\
= & \lim_{i\to\infty}\int_{\psi_{i}(L)}\left|\int_{\mathcal{R}_{t_{\ast}}}v_{\rho,\epsilon,t_{\ast}}(y)K(\psi_{i}^{-1}(x);y)dg_{t_{\ast}}(y)\right|^{p}K^{i}(x_{i},0;x,t)dg_{i,t}(x)\\
= & \lim_{i\to\infty}\int_{\psi_{i}(L)}\left|\int_{M_{i}}v_{i,\rho,\epsilon,t_{\ast}}(y)K^{i}(x;y)dg_{i,t_{\ast}}(y)\right|^{p}K^{i}(x_{i},0;x,t)dg_{i,t}(x)\\
= & \lim_{i\to\infty}\int_{\psi_{i}(L)}|v_{i,\rho,\epsilon,t}|^{p}d\nu_{t}^{i}\\
\leq & \liminf_{i\to\infty}\int_{M_{i}}|v_{i,\rho,\epsilon,t_{\ast}}|^{2}d\nu_{t_{\ast}}^{i}\\
= & \int_{\mathcal{R}_{t_{\ast}}}|v_{\rho,\epsilon,t_{\ast}}|^{2}d\nu_{t_{\ast}}\\
\leq & \int_{\mathcal{R}_{t_{\ast}}}|u|^{2}d\nu_{t_{\ast}}.
\end{align*}
For any $x\in L$, Lemma \ref{goodgrad} gives
\begin{align*}
|v_{\rho,\epsilon,t}|(x)\leq & \int_{\mathcal{R}_{t_{\ast}}}|u_{t_{\ast}}|d\nu_{x;t_{\ast}}\leq C(L,t_{\ast})\int_{\mathcal{R}_{t_{\ast}}}|u_{t_{\ast}}|e^{\frac{1}{8}f_{t_{\ast}}}d\nu_{t_{\ast}}\\
\leq & C(L,t_{\ast})\left(\int_{\mathcal{R}_{t_{\ast}}}|u_{t_{\ast}}|^{2}d\nu_{t_{\ast}}\right)^{\frac{1}{2}}\left(\int_{\mathcal{R}_{t_{\ast}}}e^{\frac{1}{4}f_{t_{\ast}}}d\nu_{t_{\ast}}\right)^{\frac{1}{2}}\\
\leq & C(L,t_{\ast})\left(\int_{\mathcal{R}_{t_{\ast}}}|u_{t_{\ast}}|^{2}d\nu_{t_{\ast}}\right)^{\frac{1}{2}}.
\end{align*}
Moreover, $u_{t_{\ast}}\in L^{1}(\mathcal{R}_{t_{\ast}},\nu_{x;t_{\ast}})$
and the dominated convergence theorem imply
\[
\lim_{\rho\to\infty}\lim_{\epsilon\searrow0}v_{\rho,\epsilon,t}(x)=\lim_{\rho\to\infty}\lim_{\epsilon\searrow0}\int_{\mathcal{R}_{t_{\ast}}}\chi_{\rho}\eta_{\epsilon}ud\nu_{x;t_{\ast}}=\int_{\mathcal{R}_{t_{\ast}}}ud\nu_{x;t_{\ast}}=u_{t}(x)
\]
for all $x\in L$. Another application of the dominated convergence
theorem then yields
\[
\int_{L}|u_{t}|^{p}d\nu_{t}=\lim_{\rho\to\infty}\lim_{\epsilon\searrow0}\int_{L}|v_{\rho,\epsilon,t}|^{p}d\nu_{t}\leq\int_{\mathcal{R}_{t_{\ast}}}|u|^{2}d\nu_{t_{\ast}}.
\]
Because this holds for arbitrary compact subsets $L\subseteq\mathcal{R}_{t}$,
the claim follows. 

\noindent $(iii)$ Because $u_{t}$ is locally bounded, we know
\begin{align*}
\int_{\mathcal{R}_{t}}|u_{t}-u_{\rho,\epsilon,t}|^{p}d\nu_{t}\leq & \left(\text{sup}_{B(x_{\ast},2\rho)}|u_{t}|\right)^{p}\int_{\mathcal{R}_{t}\cap B(x_{\ast},2\rho)}(1-\eta_{\epsilon})^{p}d\nu_{t}\\
 & +\left(\int_{\mathcal{R}_{t}}|u_{t}|^{2p}d\nu_{t}\right)^{\frac{1}{2}}\left(\nu_{t}(\mathcal{R}_{t}\setminus B(x_{\ast},\rho))\right)^{\frac{1}{2}}.
\end{align*}
The claim then follows from $(ii)$ and $\lim_{r\to\infty}\nu_{t}(\mathcal{R}_{t}\setminus B(x_{\ast},\rho))=0$. 

\noindent $(iv)$ We compute
\begin{align*}
(\partial_{t}-\Delta)u_{\rho,\epsilon}= & -(\Delta\chi_{\rho})\eta_{\epsilon}u-\chi_{\rho}(\Delta\eta_{\epsilon})u-2\left(\text{Re}\langle\nabla\chi_{\rho},\overline{\nabla}\eta_{\epsilon}\rangle\right)u\\
 & -\chi_{\rho}\langle\overline{\nabla}\eta_{\epsilon},\nabla u\rangle-\eta_{\epsilon}\langle\overline{\nabla}\chi_{\rho},\nabla u\rangle,
\end{align*}
which allows us to estimate
\begin{align*}
|(\partial_{t}-\Delta)u_{r,\epsilon}|\leq & 2\left(\frac{C}{\rho^{2}}+\frac{C}{\rho}|\nabla\eta_{\epsilon}|+\chi_{\rho}|\Delta\eta_{\epsilon}|\right)|u|\\
 & +\chi_{\rho}|\nabla\eta_{\epsilon}|\cdot|\nabla u|+\frac{C}{\rho}|\nabla u|.
\end{align*}
Integrating in spacetime and appealing to Lemmas \ref{LocLipschitz},
\ref{L2gradient} yields
\begin{align*}
\int_{t_{0}}^{t_{1}}\int_{\mathcal{R}_{t}}|(\partial_{t}-\Delta)u_{\rho,\epsilon}|^{\frac{3}{2}}d\nu_{t}\leq & C\left(\sup_{t\in[t_{0},t_{1}]}\sup_{B(x_{\ast},2\rho)}(|u_{t}|+|\nabla u_{t}|)\right)^{\frac{3}{2}}\int_{\mathcal{R}_{t}\cap\{q\leq4r^{2}\}}(|\Delta\eta_{\epsilon}|^{\frac{3}{2}}+|\nabla\eta_{\epsilon}|^{\frac{3}{2}})d\nu_{t}+\Psi(\rho^{-1})\\
\leq & \Psi(\epsilon|\rho)+\Psi(\rho^{-1}).
\end{align*}
\end{proof}
We now establish a sufficient criterion for a cone to be regular. 
\begin{lem}
\label{conereplacement} Suppose $\mathcal{X}$ is a static metric
flow corresponding to a Ricci-flat cone, obtained as a limit of closed
Kähler-Ricci flows as in (\ref{F}) with $\mathcal{N}_{x_{i},0}^{g_{i}}(\epsilon_{i}^{-1})\ge-Y$,
and let $(\nu_{t})_{t\in(-\infty,0]}$ be the conjugate heat flow
corresponding to the vertex $x_{\ast}\in\mathcal{X}_{0}$, so that
$d\nu_{t}=(2\pi\tau)^{-n}e^{-\frac{d^{2}(x_{\ast},\cdot)}{2\tau}+a}dg_{t}$
for some $a\in\mathbb{R}$. Suppose $u^{1},...,u^{n}\in C^{\infty}(\mathcal{R})\cap C(\mathcal{X})$
are holomorphic functions satisfying the following:

$(i)$ $\int_{\mathcal{R}_{-2}}u^{\alpha}\overline{u}^{\beta}d\nu_{-2}=\delta_{\alpha\beta}$
for $1\leq\alpha,\beta\leq n$,

$(ii)$ $\int_{\mathcal{R}_{-1}}u^{\alpha}\overline{u}^{\beta}d\nu_{-1}=0$
for $1\leq\alpha\neq\beta\leq n$,

$(iii)$ $\int_{\mathcal{R}_{-2}}|u^{\alpha}|^{2}d\nu_{-2}\leq2\int_{\mathcal{R}_{-1}}|u^{\alpha}|^{2}d\nu_{-1}$
for $1\leq\alpha\leq n$,

$(iv)$ $u^{\alpha}(x_{0})=0$ for $1\leq\alpha\leq n$,

$(v)$ $\int_{\mathcal{R}_{-N}}|u^{\alpha}|^{2}d\nu_{-N}<\infty$.

\noindent Then $x_{\ast}\in\mathcal{R}_{0}$, hence $\mathcal{X}$
is regular (and thus isometric to $\mathbb{C}^{n}$). 
\end{lem}

\begin{proof}
Set $B:=\int_{\mathcal{R}_{-N}}|u^{\alpha}|^{2}d\nu_{-N}<\infty$.
Let $(U_{i})$ be a precompact exhaustion of $\mathcal{R}$, with
diffeomorphisms $\psi_{i}:U_{i}\to M_{i}$ realizing smooth convergence
on $\mathcal{R}.$ For any $\rho\in[1,\infty)$ and $\epsilon\in(0,1]$,
we know $u_{\rho,\epsilon}^{\alpha}:=\chi_{\rho}\eta_{\epsilon}u^{\alpha}$
satisfies $\text{supp}(u_{\rho,\epsilon}^{\alpha})\subset\subset\mathcal{R}$,
so we can define 
\[
u_{i,\rho,\epsilon}^{\alpha}:=(\psi_{i}^{-1})^{\ast}u_{\rho,\epsilon}^{\alpha}
\]
for $i=i(\rho,\epsilon)\in\mathbb{N}$ sufficiently large. Fix $\delta>0$
to be determined. Our hypotheses and Lemma \ref{boringestimates}
imply the following when $i=i(\rho,\epsilon,B,\delta)\in\mathbb{N}$
is sufficiently large:

$(i)$ $\left|\int_{M_{i}}u_{i,\rho,\epsilon}^{\alpha}\overline{u}_{i,\rho,\epsilon}^{\beta}d\nu_{x_{i},0;-2}^{i}-\delta_{\alpha\beta}\right|\leq\Psi(i^{-1}|\rho,\epsilon,B)+\Psi(\epsilon|\rho,B)+\Psi(\rho^{-1}|B)$
for $1\leq\alpha,\beta\leq n$,

$(ii)$ $\left|\int_{M_{i}}u_{i,\rho,\epsilon}^{\alpha}\overline{u}_{i,\rho,\epsilon}^{\beta}d\nu_{x_{i},0;-1}^{i}\right|\leq\Psi(i^{-1}|\rho,\epsilon,B)+\Psi(\epsilon|\rho,B)+\Psi(\rho^{-1}|B)$
for $1\leq\alpha\neq\beta\leq n$,

$(iii)$ $\int_{M_{i}}|u_{i,\rho,\epsilon}^{\alpha}|^{2}d\nu_{x_{i},0;-2}^{i}\leq\left(2+\Psi(i^{-1}|\rho,\epsilon,B)+\Psi(\epsilon|\rho,B)+\Psi(\rho^{-1}|B)\right)\int_{M_{i}}|u_{i,r,\epsilon}^{\alpha}|^{2}d\nu_{x_{i},0;-1}^{i}$
for $1\leq\alpha\leq n$,

$(iv)$ $\sup_{t\in[-N,-\delta]}\int_{M_{i}}u_{i,\rho,\epsilon}^{\alpha}d\nu_{x_{i},0;t}\leq\Psi(i^{-1}|\rho,\epsilon,B)+\Psi(\epsilon|\rho,B)+\Psi(\rho^{-1}|B)$
for $1\leq\alpha\leq n$,

$(v)$ $\sup_{t\in[-N,-\delta]}\int_{M_{i}}|u_{i,r,\epsilon}^{\alpha}|^{2}d\nu_{x_{i},0;t}^{i}\leq C(B)$,

$(vi)$ $\int_{-\frac{1}{2}N}^{-\delta}\int_{M_{i}}|(\partial_{t}-\Delta)u_{i,\rho,\epsilon}^{\alpha}|^{\frac{3}{2}}d\nu_{t}^{i}\leq\Psi(i^{-1}|\rho,\epsilon,B,\delta)+\Psi(\epsilon|\rho,B,\delta)+\Psi(\rho^{-1}|B,\delta),$

$(vii)$ $\sup_{t\in[-\frac{1}{2}N,-\delta]}\int_{M_{i}}|\overline{\partial}u_{i,\rho,\epsilon}^{\alpha}|^{2}d\nu_{t}^{i}\leq\Psi(i^{-1}|\rho,\epsilon,B,\delta)+\Psi(\epsilon|\rho,B,\delta)+\Psi(\rho^{-1}|B,\delta)$. 

\noindent Let $\widehat{u}_{i,\rho,\epsilon}^{\alpha}$ be the heat
flow starting at $\widehat{u}_{i,\rho,\epsilon,-50}^{\alpha}$, so
that if we write $\nu_{t}^{i}:=\nu_{x_{i},0;t}^{i}$, then 
\begin{align*}
\frac{d}{dt}\int_{M_{i}}|u_{i,\rho,\epsilon}^{\alpha}-\widehat{u}_{i,\rho,\epsilon}^{\alpha}|^{2}d\nu_{t}^{i}= & 2\int_{M_{i}}(u_{i,\rho,\epsilon}^{\alpha}-\widehat{u}_{i,\rho,\epsilon}^{\alpha})(\partial_{t}-\Delta)u_{i,\rho,\epsilon}^{\alpha}d\nu_{t}^{i}\\
 & -\int_{M}|\nabla^{\mathbb{R}}(u_{i,\rho,\epsilon}^{\alpha}-\widehat{u}_{i,\rho,\epsilon}^{\alpha})|^{2}d\nu_{t}^{i},
\end{align*}
hence we may use $(v),(vi)$, and Theorem 12.1 of \cite{bamlergen1}
to estimate
\begin{align*}
\sup_{t\in[-20,-\delta]}\int_{M_{i}}|u_{i,\rho,\epsilon}^{\alpha}-\widehat{u}_{i,\rho,\epsilon}^{\alpha}|^{2}d\nu_{t}^{i}+ & \int_{-20}^{-\delta}\int_{M_{i}}|\nabla^{\mathbb{R}}(u_{i,\rho,\epsilon}^{\alpha}-\widehat{u}_{i,\rho,\epsilon}^{\alpha})|^{2}d\nu_{t}^{i}\\
\le & C\left(\left(\int_{M_{i}}|u_{i,\rho,\epsilon}^{\alpha}|^{3}d\nu_{t}^{i}\right)^{\frac{1}{3}}+\left(\int_{M_{i}}|\widehat{u}_{i,\rho,\epsilon}^{\alpha}|^{3}d\nu_{t}^{i}\right)^{\frac{1}{3}}\right)\left(\int_{M_{i}}|(\partial_{t}-\Delta)u_{i,\rho,\epsilon}^{\alpha}|^{\frac{3}{2}}d\nu_{t}^{i}\right)^{\frac{2}{3}}\\
\le & \left(1+\left(\int_{M_{i}}|\widehat{u}_{i,\rho,\epsilon}^{\alpha}|^{3}d\nu_{t}^{i}\right)^{\frac{1}{3}}\right)\left(\Psi(i^{-1}|\rho,\epsilon,B,\delta)+\Psi(\epsilon|\rho,B,\delta)+\Psi(\rho^{-1}|B,\delta)\right)\\
\leq & \Psi(i^{-1}|\rho,\epsilon,B,\delta)+\Psi(\epsilon|\rho,B,\delta)+\Psi(\rho^{-1}|B,\delta),
\end{align*}
where the last inequality follows from estimate $(ii)$ of Lemma \ref{boringestimates}
since
\[
\lim_{i\to\infty}\int_{M_{i}}|\widehat{u}_{i,\rho,\epsilon}^{\alpha}|^{3}d\nu_{t}^{i}=\int_{\mathcal{R}_{t}}|u_{\rho,\epsilon}^{\alpha}|^{3}d\nu_{t}\leq\int_{\mathcal{R}_{t}}|u^{\alpha}|^{3}d\nu_{t}.
\]
By choosing $\rho\geq\underline{\rho}(\delta)$, $\epsilon\leq\overline{\epsilon}(\rho,\delta)$,
and $i\geq\underline{i}(\epsilon,\rho,\delta)$, we therefore obtain
$\mathbb{C}$-valued heat flows $v^{1},...,v^{n}$ (which are restrictions
of $\widehat{u}_{i,\rho,\epsilon}^{\alpha}+\lambda_{\alpha}$ for
$|\lambda_{\alpha}|\leq\delta$) which satisfy the following:

$(i)'$ $\left|\int_{M_{i}}v^{\alpha}\overline{v}^{\beta}d\nu_{-2}^{i}-\delta_{\alpha\beta}\right|\leq\delta$
for $1\leq\alpha,\beta\leq n$,

$(ii)$' $\left|\int_{M_{i}}v^{\alpha}\overline{v}^{\beta}d\nu_{-1}^{i}\right|\leq\delta$
for $1\leq\alpha\neq\beta\leq n$,

$(iii)'$ $\int_{M_{i}}|v^{\alpha}|^{2}d\nu_{-2}^{i}\leq2^{1+\delta}\int_{M_{i}}|v^{\alpha}|^{2}d\nu_{-1}^{i}$
for $1\leq\alpha\leq n$,

$(iv)'$ $v^{\alpha}(x_{0})=0$ for $1\leq\alpha\leq n$,

$(v)'$ $\sup_{t\in[-\frac{1}{10}N,-\delta]}\int_{M_{i}}|v^{\alpha}|^{2}d\nu_{t}^{i}\leq C(B)$,

$(vii)'$ $\sup_{t\in[-\frac{1}{10}N,-\delta]}\int_{M_{i}}|\overline{\partial}v^{\alpha}|^{2}d\nu_{t}^{i}\leq\delta$.

\noindent After replacing each $v^{\alpha}$ with $e^{\sqrt{-1}\theta_{\alpha}}v^{\alpha}$
for appropriate $\theta_{\alpha}\in\mathbb{R}$, we may moreover assume
that $y^{2\alpha-1}:=\sqrt{2}\text{Re}(v^{\alpha})$, $y^{2\alpha}:=\sqrt{2}\text{Im}(v^{\alpha})$
satisfy
\[
\left|\int_{M_{i}}(y^{2\alpha})^{2}d\nu_{-2}-1\right|<\delta,
\]
\[
\left|\int_{M_{i}}(y^{2\alpha-1})^{2}d\nu_{-2}-1\right|<\delta,
\]
\[
\left|\int_{M_{i}}y^{2\alpha}y^{2\alpha-1}d\nu_{-2}\right|<\delta,
\]

\noindent so that by $(i)',(ii)',(iv)',(v)'$, we have the following
for $i=i(\delta,B)\in\mathbb{N}$ sufficiently large:

$(i)''$ $\left|\int_{M_{i}}y^{\alpha}y^{\beta}d\nu_{-2}^{i}-\delta_{\alpha\beta}\right|\leq\delta$
for $1\leq\alpha,\beta\leq2n$,

$(ii)''$ $\left|\int_{M_{i}}y^{\alpha}y^{\beta}d\nu_{-1}^{i}\right|\leq\delta$
for $1\leq\alpha\neq\beta\leq2n$,

$(iv)''$ $y^{\alpha}(x_{0})=0$ for $1\leq\alpha\leq n$,

$(v)''$ $\sup_{t\in[-20,-\delta]}\int_{M_{i}}(y^{\alpha})^{2}d\nu_{t}^{i}\leq C(B)$,

\noindent For $1\leq\alpha\leq n$, $\frac{d}{dt}\left(\frac{1}{\tau}\int_{M_{i}}|y|^{2}d\nu_{t}\right)\leq0$
implies
\[
\int_{M_{i}}(y^{2\alpha-1})^{2}d\nu_{-2}^{i}\geq2\int_{M}(y^{2\alpha-1})^{2}d\nu_{-1}^{i},
\]
\[
\int_{M_{i}}(y^{2\alpha})^{2}d\nu_{-2}^{i}\geq2\int_{M}(y^{2\alpha-1})^{2}d\nu_{-1}^{i},
\]
so that $|v^{\alpha}|^{2}=\frac{1}{2}(y^{2\alpha})^{2}+\frac{1}{2}(y^{2\alpha-1})^{2}$
and $(iii)'$ imply
\[
(iii)''\,\,\,\,\,\,\,\,\,\,\,\,\,\,\,\,\,\,\,\,\int_{M_{i}}(y^{\beta})^{2}d\nu_{-2}^{i}\leq2^{1+\delta}\int_{M_{i}}(y^{\beta})^{2}d\nu_{-1}^{i},
\]
for $1\leq\beta\leq2n$. Moreover, $(vii)'$ implies
\[
(vii)''\,\,\,\,\,\,\,\,\,\,\,\,\,\,\,\,\,\,\,\sup_{t\in[-20,-\delta]}\int_{M_{i}}|\nabla^{\mathbb{R}}y^{2\alpha-1}-J\nabla^{\mathbb{R}}y^{2\alpha}|^{2}d\nu_{t}^{i}\leq C\delta
\]
for some universal $C<\infty$. In particular, we have $\int_{M_{i}}|\langle\nabla^{\mathbb{R}}y^{2\alpha-1},\nabla^{\mathbb{R}}y^{2\alpha}\rangle|^{2}d\nu_{t}^{i}\leq C\delta$
for $1\leq\alpha\leq n$ and $t\in[-20,-\delta]$. Set $a_{\alpha}:=\frac{1}{2}\int_{M_{i}}(y^{\alpha})^{2}d\nu_{-1}^{i}\geq\frac{1}{4}$
and $\widehat{y}^{\alpha}:=a^{-\frac{1}{2}}y^{\alpha}$. By $(iii)'',(iv)'',(v)''$,
and Lemma \ref{frequency}, we know that for any $\theta\leq\overline{\theta}$,
we have
\begin{equation}
\int_{-1}^{-\theta}\int_{M_{i}}\left||\nabla^{\mathbb{R}}y^{\alpha}|^{2}-1\right|d\nu_{t}dt+\int_{-1}^{-\theta}\int_{M_{i}}|\nabla^{\mathbb{R}}\nabla^{\mathbb{R}}y^{\alpha}|^{2}d\nu_{t}dt\leq\Psi(\delta|\theta)\label{offdiagonal}
\end{equation}

\noindent for $1\leq\alpha\leq2n$. Integrating
\[
\frac{d}{dt}\int_{M_{i}}(y^{\alpha}-y^{\beta})^{2}d\nu_{t}^{i}=-\int_{M_{i}}|\nabla^{\mathbb{R}}(y^{\alpha}-y^{\beta})|^{2}d\nu_{t}^{i}
\]
from $t=-2$ to $t=-1$ yields (for $1\leq\alpha\neq\beta\leq2n$)
\begin{align*}
\int_{-2}^{-1}\int_{M_{i}}|\nabla^{\mathbb{R}}(y^{\alpha}-y^{\beta})|^{2}d\nu_{t}^{i}dt= & -\left.\left(\int_{M_{i}}(y^{\alpha})^{2}d\nu_{t}+\int_{M_{i}}(y^{\beta})^{2}d\nu_{t}^{i}\right)\right|_{t=-2}^{t=-1}+\Psi(\delta)\\
= & 1+\Psi(\delta),
\end{align*}
where we used $(i)'',(ii)'',(iii)''$. Because $\frac{d}{dt}\int_{M_{i}}|\nabla^{\mathbb{R}}(y^{\alpha}-y^{\beta})|^{2}d\nu_{t}^{i}\leq0$,
this implies
\[
\sup_{t\in[-1,\theta]}\int_{M_{i}}|\nabla^{\mathbb{R}}(y^{\alpha}-y^{\beta})|^{2}d\nu_{t}\leq1+\Psi(\delta).
\]
On the other hand, the $L^{2}$-Poincare inequality gives
\[
\int_{M_{i}}|\nabla(y^{\alpha}-y^{\beta})|^{2}d\nu_{-1}\geq\int_{M_{i}}(y^{\alpha}-y^{\beta})^{2}d\nu_{-1}\geq1-\Psi(\delta),
\]

\noindent so combining with (\ref{offdiagonal}) gives
\begin{align*}
\left|\int_{-1}^{-\theta}\int_{M_{i}}\langle\nabla y^{\alpha},\nabla y^{\beta}\rangle d\nu_{t}^{i}dt\right|= & \frac{1}{2}\left|\int_{-1}^{-\theta}\int_{M_{i}}\left(|\nabla(y^{\alpha}-y^{\beta})|^{2}-|\nabla y^{\alpha}|^{2}-|\nabla y^{\beta}|^{2}\right)d\nu_{t}^{i}dt\right|\\
\leq & \frac{1}{2}\left|\int_{-1}^{-\theta}\int_{M_{i}}\left(|\nabla(y^{\alpha}-y^{\beta})|^{2}-1\right)d\nu_{t}^{i}dt\right|\\
\leq & \Psi(\delta|\theta).
\end{align*}
Because 
\[
\left|\frac{d}{dt}\int_{M_{i}}\langle\nabla y^{\alpha},\nabla y^{\beta}\rangle d\nu_{t}^{i}\right|\leq2\int_{M_{i}}|\nabla^{2}y^{\alpha}|\cdot|\nabla^{2}y^{\beta}|d\nu_{t}^{i},
\]
we can upgrade this to 
\[
\int_{-1}^{-\theta}\left|\int_{M_{i}}\langle\nabla\widehat{y}^{\alpha},\nabla\widehat{y}^{\beta}\rangle d\nu_{t}^{i}\right|dt\leq\Psi(\delta|\theta),
\]
so the $L^{1}$-Poincare inequality (Theorem 11.1 of \cite{bamlergen1})
gives
\[
\int_{-1}^{-\theta}\int_{M_{i}}|\langle\nabla\widehat{y}^{\alpha},\nabla\widehat{y}^{\beta}\rangle|d\nu_{t}^{i}\leq\Psi(\delta|\theta).
\]
In other words, $\widehat{y}:=(\widehat{y}^{1},...,\widehat{y}^{2n})$
is a strong $(2n,\epsilon,\epsilon^{-1}\theta)$-splitting map when
$\theta\leq\overline{\theta}(\epsilon)$ and $\delta\leq\overline{\delta}(\epsilon)$,
and $(x_{0},0)$ is $(\epsilon,\theta)$-static. By choosing $\epsilon=\epsilon(Y,n)>0$
as in Proposition 14.1 of \cite{bamlergen3}, we obtain $r_{Rm}^{g_{i}}(x_{i},0)\geq\epsilon\theta$
for sufficiently large $i\in\mathbb{N}$, which implies the claim. 
\end{proof}
We will use Lemma \ref{conereplacement} along with a blowup argument
to obtain a criterion for a point in $\mathcal{X}$ to be regular.
The following is a technical preliminary step.
\begin{lem}
\label{needforconvergence} Suppose $t_{0}\in(-\infty,0)$, $x_{0}\in\mathcal{X}_{t_{0}}$,
$r_{0}>0$, and $D<\infty$. If $\tau\leq\overline{\tau}(x_{0},r_{0},\epsilon)$,
then for any holomorphic function $u\in L^{2}(B(x_{0},2r_{0})\cap\mathcal{R}_{t_{0}-N\tau},d\nu_{t_{0}-N\tau})$
and any $a\in[1,2]$, the following hold:

$(i)$ $\int_{\left(B(x_{0},2r_{0})\setminus B(x_{0},D\sqrt{\tau})\right)\cap\mathcal{R}_{t_{0}-a\tau}}|u|^{2}d\nu_{t_{0}-a\tau}\leq\Psi(D^{-1}|Y,r_{0})\int_{\mathcal{R}_{t_{0}-N\tau}\cap B(x_{0},2r_{0})}|u|^{2}d\nu_{t_{0}-N\tau},$

$(ii)$ $\sup_{B(x_{0},D\sqrt{\tau})}|u|^{2}\leq C(Y,D,r_{0})\int_{B(x_{0},2r_{0})\cap\mathcal{R}_{t_{0}-N\tau}}|u|^{2}d\nu_{t_{0}-N\tau},$

$(iii)\left|\int_{B(x_{0},2r_{0})\cap\mathcal{R}_{t_{0}-\tau}}ud\nu_{t_{0}-\tau}\right|\leq|u(x_{0})|+\Psi(\tau|Y,r_{0})\int_{B(x_{0},2r_{0})\cap\mathcal{R}_{t_{0}-N\tau}}|u|^{2}d\nu_{t_{0}-N\tau}.$
\end{lem}

\begin{proof}
Choose $\varphi\in C_{c}(B(x_{0},2r_{0}))$ such that $\varphi|_{B(x_{0},r_{0})}\equiv1$
and $|\nabla\varphi|\leq2r_{0}^{-1}$. Then for all $x\in B(x_{0},r_{0})$,
we can compute
\begin{align*}
\frac{d}{dt}\int_{\mathcal{R}_{t}}\varphi(y)\eta_{\epsilon}(y)K(x;y)u(y)dg_{t}(y)= & -\int_{\mathcal{R}_{t}}\varphi(y)\eta_{\epsilon}(y)\Delta_{y}K(x;y)u(y)dg_{t}(y)\\
= & 2\int_{\mathcal{R}_{t}}u(y)\text{Re}\langle\nabla_{y}K(x;y),\varphi(y)\overline{\nabla}\eta_{\epsilon}(y)+\eta_{\epsilon}(y)\overline{\nabla}\varphi(y)\rangle dg_{t}(y).
\end{align*}
Letting $A(x_{0},r_{1},r_{2})$ denote the annulus $B(x_{0},r_{2})\setminus\overline{B}(x_{0},r_{1})$
and arguing as in Lemma \ref{LocLipschitz} yields
\begin{align*}
\left|u^{\alpha}(x)-\int_{\mathcal{R}_{t}}\varphi(y)K(x;y)u^{\alpha}(y)dg_{t}(y)\right|\leq & C(Y)\int_{t}^{0}\int_{A(x_{0},r_{0},2r_{0})\cap\mathcal{R}_{s}}|\nabla\log_{y}K(x;y)|\cdot|u^{\alpha}(y)|K(x;y)dg_{s}(y)ds.
\end{align*}
In particular, taking $x=x_{0}$ gives
\begin{align*}
\left|\int_{\mathcal{R}_{t}}\varphi u^{\alpha}d\nu_{x_{0};t}\right|\leq & |u^{\alpha}(x_{0})|+C(Y)\int_{t}^{0}\left(\int_{A(x_{0},r_{0},2r_{0})\cap\mathcal{R}_{s}}|\nabla\log_{y}K(x_{0};y)|^{2}d\nu_{x_{0};s}(y)ds\right)^{\frac{1}{2}}\\
 & \times\left(\int_{A(x_{0},r_{0},2r_{0})\cap\mathcal{R}_{s}}|u^{\alpha}|^{2}d\nu_{x_{0};s}\right)^{\frac{1}{2}}ds.
\end{align*}
For any $s\in[t,0)$, we can estimate
\begin{align*}
\int_{A(x_{0},r_{0},2r_{0})\cap\mathcal{R}_{s}}|u^{\alpha}|^{2}d\nu_{x_{0};s}\leq & \frac{C(Y)}{|s|^{n}}\int_{A(x_{0},r_{0},2r_{0})\cap\mathcal{R}_{s}}|u^{\alpha}|^{2}\exp\left(-\frac{d^{2}(x_{0}(s),y)}{10|s|}\right)dg_{s}(y)\\
\le & \frac{C(Y)}{r_{0}^{2n}}\frac{r_{0}^{2n}}{|s|^{n}}\exp\left(-\frac{r_{0}^{2}}{20|s|}\right)\int_{A(x_{0},r_{0},2r_{0})\cap\mathcal{R}_{Nt}}|u^{\alpha}|^{2}\exp\left(-\frac{d^{2}(x_{0}(Nt),y)}{N|t|}\right)dg_{Nt}(y)\\
\leq & C(Y,r_{0})\frac{|t|^{n}}{r_{0}^{2n}}\exp\left(-\frac{r_{0}^{2}}{20|t|}\right)\int_{B(x_{0},2r_{0})\cap\mathcal{R}_{Nt}}|u^{\alpha}|^{2}d\nu_{x_{0};Nt},
\end{align*}
so combining estimates and using Lemma \ref{goodgrad} yields
\begin{align*}
\left|\int_{\mathcal{R}_{t}}\varphi u^{\alpha}d\nu_{x_{0};t}\right|\leq & |u(x_{0})|+C(Y,r_{0})\frac{|t|^{n}}{r_{0}^{2n}}\exp\left(-\frac{r_{0}^{2}}{20|t|}\right)\left(\int_{B(x_{0},2r_{0})\cap\mathcal{R}_{Nt}}|u^{\alpha}|^{2}d\nu_{x_{0};Nt}\right)^{\frac{1}{2}}\int_{t}^{0}\frac{1}{\sqrt{|s|}}ds\\
\leq & |u(x_{0})|+C(Y,r_{0})\sqrt{|t|}\left(\int_{B(x_{0},2r_{0})\cap\mathcal{R}_{Nt}}|u^{\alpha}|^{2}d\nu_{x_{0};Nt}\right)^{\frac{1}{2}},
\end{align*}
hence $(iii)$ will follow from $(i)$. For any $D<\infty$ we have
\begin{align*}
\int_{A(x_{0},D\sqrt{|t|},r_{0})\cap\mathcal{R}_{t}}|u^{\alpha}|^{2}d\nu_{x_{0};t}\leq & \int_{A(x_{0},D\sqrt{|t|},r_{0})\cap\mathcal{R}_{t}}|u^{\alpha}|^{2}(y)\frac{C(Y)}{|t|^{n}}\exp\left(-\frac{d^{2}(x_{0}(t),y)}{10|t|}\right)dg_{t}(y)\\
\leq & C(Y)\exp\left(-\frac{D^{2}}{20}\right)\int_{A(x_{0},D\sqrt{|t|},r_{0})\cap\mathcal{R}_{t}}|u^{\alpha}|^{2}(y)\frac{1}{|Nt|^{n}}\exp\left(-\frac{d^{2}(x_{0}(Nt),y)}{|Nt|}\right)dg_{Nt}(y)\\
\leq & C(Y)\exp\left(-\frac{D^{2}}{20}\right)\int_{B(x_{0},r_{0})\cap\mathcal{R}_{Nt}}|u^{\alpha}|^{2}d\nu_{x_{0};Nt},
\end{align*}
so $(i)$ follows. For any $x\in B(x_{0},D\sqrt{|t|})$, we can also
estimate
\begin{align*}
|u^{\alpha}(x)|\leq & C\int_{t}^{0}\int_{A(x_{0},r_{0},2r_{0})\cap\mathcal{R}_{s}}|\nabla_{y}\log K(x;y)|\cdot|u^{\alpha}(y)|K(x;y)dg_{s}(y)ds\\
 & +\int_{B(x_{0},r_{0})\cap\mathcal{R}_{t}}|u^{\alpha}|d\nu_{x;t}\\
\leq & C\int_{t}^{0}\left(\int_{B(x_{0},2r_{0})\cap\mathcal{R}_{s}}|\nabla_{y}\log K(x,y)|^{2}d\nu_{x;s}(y)\right)^{\frac{1}{3}}\\
 & \times\left(\int_{A(x_{0},r_{0},2r_{0})\cap\mathcal{R}_{s}}|u^{\alpha}|^{2}d\nu_{x;s}\right)^{\frac{1}{2}}ds\\
 & +\int_{B(x_{0},r_{0})\cap\mathcal{R}_{t}}|u^{\alpha}|d\nu_{x;t}.
\end{align*}
For $x\in B(x_{0},D\sqrt{|t|})$ and $y\in A(x_{0},r_{0},2r_{0})$,
we have
\begin{align*}
K(x;y)\leq & \frac{C(Y)}{|s|^{n}}\exp\left(-\frac{d^{2}(x(s),y)}{10|s|}\right)\\
\leq & \frac{C(Y)}{r_{0}^{2n}}\frac{r_{0}^{2n}}{|s|^{n}}\exp\left(-\frac{r_{0}^{2}}{40|s|}\right)K(x_{0};y(Nt)),
\end{align*}
so that
\[
\int_{A(x_{0},r_{0},2r_{0})\cap\mathcal{R}_{s}}|u^{\alpha}|^{2}d\nu_{x;s}\leq C(Y,r_{0})\left(\sup_{\sigma\in(r_{0}^{2}|t|^{-1},\infty)}\sigma^{2n}e^{-\frac{1}{40}\sigma^{2}}\right)\int_{A(x_{0},r_{0},2r_{0})}|u^{\alpha}|^{2}d\nu_{x_{0};Nt}.
\]
Combining this with estimate $(ii)$ of Lemma \ref{goodgrad} yields
\[
|u^{\alpha}(x)|\leq\Psi(\tau|Y,r_{0})\int_{B(x_{0},2r_{0})}|u^{\alpha}|^{2}d\nu_{x_{0};Nt}+\int_{B(x_{0},r_{0})\cap\mathcal{R}_{t}}|u^{\alpha}|d\nu_{x;t}.
\]
For any $x\in B(x_{0},D\sqrt{|t|})$, estimate $(i)$ of Lemma \ref{goodgrad}
gives (writing $d\nu_{x_{0};t}=(4\pi\tau)^{-n}e^{-f_{x_{0},t}}dg_{t}$)
\begin{align*}
\int_{B(x_{0},r_{0})\cap\mathcal{R}_{t}}|u^{\alpha}|d\nu_{x;t}\leq & C(Y,D)\int_{B(x_{0},r_{0})\cap\mathcal{R}_{t}}|u^{\alpha}|e^{\frac{1}{4}f_{x_{0};t}}d\nu_{x_{0};t}\\
\leq & C(Y,D)\left(\int_{B(x_{0},r_{0})\cap\mathcal{R}_{t}}|u^{\alpha}|^{2}d\nu_{x_{0};t}\right)^{\frac{1}{2}}\left(\int_{B(x_{0},r_{0})\cap\mathcal{R}_{t}}e^{\frac{1}{2}f_{x_{0};t}}d\nu_{x_{0};t}\right)^{\frac{1}{2}}\\
\leq & C(Y,D)\left(\int_{B(x_{0},r_{0})\cap\mathcal{R}_{t}}|u^{\alpha}|^{2}d\nu_{x_{0};t}\right)^{\frac{1}{2}}.
\end{align*}
The claim then follows from Proposition \ref{heatkernel}.
\end{proof}
Next, we show that the normalized $L^{2}$-norm of a holomorphic function
on an open set is almost-monotone by a blowup argument.
\begin{lem}
\label{almostmonotone} For any $\epsilon>0$, $t_{0}\in(-\infty,0)$,
$x_{0}\in\mathcal{X}_{t_{0}}$, and $r_{0}>0$, the following holds
whenever $\tau\leq\overline{\tau}(x_{0},r_{0},\epsilon)$. For any
holomorphic $u\in L^{2}(B(x_{0},2r_{0})\cap\mathcal{R}_{t_{0}-2\tau},d\nu_{t_{0}-2\tau})$
satisfying $u(x_{0})=0$, we have
\[
\int_{\mathcal{R}_{t_{0}-2\tau}\cap B(x_{0},2r_{0})}|u|^{2}d\nu_{t_{0}-2\tau}\geq2^{1-\epsilon}\int_{\mathcal{R}_{t_{0}-\tau}\cap B(x_{0},2r_{0})}|u|^{2}d\nu_{t_{0}-\tau}.
\]
\end{lem}

\begin{proof}
By parabolic rescaling and time translation, we may assume $r_{0}=\frac{1}{2}$
and $t_{0}=0$. Suppose by way of contradiction that there are $\tau_{j}\searrow0$
and holomorphic $u_{j}\in L^{2}(B(x_{0},1)\cap\mathcal{R}_{-2\tau_{j}},d\nu_{-2\tau_{j}})$
satisfying
\[
\int_{B(x_{0},1)\cap\mathcal{R}_{-\tau_{j}}}u_{j}d\nu_{-\tau_{j}}=0,
\]
\[
\int_{B(x_{0},1)\cap\mathcal{R}_{-2\tau_{j}}}|u_{j}|^{2}d\nu_{-2\tau_{j}}\leq2^{1-\epsilon}\int_{B(x_{0},1)\cap\mathcal{R}_{-\tau_{j}}}|u_{j}|^{2}d\nu_{-\tau_{j}}.
\]
By normalizing $u_{j}$, we may assume that $\int_{\mathcal{R}_{-\tau_{j}}\cap B(x_{0},1)}|u_{j}|^{2}d\nu_{-\tau_{j}}=1$
for all $j\in\mathbb{N}$. Consider the parabolic rescalings $\mathcal{X}^{j}$
of $\mathcal{X}$ centered at $x_{0}$, so that $(\mathcal{X}_{t}^{j},d_{t}^{j})=(\mathcal{X}_{\tau_{j}t},\tau_{j}^{-\frac{1}{2}}d_{\tau_{j}t})$
for $t\in(-\infty,0]$. By passing to a subsequence, we may assume
that the metric flow pairs $(\mathcal{X}^{j},(\nu_{x_{0};t}^{j})_{t\in(-\infty,0]})$
(where $\nu_{t}^{j}=\nu_{x_{0};\tau_{j}t}$) $\mathbb{F}$-converge
to some static metric flow $(\mathcal{Y},(\mu_{t}^{\mathcal{Y}})_{t\in(-\infty,0]})$
corresponding to a Ricci-flat metric cone, with smooth convergence
on the regular part. Let $\mathcal{R}^{\mathcal{Y}}$ be the regular
set of $\mathcal{Y}$, and let $(V_{i})$ be a precompact exhaustion
of $\mathcal{R}^{\mathcal{Y}}$ so that there are time-preserving
diffeomorphisms $\zeta_{i}:V_{i}\to\mathcal{R}$ realizing the smooth
convergence (for more details, see Theorem 9.31 of \cite{bamlergen2}).
By the lower heat kernel bounds \ref{heatkernel} and local elliptic
regularity, we can pass to a further subsequence to find $u_{\infty}\in L^{2}(\mathcal{R}_{-2},d\nu_{-2})$
such that $\zeta_{j}^{\ast}u_{j}\to u_{\infty}$ in $C_{\text{loc}}^{\infty}(\mathcal{R})$
and 
\[
\int_{\mathcal{R}_{-2}}|u_{\infty}|^{2}d\nu_{-2}\leq2^{1-\epsilon}.
\]
We also know from Lemma \ref{needforconvergence} that
\[
\int_{\mathcal{R}_{-1}}|u_{\infty}|^{2}d\nu_{-1}=\lim_{j\to\infty}\int_{\mathcal{R}_{-\tau_{j}}\cap B(x_{0},1)}|u_{j}|^{2}d\nu_{-\tau_{j}}=1,
\]
\[
\int_{\mathcal{R}_{-1}}u_{\infty}d\nu_{-1}=0,
\]
hence by Lemma \ref{LocLipschitz}, we have $\int_{\mathcal{R}_{-1}}u_{\infty}d\nu_{t}=0$
for all $t\in(-\infty,0)$. 

On the other hand, the proof of Lemma \ref{L2gradient} gives the
following for $\epsilon,\sigma>0$ and $\rho<\infty$:
\begin{align*}
\frac{d}{dt}\int_{\mathcal{R}_{t}}\eta_{\epsilon}^{2}\chi_{\rho}^{2}|u|^{2}d\nu_{t}\leq & -(1-4\sigma)\int_{\mathcal{R}_{t}}\eta_{\epsilon}^{2}\chi_{\rho}^{2}|\nabla u|^{2}d\nu_{t}+\Psi(\rho^{-1}|\sigma)+\Psi(\epsilon|\rho,\sigma)\\
\leq & -(1-4\sigma)\int_{\mathcal{R}_{t}}|\nabla(\eta_{\epsilon}\chi_{\rho}u)|^{2}d\nu_{t}+\Psi(\rho^{-1}|\sigma)+\Psi(\epsilon|\rho,\sigma).
\end{align*}
The $L^{2}$-Poincare inequality yields
\[
\int_{\mathcal{R}_{t}}|\nabla(\eta_{\epsilon}\chi_{\rho}u)|^{2}d\nu_{t}\leq\frac{1}{\tau}\int_{\mathcal{R}_{t}}|\eta_{\epsilon}\chi_{\rho}u|^{2}d\nu_{t}+\Psi(\rho^{-1}|\sigma)+\Psi(\epsilon|\rho,\sigma),
\]
so that after combining estimates, we have
\[
\frac{d}{dt}\left(\frac{1}{\tau^{1-4\sigma}}\int_{\mathcal{R}_{t}}\eta_{\epsilon}^{2}\chi_{\rho}^{2}|u|^{2}d\nu_{t}\right)\leq\Psi(\rho^{-1}|\sigma)+\Psi(\epsilon|\rho,\sigma).
\]
Integrating in time, and taking $\epsilon\searrow0$, $\rho\nearrow\infty$,
then $\sigma\searrow0$ gives
\[
\int_{\mathcal{R}_{-1}}|u|^{2}d\nu_{-1}\leq\frac{1}{2}\int_{\mathcal{R}_{-2}}|u|^{2}d\nu_{-2}.
\]
\end{proof}
We combine these technical results with Proposition \ref{conereplacement}
to obtain the following.
\begin{prop}
For any $x_{0}\in\mathcal{X}_{t_{0}}$ corresponding to a regular
point of the analytic variety, we have $x_{0}\in\mathcal{R}_{t_{0}}$. 
\end{prop}

\begin{proof}
We can assume $t_{0}=0$. By hypothesis, we can find a holomorphic
chart $(U,(u_{0}^{\alpha})_{\alpha=1}^{n})$ for $\mathcal{X}_{0}$
centered at $x_{0}$ (so that $u^{\alpha}(x_{0})=0$ for $1\leq\alpha\leq n$).
Choose $r_{0}>0$ such that $B(x_{0},2r_{0})\subseteq U$. 

\noindent \textbf{Claim 1:} After possibly shrinking $r_{0}>0$, there
exists $b>0$ such that for all $r\in(0,r_{0}]$ and $x\in\partial B(x_{0},r)$,
there exists holomorphic $h\in L^{2}(\mathcal{X}_{0},d\nu_{0})$ such
that $h(x_{0})=0$, $\sup_{B(x_{0},2r)}|h|=1$, and $|h(x)|>b$.

The proof proceeds as in Claim 4.1 of \cite{szek1}. Suppose by way
of contradiction that there exist $r_{i}\searrow0$, $b_{i}\searrow0$
and $x_{i}\in\partial B(x_{0},r_{i})$ such that for all $h\in\mathcal{O}(\mathcal{R}_{0})\cap L^{2}(\mathcal{R}_{0},d\nu_{0})$
with $h(x_{0})=0$ and $\sup_{B(x_{0},2r_{i})}|h|=1$, we have $|h(x_{i})|<b_{i}$.
By passing to a subsequence, we may moreover assume $(\mathcal{X}_{0},r_{i}^{-1}d_{0},x_{0})$
converges in the pointed Gromov-Hausdorff sense to a metric cone $(C(Y),d_{C(Y)},y_{\ast})$,
and $x_{i}\to x_{\infty}\in\partial B(y_{\ast},1)$ with respect to
the Gromov-Hausdorff convergence. Moreover, the convergence is polarized
in the sense that the sequence
\[
(\mathcal{X}_{0},r_{i}^{-\frac{1}{2}}d_{0},x_{0},\mathcal{R}_{0},r_{i}^{-1}g_{0},J_{0},\mathcal{O}_{\mathcal{X}_{0}},(2\pi)^{-n}e^{-r_{i}^{-1}f_{0}})
\]
of polarized Kähler singular spaces converge to $C(Y)$ equipped with
a holomorphic Hermitian line bundle $(L_{\infty},h_{\infty})$. By
arguing as in the previous section, we can find an almost-holomorphic
section of $L_{\infty}$ vanishing at $y_{\ast}$ but nonzero at $x_{\infty}$,
use a cutoff function to transplant this section to $(\mathcal{X}_{0},r_{i}^{-2}g_{t_{0}})$,
and then use Theorem \ref{hormandersection} to solve the $\overline{\partial}$-equation,
thus obtaining holomorphic functions on $\mathcal{X}_{0}$ with the
desired property if $i\in\mathbb{N}$ is sufficiently large, a contradiction.
$\square$

\noindent \textbf{Claim 2: }There exist $\beta=\beta>0$ and $C<\infty$
such that
\[
C^{-1}d^{\beta}(x,x_{0})\leq|u(x)|\leq Cd(x,x_{0})
\]
for all $x\in\mathcal{X}_{0}$ sufficiently close to $x_{0}$, where
$u=(u^{1},...,u^{n})$. 

The proof is identical to Claim 4.2 of \cite{szek2} given Claim 1,
so we omit it. $\square$

Suppose by way of contradiction that $x_{0}\notin\mathcal{R}_{0}$.
Let $\mathcal{V}$ be the $n$-dimensional $\mathbb{C}$-linear vector
space spanned by $u^{1},...,u^{n}$, and define Hermitian metrics
$\Omega_{t}$ on $\mathcal{V}$ by 
\[
\Omega_{\tau}(u,v):=\frac{1}{\tau}\int_{B(x_{0},2r_{0})\cap\mathcal{R}_{-\tau}}u\overline{v}d\nu_{x_{0};-\tau}
\]
for $u,v\in\mathcal{V}$. For $\delta>0$ to be determined, Lemma
\ref{almostmonotone} gives $\Omega_{\tau}\leq2^{\Psi(\tau)}\Omega_{2\tau}$
in the sense of Hermitian forms on $\mathcal{V}$. 

\noindent \textbf{Claim 3: }If $\tau_{0}>0$ is sufficiently small,
then there exists $\epsilon>0$ such that for any $\tau\in(0,\tau_{0}]$,
we have
\[
\text{tr}_{\Omega_{N\tau}}(\Omega_{\tau})<N^{-100n\epsilon}.
\]

Suppose by way of contradiction that there are $\tau_{j}\searrow0$
and $\epsilon_{j}\searrow0$ such that $\text{tr}_{\Omega_{N\tau_{j}}}(\Omega_{\tau_{j}})\geq N^{-\epsilon_{j}}$
for all $j\in\mathbb{N}$. Let $(u_{j}^{\alpha})_{\alpha=1}^{n}$
be unitary bases for $(\mathcal{V},\Omega_{2\tau_{j}})$ such that
$\Omega_{\tau_{j}}(u_{j}^{\alpha},u_{j}^{\beta})=0$ for $1\leq\alpha\neq\beta\le n$.
This means $v_{j}^{\alpha}:=\sqrt{2\tau}u_{j}^{\alpha}$ satisfy the
following:

$(i)$ $\int_{B(x_{0},2r_{0})\cap\mathcal{R}_{-2\tau_{j}}}v_{j}^{\alpha}\overline{v_{j}^{\beta}}d\nu_{x_{0};-2\tau_{j}}=\delta_{\alpha\beta}$
for $1\leq\alpha,\beta\leq n$,

$(ii)$ $\int_{B(x_{0},2r_{0})\cap\mathcal{R}_{-\tau_{j}}}v_{j}^{\alpha}\overline{v_{j}^{\beta}}d\nu_{x_{0};-\tau_{j}}=0$
for $1\leq\alpha\neq\beta\leq n$.

\noindent Because $\Omega_{2^{k-1}\tau_{j}}\leq2^{\Psi(j^{-1})}\Omega_{2^{k}\tau_{j}}$
for $k=1,...,10$, the assumption $\text{tr}_{\Omega_{N\tau_{j}}}(\Omega_{\tau_{j}})\geq2^{-\epsilon_{j}}$
implies that (after replacing $\epsilon_{j}$ with a slightly larger
sequence $\epsilon_{j}'\searrow0$) $\text{tr}_{\Omega_{2\tau_{j}}}(\Omega_{\tau_{j}})\geq2^{-\epsilon_{j}}$
as well. That is, we have
\[
\sum_{\alpha=1}^{n}\int_{B(x_{0},2r_{0})\cap\mathcal{R}_{-\tau_{j}}}|v_{j}^{\alpha}|^{2}d\nu_{x_{0};-\tau_{j}}\geq n2^{-1-\epsilon_{j}}.
\]
Because $\int_{B(x_{0},2r_{0})\cap\mathcal{R}_{-\tau_{j}}}|v_{j}^{\alpha}|^{2}d\nu_{x_{0};-\tau_{j}}\leq2^{1+\Psi(j^{-1})}$
for $1\leq\alpha\leq n$, we can thus assume

$(iii)$ $\int_{B(x_{0},2r_{0})\cap\mathcal{R}_{-\tau_{j}}}|v_{j}^{\alpha}|^{2}d\nu_{x_{0};-\tau_{j}}\geq2^{-1-\epsilon_{j}}\int_{B(x_{0},2r_{0})\cap\mathcal{R}_{-2\tau_{j}}}|v_{j}^{\alpha}|^{2}d\nu_{x_{0};-2\tau_{j}}$
for $1\leq\alpha\leq n$. 

\noindent Moreover, $\text{tr}_{\Omega_{N\tau_{j}}}(\Omega_{\tau_{j}})\geq N^{-\epsilon_{j}}$
and $\Omega_{\tau_{j}}\leq2^{\Psi(j^{-1})}\Omega_{N\tau_{j}}$ imply
$\Omega_{N\tau_{j}}\leq2\Omega_{\tau_{j}}$ for sufficiently large
$j\in\mathbb{N}$, hence we also have

$(v)$ $\int_{B(x_{0},2r_{0})\cap\mathcal{R}_{-N\tau_{j}}}|v_{j}^{\alpha}|^{2}d\nu_{x_{0};-N\tau_{j}}\leq2^{11}\int_{B(x_{0},2r_{0})\cap\mathcal{R}_{-\tau_{j}}}|v_{j}^{\alpha}|^{2}d\nu_{x_{0};-\tau_{j}}$.

\noindent As in the proof of Lemma \ref{almostmonotone}, pass to
a subsequence so that the parabolic rescalings $\mathcal{X}^{i}$
of $\mathcal{X}$ centered at $x_{0}$ $\mathbb{F}$-converge to some
static metric flow $(\mathcal{Y},(\mu_{t}^{\mathcal{Y}})_{t\in(-\infty,0]})$
corresponding to a Ricci-flat metric cone, with smooth convergence
on the regular part $\mathcal{R}^{\mathcal{Y}}$ given by time-preserving
diffeomorphisms $\zeta_{i}:V_{i}\to\mathcal{R}$. Pass to a further
subsequence to assume that $\zeta_{i}^{\ast}v_{i}^{\alpha}$ converge
in $C_{\text{loc}}^{\infty}(\mathcal{R}^{\mathcal{Y}})\cap C_{\text{loc}}^{0}(\mathcal{Y})$
to holomorphic functions $v_{\infty}$ on $\mathcal{Y}$ satisfying
the following:

$(i)'$ $\int_{\mathcal{R}_{-2}^{\mathcal{Y}}}v_{\infty}^{\alpha}\overline{v}_{\infty}^{\beta}d\mu_{-2}^{\mathcal{Y}}=\delta_{\alpha\beta}$
for $1\leq\alpha,\beta\leq n$,

$(ii)'$ $\int_{\mathcal{R}_{-1}^{\mathcal{Y}}}v_{\infty}^{\alpha}\overline{v}_{\infty}^{\beta}d\mu_{-1}^{\mathcal{Y}}=0$
for $1\leq\alpha\neq\beta\leq n$,

$(iii)'$ $\int_{\mathcal{R}_{-2}^{\mathcal{Y}}}|v_{\infty}^{\alpha}|^{2}d\mu_{-2}^{\mathcal{Y}}\leq2\int_{\mathcal{R}_{-1}^{\mathcal{Y}}}|v_{\infty}^{\alpha}|^{2}d\mu_{-1}^{\mathcal{Y}}$
for $1\leq\alpha\leq n$,

$(iv)'$ $v_{\infty}^{\alpha}(x_{0})=0$ for $1\leq\alpha\leq n$,

$(v)'$ $\int_{\mathcal{R}_{-N}^{\mathcal{Y}}}|v_{\infty}^{\alpha}|^{2}d\mu_{-N}^{\mathcal{Y}}\leq2^{13}$.

\noindent Moreover, we can write $\mu_{t}^{\mathcal{Y}}=\nu_{y_{\ast};t}$,
where $y_{\ast}\in\mathcal{Y}_{0}$ corresponds to a vertex of the
metric cone. By Lemma \ref{conereplacement}, this implies that $\mathcal{Y}_{0}$
is isometric to flat $\mathbb{C}^{n}$, so that $x_{0}\in\mathcal{R}_{0}$,
a contradiction. $\square$

Now fix $\tau_{0}>0$ such that Claim 3 holds and $\Omega_{\tau}\leq2^{\frac{\epsilon}{100}}\Omega_{N\tau}$
for all $\tau\in(0,\tau_{0}]$, where $\epsilon>0$ is as in Claim
3. By rescaling, we may assume that $\tau_{0}=1$. Let $(v^{\alpha})_{\alpha=1}^{n}$
be a unitary basis of $(\mathcal{V},\Omega_{1})$ which is also $\Omega_{N^{-1}}$-orthogonal.
Then Claim 3 gives
\[
\sum_{\alpha=1}^{n}\int_{B(x_{0},2r_{0})\cap\mathcal{R}_{-N^{-1}}}|v^{\alpha}|^{2}d\nu_{x_{0};-N^{-1}}\leq\frac{n}{N^{1+\epsilon}},
\]
so that some $\alpha_{0}\in\{1,...,n\}$ satisfies 
\[
\int_{B(x_{0},2r_{0})\cap\mathcal{R}_{-N^{-1}}}|v^{\alpha_{0}}|^{2}d\nu_{x_{0};-N^{-1}}\leq\frac{1}{N^{1+\epsilon}}.
\]
Define $\widehat{v}^{\alpha}:=N^{-2\epsilon}v^{\alpha}$ for $\alpha\neq\alpha_{0}$,
and $\widehat{v}^{\alpha_{0}}:=N^{2n\epsilon}v^{\alpha_{0}}$, so
that for any $\alpha\in\{1,...,n\}$, we have 
\[
\int_{B(x_{0},2r_{0})\cap\mathcal{R}_{-N^{-1}}}|v^{\alpha}|^{2}d\nu_{x_{0};-N^{-1}}\leq\frac{1}{N^{1+\epsilon}}.
\]
Thus, $\mathcal{T}_{1}:(\mathcal{V},\Omega_{1})\to(\mathcal{V},\Omega_{N^{-1}})$,
$v^{\alpha}\mapsto\widehat{v}^{\alpha}$ is a $\mathbb{C}$-linear
endomorphism of $\mathcal{V}$ such that $|\mathcal{T}_{1}|_{\text{op}}\leq N^{-\epsilon}$
and 
\[
\mathcal{T}_{1}v^{1}\wedge\cdots\wedge\mathcal{T}_{1}v^{n}=v^{1}\wedge\cdots\wedge v^{n}.
\]
By iterating this construction (and also composing with unitary transformations),
we obtain a sequence of endomorphisms $\mathcal{T}_{j}:(\mathcal{V},\Omega_{N^{-j+1}})\to(\mathcal{V},\Omega_{N^{-j}})$
satisfying $|\mathcal{T}_{j}|_{\text{op}}\leq N^{-\epsilon}$ and
such that $v_{j}^{\alpha}:=\mathcal{T}_{j}\cdots\mathcal{T}_{1}v$
satisfy
\[
v_{j}^{1}\wedge\cdots\wedge v_{j}^{n}=v^{1}\wedge\cdots\wedge v^{n}.
\]
Let $\mathcal{T}_{0}:(\mathcal{V},\Omega_{1})\to(\mathcal{V},\Omega_{1})$
be given by $\mathcal{T}_{0}u^{\alpha}=v^{\alpha}$, so that 
\[
v_{j}^{1}\wedge\cdots\wedge v_{j}^{n}=(\det\mathcal{T}_{0})u^{1}\wedge\cdots\wedge u^{n}.
\]
By Lemma \ref{needforconvergence}$(ii)$, we have
\[
\sup_{B(x_{0},N^{-j})}|v_{j}^{\alpha}|^{2}\leq C(Y)\int_{B(x_{0},2r_{0})}|v_{j}^{\alpha}|^{2}d\nu_{-N^{-j}}\leq C(Y)N^{-(1+\epsilon)j}.
\]
Note that Corollary 6.3 of \cite{bam1} applies for infinitely many
$j\in\mathbb{N}$ by the proof of Lemma \ref{bddimplieslip}, so that
\[
\sup_{B(x_{0},N^{-j})}|\nabla v_{j}^{\alpha}|^{2}\leq C(Y)N^{-\epsilon j}
\]
for infinitely many $j\in\mathbb{N}$. Thus, the Lelong number $\Theta(x_{0})$
of $du^{1}\wedge\cdots\wedge du^{n}$ at $x_{0}$ satisfies (also
using Claim 2)
\[
\Theta(x_{0})=\lim_{j\to\infty}\inf_{x\in B(x_{0},N^{-j})\cap\mathcal{R}_{-N^{-j}}}\frac{\log|du^{1}\wedge\cdots\wedge du^{n}|(x)}{\log|u(x)|}\geq\liminf_{j\to\infty}\frac{\log(CN^{-\epsilon j})}{\log(C^{-1}N^{-\beta j})}=\frac{\epsilon}{\beta}>0.
\]
On the other hand, we have
\begin{align*}
\sqrt{-1}\partial\overline{\partial}\log|du^{1}\wedge\cdots\wedge du^{n}|= & \sqrt{-1}\partial\overline{\partial}\log\left(\frac{du^{1}\wedge\cdots\wedge du^{n}}{\omega_{0}^{n}}\right)\\
= & Rc(\omega_{0})=0
\end{align*}
on the subset of $\mathcal{R}_{0}$ where $|du^{1}\wedge\cdots\wedge du^{n}|\neq0$.
In particular, $\log|du^{1}\wedge\cdots\wedge du^{n}|$ is pluriharmonic
on $\mathcal{R}_{0}\cap U$. Because $U\setminus\mathcal{R}_{0}$
has singularities of codimension four with respect to the standard
flat metric on the chart $U$ (since the $u^{\alpha}$ are Lipschitz)
we know that $\log|du^{1}\wedge\cdots\wedge du^{n}|$ is actually
pluriharmonic on all of $U$ (where the complex structure is defined
using the holomorphic chart $(U,(u^{\alpha})_{\alpha=1}^{n})$). In
particular, it is bounded, so $\Theta(x_{0})=0$ a contradiction. 
\end{proof}
Next, we use ideas from \cite{donaldsun2} to show that $\mathcal{X}$
is in fact affine algebraic. 

We recall from Section 3.8 of \cite{bamlergen2} that any metric soliton
$(\mathcal{X},(\mu_{t})_{t\in(-\infty,0})$ is equipped with maps
$\psi_{\lambda}:\mathcal{X}\to\mathcal{X}$ for each $\lambda>0$
such which are metric flow isometries when viewed as mapping from
the parabolic rescaling $\mathcal{X}^{\lambda^{2}}$ of $\mathcal{X}$.
Moreover, the restriction of these maps to $\mathcal{R}$ form the
one-parameter group generated by the flow of $\tau(\partial_{\mathfrak{t}}-\nabla f)$. 
\begin{defn}
The group of static automorphisms $\mathcal{G}$ of $\mathcal{X}$
is the group of metric flow pair isometries $\phi:(\mathcal{X},(\mu_{t})_{t\in(-\infty,0]})\to(\mathcal{X},(\mu_{t})_{t\in(-\infty,0]})$
which commute with the maps $\psi_{\lambda}$. 
\end{defn}

\begin{rem}
\label{groupremark} Note that any $\phi\in\mathcal{G}$ is determined
by its restriction to $\mathcal{X}_{-1}$: in fact, we have 
\[
\phi_{-\lambda}=\psi_{\lambda}\circ\phi_{-1}\circ\psi_{\lambda^{-1}}.
\]
Moreover, in terms of the identification of $\mathcal{X}$ with its
soliton model $(X,d,\mu,(\nu_{x;t}')_{x\in X,t\in(-\infty,0]})$ (see
Definition 3.57 of \cite{bamlergen2}) we can identify $\mathcal{G}$
with the group of isometries $\phi$ of the metric measure space $(X,d,\mu)$
which satisfy
\[
\phi_{\ast}\nu_{x;t}'=\nu_{\phi(x);t}'
\]
for all $x\in X$ and $t\in(-\infty,0]$. 
\end{rem}

\begin{prop}
$\mathcal{G}$ equipped with the compact-open topology admits a unique
smooth structure making it a Lie group. 
\end{prop}

\begin{proof}
By Remark \ref{groupremark}, we can view $\mathcal{G}$ as a subgroup
of the group $\text{Iso}(X,d)$ of isometries of $(X,d)$ by the restriction
$\phi\mapsto\phi|_{\mathcal{X}_{-1}}$. Moreover, because $\mathcal{G}$
preserves the regular set $\mathcal{R}\subseteq\mathcal{X}$, we know
that $\phi|_{\mathcal{X}_{-1}}$ preserves $\phi|_{\mathcal{R}_{-1}}$,
so we can consider the further restriction $\widetilde{\phi}:=\phi|_{\mathcal{R}_{-1}}:(\mathcal{R}_{-1},d_{g_{-1}})\to(\mathcal{R}_{-1},d_{g_{-1}})$,
which is a metric isometry of Riemannian manifolds, hence a Riemannian
isometry. The map $\mathcal{G}\to\text{Iso}(\mathcal{R}_{-1},g_{-1})$,
$\phi\mapsto\widetilde{\phi}$ thus realizes $\mathcal{G}$ as a subgroup
of the Lie group of isometries of a Riemannian manifold (where $\text{Iso}(\mathcal{R}_{-1},g_{-1})$
is equipped with the compact-open topology), so it suffices to show
that $\mathcal{G}$ is a closed subgroup.

Suppose $(\phi_{i})$ is a sequence in $\mathcal{G}$ such that $\widetilde{\phi}_{i}\to\widetilde{\phi}$
in the compact-open topology, where $\widetilde{\phi}\in\text{Iso}(\mathcal{R}_{-1},g_{-1})$.
Then $\widetilde{\phi}$ extends to an isometry $\phi_{-1}:(\mathcal{X}_{-1},d_{-1})\to(\mathcal{X}_{-1},d_{-1})$,
and $\phi_{i}$ converge locally uniformly to $\phi$. By defining
$\phi_{-\lambda}:=\psi_{\lambda}\circ\phi_{-1}\circ\psi_{\lambda}^{-1}$,
we can view $\phi$ as a time-preserving map $\phi:\mathcal{X}\to\mathcal{X}$
which commutes with the $(\psi_{\lambda})_{\lambda>0}$ and restricts
to a metric isometry of each time slice. 

\noindent \textbf{Claim 1: }$(\phi_{t})_{\ast}\mu_{t}=\mu_{t}$ for
all $t\in(-\infty,0]$. 

For any $h\in C_{b}(\mathcal{X}_{t})$, $(\phi_{i,t})_{\ast}\mu_{t}=\mu_{t}$
implies 
\[
\int_{\mathcal{X}_{t}}hd\mu_{t}=\int_{\mathcal{X}_{t}}hd((\phi_{i,t})_{\ast}\mu_{t})=\int_{\mathcal{X}_{t}}(h\circ\phi_{i,t})d\mu_{t}.
\]
Because $h\circ\phi_{i,t}\to h\circ\phi_{t}$ pointwise and $|h\circ\phi_{i,t}|\leq||h||_{C^{0}(\mathcal{X}_{t})}$,
the dominated convergence theorem gives
\[
\int_{\mathcal{X}_{t}}hd\mu_{t}=\lim_{i\to\infty}\int_{\mathcal{X}_{t}}(h\circ\phi_{i,t})d\mu_{t}=\int_{\mathcal{X}_{t}}(h\circ\phi_{t})d\mu_{t}=\int_{\mathcal{X}_{t}}hd((\phi_{t})_{\ast}\mu_{t}).\hfill\text{\ensuremath{\square}}
\]

\noindent \textbf{Claim 2: }$(\phi_{s})_{\ast}\nu_{x;s}=\nu_{\phi_{t}(x);s}$
for all $x\in\mathcal{X}_{t}$ and $s\in(-\infty,t]$. 

For $L$-Lipschitz $h\in C_{b}(\mathcal{X}_{s})$, we similarly have
\[
\int_{\mathcal{X}_{s}}hd\nu_{\phi_{i,t}(x);s}=\int_{\mathcal{X}_{s}}(h\circ\phi_{i,s})d\nu_{x;s},
\]
and by Kantorovich duality, we know
\[
\left|\int_{\mathcal{X}_{s}}hd\nu_{\phi_{i,t}(x);s}-\int_{\mathcal{X}_{s}}hd\nu_{\phi_{t}(x);s}\right|\leq Ld_{W_{1}}^{\mathcal{X}_{s}}(\nu_{\phi_{i,t}(x);s},\nu_{\phi_{t}(x);s})\leq Ld_{t}(\phi_{i,t}(x),\phi_{t}(x)),
\]
so that taking $i\to\infty$ gives 
\[
\int_{\mathcal{X}_{s}}hd\nu_{\phi_{t}(x);s}=\int_{\mathcal{X}_{s}}(h\circ\phi_{s})d\nu_{x;s}=\int_{\mathcal{X}_{s}}hd((\phi_{s})_{\ast}\nu_{x;s})
\]
for any bounded Lipshitz function $h$ on $\mathcal{X}_{s}$. It follows
that $(\phi_{s})_{\ast}\nu_{x;s}=\nu_{\phi_{t}(x);s}$. 
\end{proof}
\begin{proof}[Completion of the proof of Theorem \ref{thm1}]
 By Proposition 5.1 of \cite{complexsplitting} and its proof, the
gradient flow of $\tau J\nabla f$ on $\mathcal{R}$ induces a 1-parameter
family $(\phi^{s})_{s\in\mathbb{R}}$ of isometries of the metric
flow pair $(\mathcal{X},(\mu_{t})_{t\in(-\infty,0]})$. Moreover,
because $\tau J\nabla f$ commutes with the generator $\tau(\partial_{\mathfrak{t}}-\nabla f)$
of $(\psi_{\lambda})_{\lambda>0}$, we know that $\phi^{s}$ commute
with $\psi_{\lambda}$, hence $(\phi^{s})_{s\in\mathbb{R}}$ is a
one-parameter subgroup of $\mathcal{G}$. We may therefore take its
closure to obtain a (real) torus $\mathbb{T}\subseteq\mathcal{G}$
of rank $r$. For any holomorphic function $(u_{t})_{t\in(-\infty,0)}$
on $\mathcal{X}$, and any $\alpha\in\mathbb{Z}^{r}$, we may therefore
define holomorphic functions
\[
(u_{t})_{\alpha}(z):=\int_{\mathbb{T}}u_{t}(e^{\sqrt{-1}\theta}\cdot z)e^{-\sqrt{-1}\langle\theta,\alpha\rangle}d\theta,
\]
where we identify $\mathbb{T}\cong[0,2\pi]^{r}$, and where $d\theta$
is the Haar measure of $\mathbb{T}$. Arguing as in Section 2.3 of
\cite{donaldsun2}, we get
\[
u_{t}(z)=\sum_{\alpha\in\mathbb{Z}}(u_{t})_{\alpha}(z),
\]
where the sum converges locally uniformly on $\mathcal{X}$, and if
$(u^{1},...,u^{N})$ give a local embedding of the vertex of $\mathcal{X}$,
then (by the inverse function theorem for holomorphic varieties, for
example) there are $m_{1},...,m_{N}\in\mathbb{N}^{r}$ such that each
\[
\widehat{F}_{t}:=\left(\sum_{|\alpha|\leq m_{1}}(u_{t}^{1})_{\alpha},...,\sum_{|\alpha|\leq m_{N}}(u_{t}^{N})_{\alpha}\right):\mathcal{X}_{t}\to\mathbb{C}^{N}
\]
is also an embedding. In particular, we know that the map $F_{t}:\mathcal{X}_{t}\to\mathbb{C}^{Q}$
obtained by taking each of the components $(u_{t}^{j})_{\alpha}$
with $1\leq j\leq N$ and $|\alpha|\leq m_{j}$ is a $\mathbb{T}$-equivariant
emebedding of $\mathcal{X}_{t}$ with respect to a linear action of
$\mathbb{T}$ on $\mathbb{C}^{Q}$. We may then argue as in Lemma
2.19 of \cite{donaldsun2} to conclude that $F(\mathcal{X}_{t})$
is an affine algebraic variety. 
\end{proof}

\section{Appendix: The Heat Kernel of a Ricci-Flat Cone}

In this appendix, we consider the heat kernel $K$ corresponding to
a static cone $(\mathcal{X},(\nu_{t})_{t\in(-\infty,0)})$ obtained
as a limit of (not necessarily Kähler) Ricci flows as in (\ref{F}),
assuming in particular the uniform noncollapsing bound $\mathcal{N}_{x_{i},0}(\epsilon_{i}^{-1})\geq-Y$.
We will relate $K$ to the static heat kernel on the Riemannian manifold
$(\mathcal{R}_{-1},g_{-1})$, and use this to prove a bound for $|\nabla_{y}K(x;y)|$
which is needed in Section 4. 

For $x,y\in\mathcal{R}_{-1}$ and $t>0$, we define 
\[
K(x,y,t):=K(x;y(-1-t)).
\]
Let $x_{0}\in\mathcal{X}_{-1}$ correspond to the vertex of the cone,
so that $f_{t}(x)=\frac{1}{2\tau}d_{-1}^{2}(\cdot,x_{0})$. 
\begin{prop}
\label{heatkernel} For any $x,y\in\mathcal{R}_{-1}$ and $t>0$,
we have the following:

\[
\frac{1}{(4\pi t)^{\frac{n}{2}}}\exp\left(-\frac{d_{-1}^{2}(x,y)}{4t}\right)\leq K(x,y,t)\leq\frac{C(Y)}{t^{\frac{n}{2}}}\exp\left(-\frac{d_{-1}^{2}(x,y)}{10t}\right),
\]
\[
|\nabla_{x}\log K_{f}(x,y,t)|\leq\frac{C(Y)}{\sqrt{t}}\left(1+\frac{d_{-1}(x,y)}{\sqrt{t}}\right).
\]
\end{prop}

\begin{proof}
By noting that the $\mathcal{L}$-length of a the path $x\mapsto x(t)$
in $\mathcal{R}$ is zero, Lemma 22.2 of \cite{bamlergen3} gives
\[
K(x;x(-1-t))\geq\frac{1}{(4\pi t)^{\frac{n}{2}}}
\]
and also 
\[
d_{W_{1}}^{\mathcal{X}_{-1-t}}(\nu_{x;-1-t},\delta_{x(-1-t)})\leq C(Y)\sqrt{t}
\]
for all $t>0$. Thus Lemma 15.9 of \cite{bamlergen3} gives
\[
K(x;y(-1-t))\leq\frac{C(Y)}{t^{\frac{n}{2}}}\exp\left(-\frac{d_{-1}^{2}(x,y)}{10t}\right)
\]
for $x,y\in\mathcal{R}_{-1}$ and $t>0$. Moreover, the integrated
form of Perelman's differential Harnack inequality gives
\[
K(x;y(-1-t))\geq\frac{1}{(4\pi t)^{\frac{n}{2}}}\exp\left(-\frac{d_{-1}^{2}(x,y)}{4t}\right).
\]
Thus we have 
\[
\frac{1}{(4\pi t)^{\frac{n}{2}}}\exp\left(-\frac{d_{-1}^{2}(x,y)}{4t}\right)\leq K(x,y,t)\leq\frac{C(Y)}{t^{\frac{n}{2}}}\exp\left(-\frac{d_{-1}^{2}(x,y)}{10t}\right)
\]
for all $x,y\in\mathcal{R}_{-1}$ and $t>0$. 

Next, passing the gradient estimate of \cite{qizhanggradient} to
the limit and combining with the above estimates gives the following
for all $x,y\in\mathcal{R}_{-1}$:
\begin{align*}
|\nabla_{x}\log K(x;y(-1-t))|\leq & \frac{1}{\sqrt{t}}\sqrt{\log\left(\frac{C(Y)}{\exp\left(-\frac{d_{-1}^{2}(x,y)}{4t}\right)}\right)}\\
\leq & \frac{C(Y)}{\sqrt{t}}\left(1+\frac{d_{-1}(x,y)}{\sqrt{t}}\right).
\end{align*}
\end{proof}
Next we show that $K$ is the heat kernel of the static manifold $(\mathcal{R}_{-1},g_{-1})$.
\begin{lem}
\label{fundsol} $(i)$ For all $(x,y,t)\in\mathcal{R}_{-1}\times\mathcal{R}_{-1}\times(0,\infty)$,
we have
\[
(\partial_{t}-\Delta_{x})K(x,y,t)=0,
\]
\[
(\partial_{t}-\Delta_{y})K(x,y,t)=0.
\]
$(ii)$ For any Lipschitz $\psi\in C_{c}(\mathcal{R}_{-1})$, we have 

\[
\lim_{t\searrow0}\int_{M}K(x,y,t)\psi(y)dg_{-1}(y)=\psi(x).
\]
\[
\lim_{t\searrow0}\int_{M}K(x,y,t)\psi(x)dg_{-1}(x)=\psi(y).
\]
\end{lem}

\begin{proof}
$(i)$ This is an easy computation using the infinitesimal symmetry
(Theorem 15.60 of \cite{bamlergen3})
\[
\partial_{\mathfrak{t}}|_{x}K(x;y)+\partial_{\mathfrak{t}}|_{y}K(x;y)=0.
\]
$(ii)$ We observe that
\begin{align*}
\int_{\mathcal{R}_{-1}}K(x,y,t)\psi(y)dg_{-1}(y)= & \int_{\mathcal{R}_{-1-t}}K(x;y)\psi(y)dg_{-1-t}(y)\\
= & \int_{\mathcal{R}_{-1-t}}\psi d\nu_{x;-1-t}.
\end{align*}
By the assumption that $\psi:(\mathcal{X}_{-1},d_{-1})\to\mathbb{R}$
is $L$-Lipschitz for some $L<\infty$, Kantorovich duality yields
\begin{align*}
\left|\int_{\mathcal{R}_{-1-t}}\psi(y)d\nu_{x;-1-t}-\psi(x)\right|\leq & Ld_{W_{1}}^{\mathcal{X}_{-1-t}}(\nu_{x;-1-t},\delta_{x(-1-t)})\leq LC(Y)\sqrt{t},
\end{align*}
so combining expressions 
\[
\lim_{t\searrow0}\int_{\mathcal{R}_{-1}}K(x,y,t)\psi(y)dg_{-1}(y)=\psi(x).
\]
Next, we use Lemma \ref{heatkernel} to estimate 

\begin{align*}
\left|\int_{\mathcal{R}_{-1}}K(x,y,t)\psi(x)dg_{-1}(x)-\psi(y)\right|\leq & \int_{\mathcal{R}_{-1}}K(x;y(-1-t))|\psi(x)-\psi(y)|dg_{-1}(x)\\
\leq & \frac{C(Y)L}{t^{\frac{n}{2}}}\int_{\mathcal{R}_{-1}}d_{-1}(x,y)\exp\left(-\frac{d_{-1}^{2}(x,y)}{10t}\right)dg_{-1}(x)\\
\leq & \frac{C(Y)L}{t^{\frac{n}{2}}}\sum_{j\in\mathbb{Z}}\int_{B(y,2^{j+1}\sqrt{t})\setminus B(y,2^{j}\sqrt{t})}d_{-1}(x,y)\exp\left(-\frac{d_{-1}^{2}(x,y)}{10t}\right)dg_{-1}(x)\\
\leq & \frac{C(Y)L\sqrt{t}}{t^{\frac{n}{2}}}\sum_{j\in\mathbb{Z}}2^{j+1}e^{-\frac{1}{10}2^{2j}}\text{Vol}(B(y,2^{j}\sqrt{t})).
\end{align*}
We recall that there exists $A=A(Y)>0$ such that for all $z\in\mathcal{R}_{-1}\cap B(x_{0},1)$
and $r\in(0,1]$, we have
\[
\text{Vol}(B(z,r))\leq Ar^{n}.
\]
Letting $\tau_{\lambda}:\mathcal{R}_{-1}\to\mathcal{R}_{-1}$ denote
the dilation map by $\lambda$ (so that $\tau_{\lambda}^{\ast}g=\lambda^{2}g$),
and taking $\lambda>2^{j}\sqrt{t}+d_{-1}(y,x_{0})$, we then have
\[
\text{Vol}(B(y,2^{j}\sqrt{t}))=\lambda^{n}\text{Vol}\left(B(\tau_{\lambda}^{-1}(y),\lambda^{-1}2^{j}\sqrt{t})\right)\leq A\lambda^{n}(\lambda^{-1}2^{j}\sqrt{t})^{n}=A2^{nj}t^{\frac{n}{2}}.
\]
Combining estimates leads to
\[
\left|\int_{\mathcal{R}_{-1}}K(x,y,t)\psi(x)dg_{-1}(x)-\psi(y)\right|\leq C(Y)L\sqrt{t}\sum_{j\in\mathbb{Z}}2^{(n+1)j}e^{-\frac{1}{10}2^{2j}}\leq C(Y)L\sqrt{t}.
\]
\end{proof}
\begin{lem}
\label{symmetry} We have $K(x,y,t)=K(y,x,t)$ for all $x,y\in\mathcal{R}_{-1}$.
\end{lem}

\begin{rem}
Note that this implies $K$ extends to a continuous function on all
of $\mathcal{X}_{-1}\times\mathcal{X}_{-1}\times(0,\infty)$, which
is also symmetric.
\end{rem}

\begin{proof}
Let $\chi_{\rho},\eta_{\epsilon}$ be as in Lemma \ref{LocLipschitz}.
For fixed $T>0$, we can compute
\begin{align*}
\frac{d}{dt}\int_{\mathcal{R}_{-1}}K(z,x,T-t) & \eta_{\epsilon}(z)\chi_{\rho}(z)K(z,y,t)dg_{-1}(z)\\
= & \int_{\mathcal{R}_{-1}}\left(-\Delta_{z}K(z,x,T-t)\eta_{\epsilon}(z)\chi_{\rho}(z)K(z,y,t)\right)dg_{-1}(z)\\
 & +\int_{\mathcal{R}_{-1}}\left(K(z,x,T-t)\eta_{\epsilon}(z)\chi_{\rho}(z)\Delta_{z}K(z,y,t)\right)dg_{-1}(z)\\
= & \int_{\mathcal{R}_{-1}}\langle\nabla_{z}K(z,x,T-t),\chi_{\rho}(z)\nabla\eta_{\epsilon}(z)+\eta_{\epsilon}(z)\nabla\chi_{\rho}(z)\rangle K(z,y,t)dg_{-1}(z)\\
 & +\int_{\mathcal{R}_{-1}}K(z,x,T-t)\langle\chi_{\rho}(z)\nabla\eta_{\epsilon}(z)+\eta_{\epsilon}(z)\nabla\chi_{\rho}(z),\nabla_{z}K(z,y,t)\rangle dg_{-1}(z)
\end{align*}
Recalling that $x,y\in\mathcal{R}_{-1}$, we can find $r_{0}>0$ such
that for all $\epsilon>0$ sufficiently small, we have 
\[
\inf_{z\in B\text{supp}(\eta_{\epsilon})}\min\{d_{-1}(y,z),d_{-1}(x,z)\}\geq r_{0}.
\]
We use Lemma \ref{heatkernel} to estimate
\begin{align*}
\biggl|\int_{\mathcal{R}_{-1}}K(z,x,T-t) & \langle\chi_{\rho}(z)\nabla\eta_{\epsilon}(z),\nabla_{z}K(z,y,t)\rangle dg_{-1}(z)\biggr|\\
\leq & \int_{B(x_{0},2\rho)}K(z,x,T-t)|\nabla\eta_{\epsilon}(z)|\cdot|\nabla_{z}K(z,y,t)|dg_{-1}(z)\\
\leq & \int_{B(x_{0},2\rho)}K(z,x,T-t)|\nabla\eta_{\epsilon}(z)|\left(\frac{C(Y)}{t^{\frac{n+1}{2}}}\left(1+\frac{d_{-1}(z,y)}{\sqrt{t}}\right)\exp\left(-\frac{d_{-1}^{2}(z,y)}{10t}\right)\right)dg_{-1}(z)\\
\leq & C(Y,T,\rho)\sup_{z\in B(x_{0},2r)\cap\text{supp}(\eta_{\epsilon})}\frac{1}{t^{\frac{n+1}{2}}}\left(1+\frac{d_{-1}(z,y)}{\sqrt{t}}\right)\exp\left(-\frac{d_{-1}^{2}(z,y)}{20t}\right)\\
 & \times\int_{B(x_{0},2\rho)}|\nabla\eta_{\epsilon}(z)|dg_{-1}(z)\\
\leq & C(Y,T,\rho,r_{0})\Psi(\epsilon|\rho)\sup_{t\in(0,\frac{T}{2}]}\frac{1}{t^{\frac{n+2}{2}}}\exp\left(-\frac{r_{0}^{2}}{40t}\right)\\
\leq & \Psi(\epsilon|\rho,r_{0},Y,T)
\end{align*}
whenever $t\in(0,\frac{T}{2}]$. If instead $t\in[\frac{T}{2},T)$,
we estimate
\begin{align*}
\left|\int_{\mathcal{R}_{-1}}K(z,x,T-t)\right. & \left.\langle\chi_{\rho}(z)\nabla\eta_{\epsilon}(z),\nabla_{z}K(z,y,t)\rangle dg_{-1}(z)\right|\\
\leq & C(T,\rho)\int_{B(x_{0},2\rho)\cap\mathcal{R}_{-1}}|\nabla\eta_{\epsilon}(z)|K(z,x,T-t)dg_{-1}(z)\\
\leq & C(T,\rho)\sup_{z\in B(x_{0},2\rho)\cap\text{supp}(\eta_{\epsilon})}\frac{C(Y,D)}{(T-t)^{\frac{n}{2}}}\exp\left(-\frac{d_{-1}^{2}(z,x)}{10(T-t)}\right)\int_{B(x_{0},2\rho)}|\nabla\eta_{\epsilon}(z)|dg_{-1}(z)\\
\leq & \Psi(\epsilon|\rho,r_{0},Y,T).
\end{align*}
Next, for any $t\in(0,T)$, we have 
\begin{align*}
\biggl|\int_{\mathcal{R}_{-1}}K(z,x,T-t) & \langle\nabla\chi_{\rho}(z),\nabla_{z}K(z,y,t)\rangle\eta_{\epsilon}(z)d\nu_{-1}(z)\biggr|\\
\leq & C\rho^{-1}\int_{B(x_{0},2\rho)\setminus B(x_{0},\rho)}K(z,x,T-t)|\nabla_{z}K(z,y,t)|d\nu_{-1}(z)\\
\leq & C(Y,T)\rho^{-1}\frac{1}{t^{\frac{n+1}{2}}}\left(1+\frac{\rho}{t^{\frac{1}{2}}}\right)\exp\left(-\frac{\rho^{2}}{40t}\right)\\
\leq & \Psi(\rho^{-1}|Y,T)
\end{align*}
if $\rho\geq d(x,x_{0})+d(y,x_{0})$. Combining expressions and integrating
yields
\[
\left|\left.\int_{\mathcal{R}_{-1}}K(z,x,T-t)\eta_{\epsilon}(z)\rho_{r}(z)K(z,y,t)dg_{-1}(z)\right|_{t=t_{1}}^{t=t_{2}}\right|\leq\Psi(\epsilon|\rho,Y,T,r_{0})+\Psi(\rho^{-1}|Y,T)
\]
for all $t_{1},t_{2}\in(0,T)$. Taking $t_{1}\searrow0$, $t_{2}\nearrow T$,
if $\epsilon\leq\overline{\epsilon}(x,y)$ and $\rho\geq\underline{\rho}(x,y)$,
we can apply Lemma \ref{fundsol} to obtain
\[
|K(y,x,T)-K(x,y,T)|\leq\Psi(\epsilon|\rho,Y,T,r_{0})+\Psi(\rho^{-1}|Y,T).
\]
Finally, take $\epsilon\searrow0$ and then $\rho\nearrow\infty$
to conclude $K(x,y,T)=K(y,x,T)$. 
\end{proof}
We can now combine the symmetry of $K$ with the gradient estimate
for $x\mapsto K(x;y)$ to establish a gradient estimate for $y\mapsto K(x;y)$. 
\begin{prop}
\label{conjheatgradient} For $t_{0}\in(-\infty,0)$, $\tau>0$, $x\in\mathcal{R}_{t_{0}}$,
and $y\in\mathcal{R}_{t_{0}-\tau}$, we have
\[
|\nabla_{y}K(x;y)|=|\nabla_{x}K(y;x(t_{0}-2\tau))|.
\]
\end{prop}

\begin{proof}
We may assume $t_{0}=-1$. Writing $K(x,y,t)=K(y,x,t)$ in terms of
$K(\cdot;\cdot)$ gives
\[
K(x;z(-1-\tau))=K(z;x(-1-\tau))
\]
for $x,z\in\mathcal{R}_{-1}$ and $\tau>0$. Taking $z=y(-1)$ for
$y\in\mathcal{R}_{-1-\tau}$ then gives
\[
K(x;y)=K(y(-1);x(-1-\tau))=K(y;x(-1-2\tau))
\]
Differentiating both sides in $y$ gives
\[
|\nabla_{y}K(x;y)|=|\nabla_{y}K(y;x(-1-2\tau))|.
\]
\end{proof}
\bibliographystyle{alpha}
\bibliography{ricciflow}

\begin{thebibliography}{CDS15b}

\bibitem[Bam17]{bam1}
Richard~H. Bamler.
\newblock Structure theory of singular spaces.
\newblock {\em J. Funct. Anal.}, 272(6):2504--2627, 2017.

\bibitem[Bam18]{bam2}
Richard Bamler.
\newblock Convergence of {R}icci flows with bounded scalar curvature.
\newblock {\em Ann. of Math. (2)}, 188(3):753--831, 2018.

\bibitem[Bam20a]{bamlergen2}
Richard~H Bamler.
\newblock Compactness theory of the space of super ricci flows, 2020.

\bibitem[Bam20b]{bamlergen1}
Richard~H Bamler.
\newblock Entropy and heat kernel bounds on a ricci flow background, 2020.

\bibitem[Bam20c]{bamlergen3}
Richard~H Bamler.
\newblock Structure theory of non-collapsed limits of ricci flows, 2020.

\bibitem[BK23]{julius}
Julius Baldauf and Dain Kim.
\newblock Parabolic frequency on {R}icci flows.
\newblock {\em Int. Math. Res. Not. IMRN}, (12):10098--10114, 2023.

\bibitem[CC97]{cheegercolding1}
Jeff Cheeger and Tobias~H. Colding.
\newblock On the structure of spaces with {R}icci curvature bounded below. {I}.
\newblock {\em J. Differential Geom.}, 46(3):406--480, 1997.

\bibitem[CCT02]{cheegercoldingtian}
J.~Cheeger, T.~H. Colding, and G.~Tian.
\newblock On the singularities of spaces with bounded {R}icci curvature.
\newblock {\em Geom. Funct. Anal.}, 12(5):873--914, 2002.

\bibitem[CDS15a]{CDSI}
Xiuxiong Chen, Simon Donaldson, and Song Sun.
\newblock K\"{a}hler-{E}instein metrics on {F}ano manifolds. {I}:
  {A}pproximation of metrics with cone singularities.
\newblock {\em J. Amer. Math. Soc.}, 28(1):183--197, 2015.

\bibitem[CDS15b]{CDSII}
Xiuxiong Chen, Simon Donaldson, and Song Sun.
\newblock K\"{a}hler-{E}instein metrics on {F}ano manifolds. {II}: {L}imits
  with cone angle less than {$2\pi$}.
\newblock {\em J. Amer. Math. Soc.}, 28(1):199--234, 2015.

\bibitem[CDS15c]{CDSIII}
Xiuxiong Chen, Simon Donaldson, and Song Sun.
\newblock K\"{a}hler-{E}instein metrics on {F}ano manifolds. {III}: {L}imits as
  cone angle approaches {$2\pi$} and completion of the main proof.
\newblock {\em J. Amer. Math. Soc.}, 28(1):235--278, 2015.

\bibitem[CN13]{cheegernaberquant}
Jeff Cheeger and Aaron Naber.
\newblock Lower bounds on {R}icci curvature and quantitative behavior of
  singular sets.
\newblock {\em Invent. Math.}, 191(2):321--339, 2013.

\bibitem[CN15]{cheegernaberdim4}
Jeff Cheeger and Aaron Naber.
\newblock Regularity of {E}instein manifolds and the codimension 4 conjecture.
\newblock {\em Ann. of Math. (2)}, 182(3):1093--1165, 2015.

\bibitem[CSW18]{chensunwang}
Xiuxiong Chen, Song Sun, and Bing Wang.
\newblock K\"{a}hler-{R}icci flow, {K}\"{a}hler-{E}instein metric, and
  {K}-stability.
\newblock {\em Geom. Topol.}, 22(6):3145--3173, 2018.

\bibitem[CW17]{chenwang2a}
Xiuxiong Chen and Bing Wang.
\newblock Space of {R}icci flows ({II})---{P}art {A}: {M}oduli of singular
  {C}alabi-{Y}au spaces.
\newblock {\em Forum Math. Sigma}, 5:e32, 103, 2017.

\bibitem[CW20]{chenwang2b}
Xiuxiong Chen and Bing Wang.
\newblock Space of {R}icci flows ({II})---{P}art {B}: {W}eak compactness of the
  flows.
\newblock {\em J. Differential Geom.}, 116(1):1--123, 2020.

\bibitem[DS14]{donaldsun1}
Simon Donaldson and Song Sun.
\newblock Gromov-{H}ausdorff limits of {K}\"{a}hler manifolds and algebraic
  geometry.
\newblock {\em Acta Math.}, 213(1):63--106, 2014.

\bibitem[DS17]{donaldsun2}
Simon Donaldson and Song Sun.
\newblock Gromov-{H}ausdorff limits of {K}\"{a}hler manifolds and algebraic
  geometry, {II}.
\newblock {\em J. Differential Geom.}, 107(2):327--371, 2017.

\bibitem[DS20]{szekdegen}
Ruadha\'{\i} Dervan and G\'{a}bor Sz\'{e}kelyhidi.
\newblock The {K}\"{a}hler-{R}icci flow and optimal degenerations.
\newblock {\em J. Differential Geom.}, 116(1):187--203, 2020.

\bibitem[EGZ09]{singKE}
Philippe Eyssidieux, Vincent Guedj, and Ahmed Zeriahi.
\newblock Singular {K}\"{a}hler-{E}instein metrics.
\newblock {\em J. Amer. Math. Soc.}, 22(3):607--639, 2009.

\bibitem[HJ22]{complexsplitting}
Max Hallgren and Wangjian Jian.
\newblock Tangent flows of k\"ahler metric flows, 2022.

\bibitem[HJST23]{JSTII}
Max Hallgren, Wangjian Jian, Jian Song, and Gang Tian.
\newblock Geometric regularity of blow-up limits of the k\"ahler-ricci flow,
  2023.

\bibitem[HL20]{lihanunique}
Jiyuan Han and Chi Li.
\newblock Algebraic uniqueness of k\"{a}hler-ricci flow limits and optimal
  degenerations of fano varieties, 2020.

\bibitem[HN14]{hein}
Hans-Joachim Hein and Aaron Naber.
\newblock New logarithmic {S}obolev inequalities and an {$\epsilon$}-regularity
  theorem for the {R}icci flow.
\newblock {\em Comm. Pure Appl. Math.}, 67(9):1543--1561, 2014.

\bibitem[JST23]{JSTI}
Wangjian Jian, Jian Song, and Gang Tian.
\newblock Finite time singularities of the k\"ahler-ricci flow, 2023.

\bibitem[Kin71]{king}
James~R. King.
\newblock The currents defined by analytic varieties.
\newblock {\em Acta Math.}, 127(3-4):185--220, 1971.

\bibitem[KL17]{klott1}
Bruce Kleiner and John Lott.
\newblock Singular {R}icci flows {I}.
\newblock {\em Acta Math.}, 219(1):65--134, 2017.

\bibitem[KS23]{kunikawaII}
Keita Kunikawa and Yohei Sakurai.
\newblock Almost splitting and quantitative stratification for super ricci
  flow, 2023.

\bibitem[Liu18]{gangliucone}
Gang Liu.
\newblock On the tangent cone of {K}\"{a}hler manifolds with {R}icci curvature
  lower bound.
\newblock {\em Math. Ann.}, 370(1-2):649--667, 2018.

\bibitem[LS20]{szek1}
Gang Liu and G\'{a}bor Sz\'{e}kelyhidi.
\newblock Gromov-hausdorff limits of k\"ahler manifolds with ricci curvature
  bounded below, 2020.

\bibitem[LS21]{szek2}
Gang Liu and G\'{a}bor Sz\'{e}kelyhidi.
\newblock Gromov-{H}ausdorff limits of {K}\"{a}hler manifolds with {R}icci
  curvature bounded below {II}.
\newblock {\em Comm. Pure Appl. Math.}, 74(5):909--931, 2021.

\bibitem[LTZ20]{tianGman}
Yan Li, Gang Tian, and Xiaohua Zhu.
\newblock Singular limits of k\"ahler-ricci flow on fano $g$-manifolds, 2020.

\bibitem[LWX21]{liuniquecones}
Chi Li, Xiaowei Wang, and Chenyang Xu.
\newblock Algebraicity of the metric tangent cones and equivariant
  {K}-stability.
\newblock {\em J. Amer. Math. Soc.}, 34(4):1175--1214, 2021.

\bibitem[Tia15]{TianKStable}
Gang Tian.
\newblock K-stability and {K}\"{a}hler-{E}instein metrics.
\newblock {\em Comm. Pure Appl. Math.}, 68(7):1085--1156, 2015.

\bibitem[TZ16]{tianzhen}
Gang Tian and Zhenlei Zhang.
\newblock Regularity of {K}\"{a}hler-{R}icci flows on {F}ano manifolds.
\newblock {\em Acta Math.}, 216(1):127--176, 2016.

\bibitem[Zha06]{qizhanggradient}
Qi~S. Zhang.
\newblock Some gradient estimates for the heat equation on domains and for an
  equation by {P}erelman.
\newblock {\em Int. Math. Res. Not.}, pages Art. ID 92314, 39, 2006.

\end{thebibliography}
 
\end{document}